\documentstyle[11pt,amssymb,mysection]{article}
\voffset = - 2.3 truecm
\def\theequation{\arabic{section}.\arabic{equation}}
\newcounter{theo}[section]
\newcounter{lemma}[section]
\newcounter{defi}[section]
\newcounter{rem}[section]
\newcounter{cor}[section]
\def\thetheo{\thesection.\arabic{theo}}
\def\thelemma{\thesection.\arabic{lemma}}
\def\thedefi{\thesection.\arabic{defi}}
\def\therem{\thesection.\arabic{rem}}
\def\thecor{\thesection.\arabic{cor}}
\newenvironment{theo}{%
\refstepcounter{theo}
        { \noindent  THEOREM \thetheo. }
            \bgroup\rm}{\egroup}
\newenvironment{lemma}{%
\refstepcounter{lemma}
         {\noindent  LEMMA \thelemma. } 
             \bgroup\rm}{\egroup}
\newenvironment{defi}{%
\refstepcounter{defi}
      {\noindent  DEFINITION \thedefi. }
          \bgroup\rm}{\egroup}

\newenvironment{rem}{%
\refstepcounter{rem}
      {\noindent  REMARK \therem. } \bgroup\rm}{\egroup}
\newenvironment{cor}{%
\refstepcounter{cor}
         {\vspace*{1.0em} \noindent  COROLLARY \thecor. }
             \bgroup\rm}{\egroup}
\newfont{\cmss}{cmss12 scaled 1000}
\newfont{\cmex}{cmex10 scaled 1440}
\newfont{\eufm}{eufm10 scaled 1200}
\newfont{\bib}{cmcsc10 scaled 1000}

\newcommand{\Mgot}{\mbox{\eufm\symbol{77}}}

\begin{document}

\noindent
\begin{flushright}
{\it Dedicated to Moshe Liv\v{s}ic, Morenu and Rabenu  \hspace{4em}\\}
\end{flushright}
\ \vskip 3cm
\noindent
\hspace{5ex}
\begin{minipage}{14.0cm}
{\bf RIGHT AND LEFT JOINT SYSTEM REPRESENTATION OF\\ A RATIONAL
MATRIX FUNCTION IN GENERAL POSITION \\[0.1cm]
{\small\bf (SYSTEM REPRESENTATION THEORY FOR DUMMIES)}\\[0.4cm]
V.E. Katsnelson  
}
\end{minipage}

\vskip .5cm
\begin{abstract}
For a rational $k\times k$ matrix function \(R\) of one variable in general
 position, the matrix functions $R(x)\cdot R^{-1}(y)$ and 
$R^{-1}(x)\cdot R(y)$ of two variables are considered.
 For these matrix functions of two
variables, the representations  which are analogous to the 
 system representation (or realization) of a  rational matrix function of one
variable are constructed.
 This representation of the function $R(x)\cdot R^{-1}(y)$ 
(of the function $R^{-1}(x)\cdot R(y)$) is said to be the joint right
(respectively the joint left) system representation of the matrix functions
$R,\, R^{-1}$.
In these representations there appear diagonal $n\times n$ matrices,
 $A_{\cal P}=\mbox{\rm diag}\, ( \lambda _1, \dots ,\lambda _n)$ 
(called the pole matrix for $R$) and 
$A_{\cal N}=\mbox{\rm diag}\, ( \mu_1, \dots ,\mu _n)$
 (called the zero matrix for $R$), where
  $\lambda _l, \dots , \lambda _n$ are poles of $R$,
 $\mu _l, \dots, \mu _n$ are poles of $R^{-1};$ 
and $k\times n$ matrices $F_{\cal P}$ and $F_{\cal N}$
 (called the left pole and  zero semi-residual matrices) 
and $n\times k$ matrices $G_{\cal P}$ and $G_{\cal N}$
 (called the right pole and  zero semi-residual matrices) 
which can be introduced from the additive decompositions
 $R(z)=\mbox{$ R(\infty)+F_{\cal P}(zI-A_{\cal P})^{-1}G_{\cal P},$}\, \
R^{-1}(z)=R(\infty)^{-1}+F_{\cal N}(zI-A_{\cal N})^{-1}G_{\cal N}.$
The right joint system representation has the form
\mbox{$R(x)\cdot R^{-1}(y)=
I+(x-y)F_{\cal P}(xI-A_{\cal P})^{-1}(S^r)^{-1}
(yI-A_{\cal N})^{-1}G_{\cal N}$,}
the left one has the form
\mbox{\(R^{-1}(x) \cdot R(y)=
I+(x-y)F_{\cal N}(xI-A_{\cal N})^{-1}(S^l)^{-1}(yI-A_{\cal P})^{-1}
G_{\cal P}\).}
The \(n\times n\) matrices $S^r$ and $S^l$ (the so-called right and left
zero-pole coupling matrices for $R$)
are solutions of the appropriate Sylvester-Lyapunov equations. 
These matrices are mutually inverse: \(S^r\cdot S^l=S^l\cdot S^r=I .\)

These results are essentially not new: they could be easily derived from known
results on realization of a rational matrix functions
(for example, from results by L. Sakhnovich or J. Ball, I. Gohberg, L. Rodman),
however the method is new, as well as the emphasis on ``the left, the right
and their relationships''.
 The presentation is oriented to a
``traditional'' analyst. No previous knowledge in realization theory
of matrix functions
or its ideology is assumed. One of the purposes of this paper
is to provide a realization theory background for investigations
of the deformation theory of Fuchsian differential system and
of rational solutions of the Schlesinger system. As an application
we also consider the spectral (Wiener-Hopf) factorization.

The concluding Section 5 contains some historical remarks highlighting
the role of M.S.\,Liv\v{s}ic as the forefather of the system realization theory.
\end{abstract}
\setcounter{section}{-1} 

\noindent
NOTATIONS.\\
\(\bullet\) $\Bbb C$  stands for the complex plane;
$\overline{\Bbb C}$ is the extended complex plane:
$\overline{\Bbb C}\stackrel{\rm\tiny def}{=}{\Bbb C}\cup{\infty};$\\
\(\bullet\) ${\Mgot}_k$ stands for 
the set of all $k\times k$ matrices with complex entries;\\
\(\bullet\) \(I\)  stands for the unity
matrix of the appropriate dimension;\\
\(\bullet\) ${\cal R} ({\Mgot}_k)$ stands for the set 
of all rational ${\Mgot}_k$-valued  functions $R$ with 
 $\mbox{\rm det}R(z)\not\equiv 0$.\\
\(\bullet\)
${\cal P}(R)$ stands for the set of all poles of the function $R$,
${\cal N}(R)$ stands for the set of all poles of the function $R^{-1}$;
${\cal P}(R)$ is said to be {\sf the pole set of the function $R$},
${\cal N}(R)$ is said to be {\sf the zero set of the function $R$}.
\newpage


\begin{minipage}{15.0cm}
\section{\hspace{-0.4cm}.\hspace{0.2cm} 
PREFACE}
\end{minipage}\\[0.18cm]
\setcounter{equation}{0}
 
The problem which we set as a goal in this paper for scalar
(i.e. complex valued) functions means 
{\sl to restore a rational function
from its poles and zeros}. The traditional solution of this problem
uses products constructed from the poles and zeros of the function.
For rational functions {\it in general position}, this solution can
be explained particularly clearly.

 Namely, let \(r\) be a rational
function in general position\,\footnote{This means that all poles and zeros of
the function \(r\) are simple and \(r(\infty)\not=0,\,\infty\).},
with the pole set \({\cal P}(r)\) and the zero set \({\cal N}(r)\).
These sets \({\cal P}(r)\) and \({\cal N}(r)\) 
do not intersect (i.e. \( {\cal P}(r)\cap{\cal N}(r)= \emptyset\) ) and
are of the same cardinality: \(\#{\cal P}(r)=\#{\cal N}(r)\).
The function \(r\) admits the  representation
\begin{equation}
r(z)=c\left( \prod\limits_{\mu_l\in {\cal N}(r)}(z-\mu_l)\right)\cdot
\left( \prod\limits_{\lambda_j\in {\cal P}(r)}(z-\lambda_j)\right)^{-1},
\label{multrecov}
\end{equation} 
where \(c=r(\infty)\). This {\it multiplicative} representation recovers
the function \(r\) from its pole and zero sets and from  the value
 \(r(\infty)\).
Inversely, given two finite non-intersecting sets
\({\cal P}\) and \({\cal N}\)
 (\({\cal P},\,{\cal N}\in\Bbb C,\,{\cal P}\cap{\cal N}=\emptyset\))
of the same cardinality
and a complex number \(c\not=0,\,\infty\), we {\sf define} the function
\(r\) by the formula (\ref{multrecov}). This function \(r\) is a rational
function in general position, the given sets \({\cal P}\) and \({\cal N}\) 
are its pole and zero sets \({\cal P}(r)\) and \({\cal N}(r)\) and 
 \(c=r(\infty)\).

However, in view of non-commutativity of the matricial multiplication,
the multiplicative representation (\ref{multrecov}) seems to be
unsuitable for generalization to  {\sf matrix} functions.
We present now such a representation of a rational matrix function
(in general position) from its poles and zeros which can be generalized to
the matricial case. This is the so-called {\sf system representation}
of a rational function.

So, let again \(r\) be a rational function in general position, with the pole
 and zero sets
\({\cal P}(r)\) and \({\cal N}(r)\). We derive its system representation.
We assume for simplicity that the function \(r\) is normalized by the condition
\(r(\infty)=1\).  We start from the additive decomposition of the matrix
function \(r\):
\begin{equation}
r(z)=1+\sum\limits_{\lambda_q\in{{\cal P}(r)}}\frac{\xi_q}{z-\lambda_q}.
\label{adddecomp}
\end{equation}
The condition 
\begin{equation}
r(\mu_p)=0 \qquad (\forall\mu_p \in {\cal N}(r))
\label{zerocond}
\end{equation}
leads to the system of linear equations
\begin{equation}
\sum\limits_{\lambda_q\in{{\cal P}(r)}}\frac{\xi_q}{\mu_p-\lambda_q}=-1  \qquad
(\forall\mu_p \in {\cal N}(r)).
\label{zeroeqw}
\end{equation}
Thus, to restore a rational function in general position from its poles and
zeros, we have to solve the linear system (\ref{zeroeqw}) with respect to
\(\xi_q\) and then to substitute these \(\xi_q\) into (\ref{adddecomp}). Since
\(\#{\cal P}(r)=\#{\cal N}(r)\stackrel{\rm\tiny def}{=}n\), the matrix \(S\)
of the  system (\ref{zeroeqw}) is square:
\begin{equation}
S=\|s_{p,q}\|_{1\leq p,q\leq _{}n},\qquad s_{p,q}=\frac{1}{\mu_p-\lambda_q}.
\label{remarkablematrix}
\end{equation}
The system (\ref{zeroeqw}) is uniquely solvable: its determinant
(which is known as the {\it Cauchy determinant}) can be calculated
explicitly (see, for example, \cite{PS}, Pt.VII:\,\S 1,\,no.\,{\bf 3}). 
 From this explicit expression for the determinant it is evident that
\(\mbox{\rm det}\,S\not=0.\)

We can formulate this method of restoring of the function \(r\) from 
\({\cal P}(r)\) and 
\({\cal N}(r)\) in the matricial form. Let \(A_{\cal P}\) and \(A_{\cal N}\)
be the diagonal matrices constructed from \({\cal P}(r)\) and \({\cal N}(r)\):
\begin{equation}
A_{\cal P}=\mbox{diag}(\lambda_1,\,\lambda_2,\,\dots\, \lambda_n), 
\qquad 
A_{\cal N}=\mbox{diag}(\mu_1,\,\mu_2,\,\dots\, \mu_n).
\label{zpm}
\end{equation}
Let \(e\) be the \(n\)-row (i.e. \(1\times n\) matrix):
\begin{equation}
e=[1,\,1,\, \dots ,\, 1].
\label{e}
\end{equation}
As usual, by \(e^{\ast}\) we denote the Hermitian conjugate to \(e\): 
\(e^{\ast}\) is a \(n\)-column (i.e. \(n\times 1\) matrix).
The representation (\ref{adddecomp}) can be put down in
 the form%
\footnote{As usual, for the matrix \(M\), \(M^T\) denotes
the transpose one.} 
\begin{equation}
r(z)=1+e\,(zI-A_{\cal P})^{-1}[\xi_1,\,\xi_2,\,\dots,\,\xi_n]^{T}.
\label{m1}
\end{equation}
The system (\ref{zeroeqw}) can be presented in matricial form:
\(S[\xi_1,\,\xi_2,\,\dots,\,\xi_n]^{T}=-e^{\ast}.\) Thus,
\[ [\xi_1,\,\xi_2,\,\dots,\,\xi_n]^{T}=-S^{-1}e^{\ast}.\]
Substituting this expression for \([\xi_1,\,\xi_2,\,\dots,\,\xi_n]^{T}\)
into (\ref{m1}), we come to the representation
\begin{equation}
r(z)=1-e\,(zI-A_{\cal P})^{-1}S^{-1}e^{\ast}.
\label{sr1}
\end{equation}
In the same way we can obtain the representation for the function \(r^{-1}\).
Starting from the additive representation
\begin{equation}
r^{-1}(z)=1+\sum\limits_{\mu_p\in{{\cal N}(r)}}\frac{\eta_p}{z-\mu_p},
\label{addreprinv}
\end{equation}
we come to the linear system (with respect to \(\eta_p\))
\[
\sum\limits_{\mu_p\in{{\cal N}(r)}}\frac{\eta_p}{\mu_p-\lambda_q}=1 \qquad
(\forall \lambda_q \in {\cal P}(r)),
\]
which can be put down in  matricial form:
\[
[\eta_1,\,\eta_2,\, \dots,\,\eta_n]S=e,
\]
or
\[
[\eta_1,\,\eta_2,\, \dots,\,\eta_n]=e\,S^{-1},
\]
Thus,
\begin{equation}
r^{-1}(z)=1+eS^{-1}(zI-A_{\cal N})^{-1}e^{\ast}.
\label{sr2}
\end{equation}
Here the matrix \(S\) is the same as in (\ref{remarkablematrix}) !\\
The matrices \(S,\,A_{\cal P},\,A_{\cal N}\) are involved in the equality
\begin{equation}
A_{\cal N}S-SA_{\cal P}=e^{\ast}e,
\label{lyapsylv}
\end{equation}
which can be directly obtained from (\ref{remarkablematrix}), (\ref{zpm}),
(\ref{e}).
{\sf The equality (\ref{lyapsylv}) and its generalization are of fundamental 
importance in the elaborated theory}.

 The equality of the form
\(AX-XB=C\) (where \(A,\,B,\,C\) are square matrices )
is known as the {\sf Sylvester-Lyapunov equality} (or as the
{\sf Sylvester-Lyapunov equation}, if it is considered as an equation with
respect to \(X\)).

Multiplying the right hand sides of the representations (\ref{sr1}) and
 (\ref{sr2}) term by term,
we obtain the equality
\begin{equation}
r(x)r^{-1}(y)=
1+(x-y)\,e\,(xI-A_{\cal P})^{-1}S^{-1}(yI-A_{\cal N})^{-1}e^{\ast }
\label{jsr}
\end{equation}
 after some calculations. These calculations are based on the 
Sylvester-Lyapunov equality (\ref{lyapsylv}).
The representations (\ref{sr1}) and (\ref{sr2})  are exactly what
we need. They are said to be {\sf the system representations} of the function
\boldmath \(r\)\unboldmath \,
and \boldmath\(r^{-1}\)\unboldmath respectively.
 The representation (\ref{jsr}) is said to be
{\sf the joint system representation of the pair of 
{\rm (mutually inverses)} functions \boldmath\(\boldmath r\)\unboldmath \,
 and  \boldmath\(r^{-1}\)\unboldmath}.

 We obtained the joint system realization (\ref{jsr}) from the system
 representations (\ref{sr1}) and (\ref{sr2}). In its turn, the representations
 (\ref{sr1}) and (\ref{sr2}) are contained
 in (\ref{jsr}):
(\ref{sr1}) is (\ref{jsr}) for \(y=\infty\), (\ref{sr2}) is (\ref{jsr}) for
 \(x=\infty \).

{\sf Now we derive the joint system representation (\ref{jsr}) in a different
 way.} Let 
\(k\) be the function 
\begin{equation}
k(x,y)=\frac{r(x)r^{-1}(y)-1}{x-y}
\label{k}
\end{equation}
of two variables. Fixing \(y\notin ({\cal P}(r)\cup{\cal N}(r)\cup\infty)\),
we  consider \(k\) as 
a function of the variable \(x\). This function is rational (with respect to
\(x\)), with simple poles located at the points \(\lambda_p\in{\cal P}(r)\),
with the residues
\({\rho}_{\lambda_p(y)}\):
\begin{equation}
\rho_{\lambda_p(y)}=\frac{1}{({r^{-1}})^{\,\,\prime}(\lambda_p)}\,\cdot\,
\frac{r^{-1}(y)}{\lambda_p-y}.
\label{rho}
\end{equation}
It is clear, that \(k(\infty,y)=0\) and that \(k(x,y)\) is holomorphic for
 \(x=y\). Thus,
\(k(x,y)\) admits the simple fraction expansion 
\begin{equation}
k(x,y)=
\sum\limits_{1\leq p\leq n}\frac{1}{x-\lambda_p}\,\cdot\,\rho_{\lambda_p(y)}\,
\cdot
\label{sfe}
\end{equation}
In its turn, the residue \(\rho_{\lambda_p(y)}\), considered as a function of
 \(y\), is a
rational one. Since \(r^{-1}(\lambda_p)=0\), the point \(\lambda_p\) is not a
 pole of the
function \(\rho_{\lambda_p}\), and the points \(\mu_q\in{\cal N}(r^{-1})\) are
 the only poles
of the function \(\rho_{\lambda_p}\). All these poles are simple, with the
 residues
\(h_{p,q}=\mbox{\rm res}\,{\rho}_{\lambda_p}(y)_{\,\,|y=\mu_q}\):
\begin{equation}
h_{p,q}=\frac{1}{({r^{-1}})^{\,\,\prime}(\lambda_p)}\cdot
 \frac{1}{\lambda_p - \mu_q}\cdot
\frac{1}{{r}^{\,\,\prime}(\mu_q)}\cdot
\label{hpq}
\end{equation}
As \(\rho_{\lambda_p}(\infty)=0\), the function \(\rho_{\lambda_p}\) admits
the simple fraction expansion 
\begin{equation}
\rho_{\lambda_p}(y)=
\sum\limits_{1\leq q\leq n} h_{p,q}\,\,\frac{1}{y-\mu_q}\cdot
\label{rhoexp}
\end{equation}
Substituting (\ref{rhoexp}) into (\ref{sfe}), we come to the formula
\begin{equation}
k(x,y)=\sum\limits_{1\leq p,q \leq n}\,\,
\frac{1}{x-\lambda_p}\cdot h_{p,q}\cdot\frac{1}{y-\mu_q},
\label{doublexp}
\end{equation}
where \(h_{p,q}\) are defined by (\ref{hpq}). The last formula can be
presented in matricial form. Let us organize the numbers \,\(h_{p,q}\)\,
into the matrix \begin{equation}
H=\|h_{p,q}\|_{1\leq p,\,q\leq n}.
\label{Hpq}
\end{equation}
Thus,
\begin{equation}
k(x,y)=e(xI-A_{\cal P})^{-1}H(yI-A_{\cal N})^{-1}e^{\ast}.
\label{krepr}
\end{equation}
and
\begin{equation}
r(x)r^{-1}(y)=1+e(xI-A_{\cal P})^{-1}H(yI-A_{\cal N})^{-1}e^{\ast}.
\label{jsrm}
\end{equation}
The comparison of the formulas (\ref{jsr}) and (\ref{jsrm}) suggests us that
 \(H=S^{-1}\). The equality
\begin{equation}
HS=I
\label{HS}
\end{equation}
may be verified starting from the formula
\begin{equation}
r^{-1}(x)=1+\sum\limits_{1\leq q \leq n}\frac{1}%
{x-\mu_q}\cdot\frac{1}{r^{\prime}(\mu_q)}\cdot
\label{int}
\end{equation}
The matrix identity (\ref{HS}) is equivalent to the system of \(n^2\) scalar
identities
\begin{equation}
\sum\limits_{1\leq q\leq n}h_{p,q}\,s_{q,t}=0
\qquad (1\leq p,t \leq n,\,\, p\not= t)
\label{pnott}
\end{equation}
and 
\begin{equation}
\sum\limits_{1\leq q \leq n}h_{p,q}\,s_{q,p}=1\qquad (1\leq p\leq n).
\label{pist}
\end{equation}
From (\ref{int}) (and of course, from  (\ref{remarkablematrix}) and 
 (\ref{hpq})) it follows that%
\footnote{\ \ \(\displaystyle\frac{1}{\lambda_p-\mu_q}\cdot \frac{1}%
{\mu_q-\lambda_t}=
\frac{1}{\lambda_p-\lambda_t}
\cdot
\left(\frac{1}{\lambda_p-\mu_q}-\frac{1}{\lambda_t-\mu_q}\right)\).}
\begin{eqnarray}
\sum\limits_{1\leq q \leq n}h_{p,q}\,s_{q,t}=\sum\limits_{1\leq q \leq n}
\frac{1}{\lambda_p-\mu_q}\cdot\frac{1}{\mu_q-\lambda_t}\cdot
\frac{1}{{(r^{-1})}^{\,\prime}(\lambda_p)}\cdot\frac{1}{r^{\prime}(\mu_q)}=
\nonumber\\
\frac{1}{{(r^{-1})}^{\,\prime}(\lambda_p)(\lambda_p-\lambda_t)}
\Big(r^{-1}(\lambda_p)-r^{-1}(\lambda_t)\Big)=0 \qquad \mbox{for}\,\, p\not= t,
\label{ver}
\end{eqnarray}
since \(r^{-1}(\lambda_k)=0\) for all
 \(\lambda_k\in{\cal P}(r)(={\cal N}(r^{-1}))\). Thus, (\ref{pnott}) is
 verified.
Analogously, the equality (\ref{pist}) may be verified using the formula
\begin{equation}
({r}^{-1})^{\,\prime}(x)=-\sum\limits_{1\leq q\leq n}\frac{1}{(x-\mu_q)^2}\,%
\cdot\,\frac{1}{r^{\prime}(\mu_q)}.
\label{der}
\end{equation}
Thus, the equality (\ref{HS}) holds. Finally, the formula (\ref{jsrm}) can
 be presented in 
the form (\ref{jsr}), with \(S\) of the form (\ref{remarkablematrix}). 

\noindent

\begin{rem}
By the way, we obtained the following rule for the inversion of the matrix 
\(S\) of the form
(\ref{remarkablematrix}), where
\(\lambda_1,\,\lambda_2,\,\dots,\,\lambda_n;\,\,
\mu_1,\,\mu_2,\,\dots,\,\mu_n\) are pairwise different
complex numbers  (this matrix is known as the Cauchy matrix).  Starting from
these numbers, we have to construct the rational function \(r\) of the form
(\ref{multrecov}) (with some \( c:\, c\neq 0,\, \infty\)).
The entries \(h_{p,q}\) of the inverse matrix \(H=S^{-1}\) are of the form
 (\ref{hpq}).
The expression (\ref{hpq}) for the inverse matrix \(S^{-1}\) may be
represented in the matricial form:

\vspace*{-3ex}
\begin{equation}
S^{-1}=D_1SD_2,\quad
 {\rm where}\ D_1={\rm diag}\left(\frac{1}{(r^{-1})^{\prime}(\lambda_l)}
 \right), \
D_2={\rm diag}\left(\frac{1}{r^{\prime}(\mu_l)} \right).
\label{matrform}
\end{equation}
 In particular, we established the invertibility of the Cauchy matrix \(S\).
 Of course, the expression (\ref{hpq}) for the entries of the
inverse matrix \(S^{-1}\) differs from the expression obtained in
 (\cite{BGR1}) in  form only.
In (\cite{BGR1}) (Lemma 6.3 there)  this expression was obtained by means of
 the Kramer  matrix
 inversion rule, using
explicit expression for the Cauchy determinant. (The determinant of the matrix
\(S\) as well as its cofactors are  Cauchy determinants).
\end{rem}
The square of the Cauchy determinant can be calculated from (\ref{matrform}):
\begin{equation}
\left({\rm det}S\right)^2=
\prod_{1\leq p \leq n} (r^{-1})^{\prime}(\lambda_p)\,\cdot\,
\prod_{1\leq q \leq n} r^{\prime}(\mu_q).
\label{squaredet}
\end{equation}
\begin{rem} \hspace{-1.5ex} This derivation
 of the joint system representation can be extended to some classes of 
{\sf meromorphic} functions.
\end{rem}

\vspace{3pt}
Obtaining ``system'' representations of ``scalar'' rational functions, we made
 no use of anything
 that can not be generalized to the matricial case.
Generalizations  of the representations (\ref{sr1}) and (\ref{sr2}) to rational
{\sf matrix} functions have been done already.
(See \cite{S3},\,\cite{GKLR}, \cite{BGR1}, \cite{BGR2}, \cite{BGRa}).  

The term ``system representation'' is related to the system theory.
 In this theory to each linear stationary 
time invariant dynamical system its transfer function is corresponded,
 which is matrix valued (and even operator valued) one.
 If the state space of the system is finite dimensional,
then this transfer function is rational. And if the system is a SISO system
(Single Input, Single Output), then its transfer function is scalar (complex
valued). This transfer function is expressed in terms of the input, output and
 state space operators of the system. Inversely, each rational matrix function
 may be represented as the transfer function of an appropriate linear time
 invariant system with finite dimensional state space.
Such a representation\,\footnote{It looks like the representation (\ref{sr1})}
is said to be {\sf the system realization} or {\sf the system representation
 of the given rational function}.

The relationship between  rational functions and linear systems can be
exploited in both directions. We may apply results from matrix functions
theory for study of linear systems. And we may  use the system representation
as a tool for study of rational matrix functions or as a tool to specify them.

It should be mentioned that it was M.\,Liv\v{s}ic who discovered
the characteristic  matrix function of linear operators and operator
colligations (=operators nodes). He was the first to relate invariant
subspaces of a linear operator and factors of its characteristic matrix
 function, \cite{L2}, \cite{L3}, \cite{LP}.
In the system theory language, the characteristic function of an operator node
is just the transfer function of the appropriate linear time invariant system.
This was shown by M.S.\,Liv\v{s}ic in \cite{L8}, \cite{L9}, \cite{BrL}.
He has also related the characteristic matrix function with the scattering
theory (\cite{L5},\,\cite{L6},\,\cite{L9},\,\cite{BrL}).

The main goal of this paper is to present some basic results on the  system
representation of rational matrix function in a self-contained form.
For the clarity of presentation, we restrict our consideration to the simple
but important case of a rational matrix function in general position. The
presentation is oriented to a ``traditional'' analyst. No previous knowledge in
realization theory of matrix functions or its ideology is assumed. 
We leave detailed historical remarks for Section 5.

\vspace{4.0ex}
\begin{minipage}{15.0cm}
\section{\hspace{-0.4cm}.\hspace{0.2cm} 
RATIONAL MATRIX FUNCTIONS IN GENERAL POSITION}
\end{minipage}\\[-0.18cm]
\setcounter{equation}{0}

\begin{defi}  
{\sl Let $M(z)$ be a $k\times k$ matrix function holomorphic in a punctured 
neighborhood of some point $a$ (i.e. the point $a$
is an isolated singularity of the function $M$). The point $a$ is said to be
{\sf a simple pole of the matrix function} $M$ if 
\begin{equation}
M(z)=\frac{M_a}{z-a}+ H(z),
\end{equation}
where $M_a$ is a constant matrix and the matrix function $H$ is holomorphic at
the point $a$.
The matrix $M_a$ is said to be {\sf the residue of the matrix function $M$ at
the point $a$}.
}
\end{defi}

\begin{defi} 
{\sl 
Let a point $a$ be a simple pole of a $k \times k$ matrix function $M$ and let
the residue $M_a$ of the function $M$ at the point $a$ be a matrix of rank
one.
As a matrix of rank one, the $k \times k$ matrix $M_a$ can be factorized in
 the form
\begin{equation}  
M_a=f_a\cdot g_a
\end{equation}
where $f_a$ is \(k\) vector-columns (i.e. \(k\times 1\) matrix), \(g_a\) is
\(k\) vector-row (i.e. \(1\times k\) matrix),
$f_a\not= 0, g_a \not= 0$.
The vectors $f_a$ and $g_a$ are said to be respectively 
{\sf the left semi-residual vector {\sl and} the right semi-residual vector at
 the point $a$}.
The vectors $f_a$ and $g_a$ are defined uniquely, up to a  constant factor:
we can represent the number $1$ in the form $1=d\cdot d^{-1}$  and then 
redistribute the factors $d$ and $d^{-1}$:
\begin{equation}
f_a\rightarrow  f_a\cdot d;\quad
g_a\rightarrow  d^{-1}\cdot g_a. 
\end{equation}
}
\end{defi}\\[0.05cm]
We emphasize that {\sf the notions of  left and right semi-residual vectors
 are 
defined only for a simple pole with residue of rank one}.

\begin{defi}\label{genpos} 
 {\sl A rational matrix function $R$
 {\rm (}$\,R\in{\cal R} ({\Mgot}_k)\, ${\rm )}\  is
said to be {\sf a rational matrix function in general position} if:}
\begin{enumerate}
\setlength{\topsep}{0.1cm}
\setlength{\itemsep}{-0.1cm}
{\rm
\item
{\sl
 The pole set ${\cal P}(R)$ and the zero set ${\cal N}(R)$ do not
intersect;
\begin{equation}
{\cal P}(R)\cap {\cal N}(R)=0;
\end{equation}}
\item
{\sl
 All poles of the function $R$ are simple, and the residues at
these poles are  matrices
of rank one;}
\item
{\sl
 All poles of the function $R^{-1}$ are simple, and the
 residues at these poles are  matrices of rank one;}
\item
{\sl
 Both functions $R$ and $R^{-1}$ are holomorphic at the point
$z=\infty$.}}
\end{enumerate}
\end{defi}
Let $R$ be a rational matrix function in general position, with the
pole set  ${\cal P}(R)$ and the zero set ${\cal N}(R)$.
 For $\lambda \in {\cal P}(R)$,  $\mu \in {\cal N}(R)$,
let
\begin{equation}
R_{\lambda}=f_{\lambda}\cdot g_{\lambda},
 \quad (f_{\lambda}\not= 0, g_{\lambda}\not= 0 );
 \qquad
R_{\mu}=f_{\mu}\cdot g_{\mu} \quad (f_{\mu}\not= 0,  g_{\mu}\not= 0 )
\label{factres}
\end{equation}
 be  factorizations of the residue  $R_{\lambda}$
of the function $R$ at the point $\lambda$ and 
the residue  $R_{\mu}$ of the function $R^{-1}$ at the point $\mu$
respectively, where $f_{\lambda}, g_{\lambda};
f_{\mu}, g_{\mu}$ are the appropriate semi-residual
 vectors.
The additive expansions
\begin{equation}
  R(z)=R(\infty)+
       \sum\limits_{\lambda\in{\cal P} (R)}\frac{R_{\lambda}}{z-\lambda}
\label{addnonfact_a} 
\end{equation}
\begin{equation}
   R^{-1}(z)=R^{-1}(\infty)+
              \sum\limits_{\mu\in{\cal N} (R)}\frac{R_{\mu}}{z-\mu}
\label{addnonfact_b}
 \end{equation}
can be rewritten in the form
\begin{equation}
R(z)=R(\infty)+
\sum\limits_{\lambda\in{\cal P} (R)}f_{\lambda}\cdot \frac{1}{z-\lambda}\cdot
g_{\lambda},
\label{addexpdir}
\end{equation}
\begin{equation}
R^{-1}(z)=R^{-1}(\infty)+
\sum\limits_{\mu\in{\cal N} (R)}f_{\mu}\cdot \frac{1}{z-\mu}\cdot
g_{\mu},
\label{addexpinv}
\end{equation}

\begin{lemma}\label{reslogder}
\hspace*{0.5em}
I.\ {\sl Let a matrix function $R$ be analytic in a punctured
neighborhood of a point
 $\lambda\in\Bbb C$, $\lambda$
 be a simple pole of the function $R$ and a holomorphicity
point of the function $R^{-1}$, and let $R_{\lambda}$ be the residue of the
 function $R$ at the point $\lambda$.
Then the point $\lambda$ is a simple pole for the 
``logarithmic derivative'' $R^{\prime}\cdot R^{-1}$, and for 
the residue $P_{\lambda}$ of this ``logarithmic derivative'' at this point
 the conditions
\begin{equation}
P_{\lambda}^2=-P_{\lambda},\quad 
\mbox{\rm rank}P_{\lambda}=\mbox{\rm rank}R_{\lambda}
\end{equation}
hold; in particular,
\begin{equation}
\mbox{\rm trace}P_{\lambda}= -\mbox{\rm rank}R_{\lambda}.
\end{equation}
}
II.\  {\sl Let a matrix function $R^{-1}$ be analytic in a punctured
neighborhood of
 a point
 $\mu\in\Bbb C$, $\mu$
 be a simple pole of the function $R^{-1}$ and a holomorphicity
point of the function $R$, and let $R_{\mu}$ be the residue of the function
 $R^{-1}$ at the point $\mu$.
Then the point $\mu$ is a simple pole for the 
``logarithmic derivative'' $R^{\prime}\cdot R^{-1}$, and for 
the residue $P_{\mu}$ of this ``logarithmic derivative'' at this point
 the conditions
\begin{equation}
P_{\mu}^2=P_{\mu},\quad 
\mbox{\rm rank}P_{\mu}=\mbox{\rm rank}R_{\mu}
\end{equation}
hold; in particular,
\begin{equation}
\mbox{\rm trace}P_{\mu}= \mbox{\rm rank}R_{\mu}.
\end{equation}
}
\end{lemma}

\noindent
PROOF. We prove only  statement I of Lemma. Statement II can be proved
analogously. Let
\begin{equation}
R(z)=\frac{R_{\lambda}}{z-\lambda}+A_0+A_1(z-\lambda )+A_2(z-\lambda)^2
+\cdots
\label{lorR}      
\end{equation}
and
\begin{equation}
R^{-1}(z)=B_0+B_1 (z-\lambda )+B_2(z-\lambda )^2 + \cdot
\label{lorR^{-1}}      
\end{equation}
be the Laurent expansions of the functions $R$ and $R^{-1}$ respectively.
Then
\begin{equation}
R^{\prime}(z)=-\frac{R_{\lambda} }{(z-\lambda )^2}+
A_1 + 2A_2(z-\lambda )+
\dots
\label{deri}
\end{equation}
Multiplying the Laurent expansions term by term, we obtain from 
 (\ref{lorR^{-1}}) and (\ref{deri})
\begin{equation}
R^{\prime}(z)\cdot R^{-1}(z)=
-\frac{R_{\lambda}B_0}{(z-\lambda)^2}
-\frac{R_{\lambda}B_1}{z-\lambda}
+(-R_{}B_2+A_1B_0)
+O(z-\lambda ).
\label{lder}
\end{equation}
Substituting the Laurent expansions (\ref{lorR}),
 (\ref{lorR^{-1}}) into the identity 
$R^{-1}(z)R(z)=I$
and multiplying these expansions term by term, we obtain
\begin{equation}
 B_0R_{\lambda}=0
\label{cons1a} 
\end{equation}
\begin{equation}
 B_1R_{\lambda}+B_0A_0=I.
\label{cons1b}
\end{equation}
Analogously, from the identity
$R(z)R^{-1}(z)=I$
we derive
\begin{equation}
 R_{\lambda}B_0=0
\label{cons2a}
\end{equation}
\begin{equation}
R_{\lambda}B_1 + A_0B_0=I.
\label{cons2b}
\end{equation}
Let us examine the expansion (\ref{lder}). According to (\ref{cons2a}),
the term
 $\displaystyle -\frac{R_{\lambda}B_0}{(z-\lambda)^2}$
vanishes. Thus, the point $\lambda$ is a simple pole for the function
$R^{\prime}(z)R^{-1}(z)$, with the residue $P_{\lambda}$,
\begin{equation}
P_{\lambda}=-R_{\lambda}B_1.
\label{res}
\end{equation}
>From (\ref{cons2b}) and (\ref{res}) it follows that $I+P_{\lambda}=A_0B_0$.
Hence,
\[(I+P_{\lambda})P_{\lambda}=
(A_0B_0)\cdot (-R_{\lambda}B_1)=
-A_0(B_0R_{\lambda})B_1.\]
According to (\ref{cons1a}), $B_0R_{\lambda}=0$. Thus
\((I+P_{\lambda})P_{\lambda}=0\), i.e. $P_{\lambda}^2 = - P_{\lambda}.$
 Finally,\\
$P_{\lambda}R_{\lambda}=-(I-A_0B_0)R_{\lambda}=
-R_{\lambda}+A_0(B_0R_{\lambda}).$ Since $B_0R_{\lambda}=0$,
\begin{equation}
R_{\lambda}=-P_{\lambda}R_{\lambda} 
\label{rank}
\end{equation}
>From (\ref{res}) and (\ref{rank}) it follows that
 $\mbox{\rm rank} P_{\lambda}=\mbox{\rm rank} R_{\lambda} .$ 
 \hfill\framebox[0.45em]{ }

\begin{rem} \label{usid}
>From (\ref{cons1b}) and (\ref{cons2a}) it follows that
 $R_{\lambda}B_1R_{\lambda}=R_{\lambda}$. Since
 $B_1=\big( R^{-1}\big)^{\prime}(\lambda)$ (see \ref{lorR^{-1}}), it can be
written as
\begin{equation}
R_{\lambda}\,  \big( R^{-1}\big)^{\prime}(\lambda)\,
 R_{\lambda}=R_{\lambda}\quad
(\lambda\in\cal P ).
\label{rel1}
\end{equation}
Analogously, we derive
\begin{equation}
R_{\mu}\,   R^{\, \prime}(\mu)\,
 R_{\mu}=R_{\mu}\quad
(\mu\in\cal N ).
\label{rel2}
\end{equation}
We shall use the relations (\ref{rel1}) and (\ref{rel2}) in the following
 section.
\hfill\framebox[0.45em]{ }
\end{rem}

\begin{lemma} \label{isocard}
{\sl For a rational matrix function $R$ in general position,
\begin{equation}
\# {\cal P} (R)=\# {\cal N} (R)\ \  \ 
\big(\stackrel{\mbox{\rm\tiny def}}{=} n(R)\big).
\label{eqcard}
\end{equation}
}
\end{lemma}
\mbox{}\vspace*{-1.4em}

\noindent
PROOF. To prove the statement of Lemma, we have to consider, in some way or
 another, the ``logarithmic derivative'' \(R^{\prime}(z)\cdot R^{-1}(z)\),
(or \(R^{-1}(z)\cdot R^{\prime}(z)\)). For a  rational matrix function $R$ in
general position, its logarithmic derivative
\(R^{\prime}\cdot R^{-1}\)
 admits the expansion
\begin{equation}
R^{\prime}(z)\cdot R^{-1}(z)
=\sum\limits_{\lambda\in{\cal P}(R)}\frac{P_{\lambda}}{z-\lambda}+
\sum\limits_{\mu\in{\cal P}(R)}\frac{P_{\mu}}{z-\mu}.
\label{logderexp}
\end{equation} 

Indeed, the (rational) function
 \(R^{\prime}(z)\cdot R^{-1}(z) \)
 may have
singularities only at the points of the set \({\cal P}(R) \cup {\cal N}(R)\).
According to item 4 of Definition \ref{genpos}, the logarithmic derivative 
vanishes at the point \(\infty\);  moreover,
\begin{equation}
R^{\prime}(z)\cdot R^{-1}(z)=O(|z|^{-2})\quad (z\rightarrow\infty )
\label{lderinf}
\end{equation}
According to Lemma \ref{reslogder}, all singularities of the logarithmic
 derivative are simple poles with  residues
 $P_{\lambda}$ and $P_{\mu}$  of rank one 
 ($\lambda \in {\cal P}(R), \, \mu \in {\cal N}(R)$ respectively).
Thus, the expansion (\ref{logderexp}) holds, with
 $\mbox{\rm trace} { P}_{\lambda}=-1,$
 $\mbox{\rm trace} { P}_{\mu}=1.$ From (\ref{logderexp}) and 
(\ref{lderinf}) it follows that 
\begin{equation}
\displaystyle
\sum_{\lambda\in {\cal P}(R)} P_{\lambda}+
\sum_{\mu\in {\cal N}(R) } P_{\mu}=0.
\label{card}
\end{equation}
Because, according to Lemma \ref{reslogder},
 \(\mbox{\rm trace} P_{\lambda} =-1,\ (\lambda \in {\cal P}(R)),\ 
\mbox{\rm trace} P_{\mu}        =1,\ (\mu \in {\cal N}(R)),\) from (\ref{card})
it follows that 
\(
\displaystyle
\sum_{\lambda\in {\cal P}(R)}(- 1) + \sum_{\mu\in {\cal N}(R)} 1 =0.\)
This equality coincides with (\ref{eqcard}).
Of course, this result could be obtained using the operator version of 
Rouch\'{e} theorem from \cite{GS}.
\hfill\framebox[0.45em]{ }


\begin{defi}
{\sl
For a rational \( k\times k\)
 matrix function \( R \) in general position, let us
order its poles and ``zeros'' somehow:
 \(
{\cal P}=\{\lambda _1, \dots, \lambda _n \},
 {\cal N}=\{\mu _1, \dots, \mu _n \}.
\)
{\rm (}We remind that \mbox{\#{ \cal P}(R)=\# {\cal N}(R)}.{\rm)} We introduce
the \(n\times n\) diagonal matrices
\begin{equation}
A_{\cal P}=\mbox{\rm diag}(\lambda _1, \dots , \lambda _n),
\quad
A_{\cal N}=\mbox{\rm diag}(\mu_1, \dots, \mu_n),
\label{polematr}
\end{equation}
 \( k\times n \) matrices \(F_{\cal P},\,F_{\cal N}\) and \( n\times k \)
matrices \(G_{\cal P},\,G_{\cal N}\):
{}\\[-2em]
\begin{eqnarray}
F_{\cal P} = [ f_{\lambda_ 1}, \dots , f_{\lambda_ n} ] ,&
\quad
&F_{\cal N} = [ f_{\mu_ 1}, \dots , f_{\mu_ n} ]
,
\label{leftsemires}
\\
G_{\cal P} = \left[
\begin{array}{c}
 g_{\lambda_1}
\\
 \vdots 
\\
 g_{\lambda_n}
\end{array}\right],&
\quad
&
G_{\cal N} = \left[
\begin{array}{c}
g_{\mu_ 1}
\\
 \vdots 
\\
 g_{\mu_ n}
\end{array}
\right],
\label{rightsemires}
\end{eqnarray}
where
 \( f_{\lambda_j},\, \, g_{\lambda_j}\ (\lambda _j \in {\cal P} ) \) 
are the left and right semi-residual vectors at the pole \( \lambda _j \) of
\( R \)\,, \
 \( f_{\mu_j},\, \, g_{\mu_j}\ (\mu _j \in {\cal N} ) \) 
are the left and right semi-residual vectors at the pole \( \mu _j \) of
\( R ^{-1}\)\,.

The matrices \( A_{\cal P}\) and \( A_{\cal N}\) are said to be 
{\sf the pole matrix} and {\sf the zero matrix} respectively
 {\sf for the matrix function} \(R. \)

The matrices \( F_{\cal P} \) and \( G_{\cal P} \) are said to be 
{\sf the left-} and {\sf the right semi-residual matrices
 corresponding to the pole set \( {\cal P}(R)\)}.

The matrices \( F_{\cal N} \) and \( G_{\cal N} \) are said to be 
{\sf the left-} and {\sf the right semi-residual matrices
corresponding to the zero set  \( {\cal N}(R) \) }.
}
\hfill\framebox[0.45em]{ }
\end{defi}

\begin{rem}\label{gauge} 
It should be mentioned that if we order somehow the poles and the zeros,
then the pole and the zero matrices \(A_{\cal P}\) and \(A_{\cal N}\) are
defined uniquely, and the semi-residual matrices 
\(F_{\cal P}, G_{\cal P} , F_{\cal N}, G_{\cal N}\)
 are defined essentially uniquely, up to 
multiplication by diagonal matrices with non-zero diagonal entries:
\begin{eqnarray}
F_{\cal P}\rightarrow F_{\cal P}\cdot  D_{\cal P}, &
G_{\cal P}\rightarrow D_{\cal P}^{-1}\cdot G_{\cal P},
\label{polegauge} \\ [-1.0em]
 &   \nonumber  \\
F_{\cal N}\rightarrow F_{\cal N}\cdot  D_{\cal N}, &
G_{\cal N}\rightarrow D_{\cal N}^{-1}\cdot G_{\cal N},
\label{zerogauge}
 \end{eqnarray}
{}\\[-2.3em]
where 
{}\\[-2.3em]
\begin{eqnarray}
D_{\cal P}\, \,  \,
 = &
\mbox{\rm diag}(d_{1,{\cal P}},\, \dots\, , d_{n,{\cal P}} )\quad &
(d_{j,{\cal P}}\not= 0,\  j=1, 2, \dots n ),\label{pdauge}\\
D_{\cal N}\, \, \, 
 = &
 \mbox{\rm diag}(d_{1,{\cal N}},\, \dots\, ,d_{n,{\cal N}} )\quad &
(d_{j,{\cal N}}\not= 0,\  j=1, 2, \dots n ).\label{ndauge}
\end{eqnarray}
This freedom in choice of these diagonal matrices 
 \( D_{\cal P},\, D_{\cal N} \) can be used to simplify some formulas.
Of course, for given left semi-residual matrix \( F_{\cal P}\), the right
semi-residual matrix \( G_{\cal P}\) is determined uniquely;
for given right  semi-residual matrix \( G_{\cal P}\), the left
 semi-residual matrix \( F_{\cal P}\) is determined uniquely,
 etc.
\hfill\framebox[0.45em]{ }

\end{rem}

 It is clear that the additive expansions (\ref{addexpdir}) and
(\ref{addexpinv}) may be rewritten in the matricial form
\begin{eqnarray}
R(z)&= &   \,\,\, \, R(\infty) \,\,   +  \, 
   F_{\cal P}\cdot\big(zI-A_{\cal P}\big)^{-1}
        \cdot G_{\cal P},
\label{matradddir}\\
R^{-1}(z)&= &R^{-1}(\infty)
+ F_{\cal N}\cdot\big(zI-A_{\cal N}\big)^{-1}
\cdot G_{\cal N}.
\label{matraddinv}
\end{eqnarray}

\noindent
\vspace*{0.2cm} 
\begin{minipage}{15.0cm}
\section{\hspace{-0.4cm}.\hspace{0.19cm}THE JOINT REPRESENTATION OF THE
KERNELS ASSOCIATED WITH A  RATIONAL MATRIX FUNCTION  IN GENERAL POSITION.}
\end{minipage}\\[-0.3cm]
\setcounter{equation}{0}

\begin{defi}
{\sl
Given a rational matrix function \(R\) of one variables, we associate with it
two matrix function of two variables,
\(K^{{}^{\,r}}_{{}_{R}}\)
 and
\(K^{{}^{\,l}}_{{}_{R}}\):
\begin{equation}
K^{{}^{\,r}}_{{}_R}(x,y)=\frac{R(x)R^{-1}(y)-I}{x-y}
\label{rkern}
\end{equation}
and
\begin{equation}
K^{{}^{\,l}}_{{}_R}(x,y)=\frac{R^{-1}(x)R(y)-I}{x-y}
\label{lkern}
\end{equation}
The function \(K^{{}^{\,r}}_{R}\) is said to be
 {\sf the right kernel associated with 
the function \(R\).} \\[0.3ex]
The function \(K^{{}^{\,l}}_{R}\) is said to be
 {\sf the left kernel associated with 
the function \(R\).}
}    
     \hfill\framebox[0.45em]{ }\\[0.3em]
\end{defi}

\vspace{-1.0ex}
\noindent
\begin{rem}
If \(P\) is a polynomial, the expression
\(\displaystyle B(x,y)=\frac{P(x)-P(y)}{x-y}\) ß
 is said to be the Bezoutiant of the polynomial
\( P\). The expressions (\ref{rkern}) and (\ref{lkern}) look like
a Bezoutiant.
\end{rem}
 \hfill\framebox[0.45em]{ }

\vspace{1.4ex}
\noindent
\begin{theo}\label{rightreprtheo}
{\sl
Let \(R\) be a rational \(k\times k\) matrix function in general position;
 \( A_{\cal P}\) and
\( A_{\cal N}\) be the pole and the zero matrices for \(R\);
\(F_{\cal P}\) and \(G_{\cal N}\) be the left pole and the right zero
semiresidual matrices respectively.\\
Then:{\rm
\begin{enumerate}
\setlength{\topsep}{0.1cm}
\setlength{\itemsep}{-0.1cm}
\item{\sl
 The right kernel \(K^{{}^ r}_{\scriptscriptstyle R}\)
 is representable in the form
\begin{equation}
K^{{}^ r}_{\scriptscriptstyle R}(x,y)=
F_{\cal P}\big(xI- A_{\cal P}\big)^{-1}
H^r\big(yI- A_{\cal N}\big)^{-1}G_{\cal N},
\label{rsysrepr}
\end{equation}
where  \(H^r\) is some \(n\times n\) matrix.}
\item{\sl
For given matrices
\footnote{\label{footgauge}We remind, that the matrices \(F_{\cal P},\,
 G_{\cal P},\,F_{\cal N},\, G_{\cal N}\)
are defined {\sf only up to transformations}
 (\ref{polegauge}), (\ref{zerogauge}) with {\sf arbitrary} diagonal
invertible matrices  \(D_{\cal P}\), \(D_{\cal N}\).
}
 \(F_{\cal P}\) and \(G_{\cal N}\),
the matrix \(H^r\) is defined uniquely:
\begin{equation}
H^r=\|h^r_{p,q}\|_{1\leq p,q\leq n}, \quad
h^r_{p,q}=\frac{g_{\lambda _p}\cdot f_{\mu _q}}{\lambda _p - \mu _q }.
\label{rmatrentry}
\end{equation}}
\item{\sl
The matrix  \(H^r=\|h^r_{p,q}\|\) is invertible.}
\end{enumerate} 
}
}
\end{theo}
The left version of this theorem holds as well.

\noindent\vspace{-1.0ex}

\begin{theo}\label{leftreprtheo} 
{\sl
Let \(R\) be a rational \(k\times k\) matrix function in general position;
\( A_{\cal P}\) and \( A_{\cal N}\)
be the pole and the zero matrices for \(R\);
\(F_{\cal N}\) and \(G_{\cal P}\)
be the left zero and the right pole semiresidual
 matrices respectively.\\
Then:
{\rm
\begin{enumerate}
\setlength{\topsep}{0.1cm}
\setlength{\itemsep}{-0.1cm}
\item
{\sl
 The left kernel \(K^l_{\scriptscriptstyle R}\)
 is representable in the form
\begin{equation}
K^l_{\scriptscriptstyle R}(x,y)=
F_{\cal N}\big(xI- A_{\cal N}\big)^{-1}
H^l\big(yI- A_{\cal P}\big)^{-1}G_{\cal P},
\label{lsysrepr}
\end{equation}
where  \(H^l\) is some \(n\times n\) matrix.}
\item
{\sl
For given matrices \({}^{\ref{footgauge}}\)
 \(F_{\cal N}\) and \(G_{\cal P}\),
the matrix \(H^l\)
 \(H^l=\|h^l_{p,q}\|_{1\leq p,q\leq n}\)
 is defined uniquely:
\begin{equation}
H^l=\|h^l_{p,q}\|_{1\leq p,q\leq n}, \quad
h^l_{p,q}=\frac{g_{\mu _p}\cdot f_{\lambda _q}}{\mu _p -\lambda _q}.
\label{lmatrentry}
\end{equation}}
\item
{\sl
The matrix  \(H^l\) is invertible.}
\end{enumerate} }
}
\end{theo}
\begin{defi}   \label{defcoup}           
{\sl
The matrices \(H^r\) and \(H^l\) which appear in the representations 
(\ref{rsysrepr}) and (\ref{lsysrepr}) of the kernels
 \(K^r_{\scriptscriptstyle R}\) and \(K^l_{\scriptscriptstyle R}\),
 are said to be \,
{\sf the right core matrix }and {\sf the left core matrix} respectively.
}
\hfill\framebox[0.45em]{ }\\[0.0em]
\end{defi}

\vspace{-2.5ex}
\begin{cor}\label{sreprcore}
{\sl 
Let \(R\) be a rational \(k\times k\) matrix function in general position;
\( A_{\cal P}\) and \( A_{\cal N}\)
be the pole and the zero matrices for \(R\);\,\,
 \(F_{\cal P}\),\,\(F_{\cal N}\),\,\(G_{\cal P}\)\,\(G_{\cal N}\)
 are the appropriate semi-residual matrices; \(H^r,H^l\) are
 the appropriate core matrices.

Then the matrices \(R(x)\cdot R^{-1}(y)\) and \(R^{-1}(x)\cdot R (y)\)
admit the representations
\begin{eqnarray}
 R(x)(R(y))^{-1}&
=&I+(x-y)\,
F_{\cal P}\,\big(xI- A_{\cal P}\big)^{-1}
H^r\big(yI- A_{\cal N}\big)^{-1}G_{\cal N},
 \label{rpsjoint}\\ 
 (R(x))^{-1}R(y)&=&I+(x-y\,)
F_{\cal N}\,\big(xI- A_{\cal N}\big)^{-1}
H^l\big(yI- A_{\cal P}\big)^{-1}G_{\cal P}.
  \label{lpsjoint}
\end{eqnarray}
Under the normalizing condition
\begin{equation}
R(\infty)=I,
\label{norm}
\end{equation}
the matrix functions \(R\),\,\(R^{-1}\) themselves admit the representations
\begin{eqnarray}
R(z)&=&I-F_{\cal P}\,(zI-A_{\cal P})^{-1}H^{r}G_{\cal N},
\label{forR}
\\
R^{-1}(z)&=&I+F_{\cal P}\,H^{r}(zI-A_{\cal N})^{-1}G_{\cal N},
\label{forR^{-1}}\\
R(z)&=&I+F_{\cal N}H^l(zI-A_{\cal P})^{-1}G_{\cal P},
\label{forRR^{-1}1}\\
R^{-1}(z)&=&I-F_{\cal N}(zI-A_{\cal N})^{-1}H^{l}G_{\cal P}.
\label{forRR^{-1}2}
\end{eqnarray} 
}
\end{cor}

\noindent
PROOF. The representations (\ref{rpsjoint}) and
 (\ref{lpsjoint}) are nothing more than the
representations (\ref{rsysrepr}) and (\ref{lsysrepr}) rewritten in the terms
 of the
functions \(R(x)\cdot R^{-1}(y)\) and \(R^{-1}(x)\cdot R(y)\).\\
Letting \(y\) tend to \(\infty\) in (\ref{rpsjoint}), we obtain (\ref{forR});
letting \(x\) tend to \(\infty\) in (\ref{rpsjoint}), we obtain
 (\ref{forR^{-1}});
Letting \(y\) tend to \(\infty\) in (\ref{lpsjoint}), we obtain
 (\ref{forRR^{-1}1});
letting \(x\) tend to \(\infty\) in (\ref{lpsjoint}), we obtain
 (\ref{forRR^{-1}2}).
\hfill\framebox[0.45em]{ }\\[0.0em]

\begin{theo}\label{LyapSylCore}
{\sl Let \(R\) be a rational \(k\times k\) matrix function in general position;
\( A_{\cal P}\) and \( A_{\cal N}\)
be the pole and the zero matrices for \(R\);\,\,
 \(F_{\cal P}\),\,\(F_{\cal N}\),\,\(G_{\cal P}\)\,\(G_{\cal N}\)
 are the appropriate semi-residual matrices; \(H^r,H^l\) are
 the appropriate core matrices.

Then these matrices are involved in the Sylvester-Lyapunov equalities}:
\begin{eqnarray}
A_{\cal P}H^r - H^rA_{\cal N}=G_{\cal P}F_{\cal N},
\label{LS1}
\\
A_{\cal N}H^l - H^lA_{\cal P}=G_{\cal N}F_{\cal P}.
\label{LS2}
\end{eqnarray}
\end{theo}
\noindent
PROOF. The matrices \(G_{\cal P}F_{\cal N}\) and \(G_{\cal N}F_{\cal P}\) are
of the form
\begin{equation}
G_{\cal P}\cdot F_{\cal N}=\|g_{\lambda_p}\cdot f_{\mu_q}\|_{1\leq p,q\leq n},
\qquad
G_{\cal N}F_{\cal P}=\|g_{\mu_p}\cdot f_{\lambda_q}\|_{1\leq p,q\leq n}.
\label{displ}
\end{equation}
The assertion of Theorem \ref{LyapSylCore} follows from the explicit
 expressions
(\ref{displ}), (\ref{rmatrentry}), (\ref{lmatrentry}), (\ref{polematr}).
\hfill\framebox[0.45em]{ }\\[0.0em]

\begin{theo}\label{couprelcore}
{\sl Let \(R\) be a rational \(k\times k\) matrix function in general
position, which satisfies the normalizing condition  (\ref{norm});
\(F_{\cal P}\),\,\(F_{\cal N}\),\,\(G_{\cal P}\)\,\(G_{\cal N}\)
are the appropriate semi-residual matrices; \(H^r,H^l\) are
the appropriate core matrices.

Then these matrices are involved in the equalities}
\begin{equation}
\mbox{\rm a)}.\,\,     H^rG_{\cal N}=-G_{\cal P}; \quad
\mbox{\rm b)}.\,\,     H^lG_{\cal P}=-G_{\cal N}; \quad
\mbox{\rm c)}.\,\,     F_{\cal P}H^r= F_{\cal N}; \quad
\mbox{\rm d)}.\,\,     F_{\cal N}H^l= F_{\cal P}
\label{coup}
\end{equation}
\end{theo}
We give two proofs of Theorem \ref{couprelcore}.\\[3pt]
PROOF I. We compare the formulas (\ref{matradddir}) and (\ref{forR}).
In the additive representation (\ref{matradddir}), for given 
left semi-residual matrix \(F_{\cal P}\), the right one
\(G_{\cal P}\) is determined uniquely. (See Remark \ref{gauge}).
Therefore, (\ref{coup}.a) holds. Analogously, comparing the formulas
(\ref{matradddir}) and (\ref{forRR^{-1}1}), we obtain (\ref{coup}.d).
Comparing  (\ref{matraddinv}) and (\ref{forR^{-1}}),  (\ref{forRR^{-1}2}),
we obtain (\ref{coup}.c) and (\ref{coup}.b).
\hfill\framebox[0.45em]{}\\[5pt]
PROOF II.  We proof only the equality (\ref{coup}.a). The equality 
(\ref{coup}.c) can be proved analogously. (\ref{coup}.b) follows from 
(\ref{coup}.a) and (\ref{invrel}), etc.
The matricial equality (\ref{coup}.a) is the same that the system of
\(n\) scalar equalities
\begin{equation}
\sum\limits_{1\leq q \leq n}h^r_{p,q}g_{\mu_q}=-g_{\lambda_p},
\qquad p=1,\,2\,\dots,\, n.
\label{coup1}
\end{equation}
Since \(f_{\lambda_p}\not=0\), the last equality is equivalent to the equality
\begin{equation}
f_{\lambda_p}\cdot
\Big(\sum\limits_{1\leq q \leq n}h^r_{p,q}\Big)\cdot g_{\mu_q}=
-f_{\lambda_p}\cdot g_{\lambda_p}. 
\label{coup2}
\end{equation}
Substituting (\ref{rmatrentry})  into (\ref{coup2}),
 we come to equality
\begin{equation}
f_{\lambda_p}\cdot \Big(
\sum\limits_{1\leq q \leq n}
\frac{g_{\lambda_p}f_{\mu_q}}{\lambda_p-\mu_q}\Big)
\cdot g_{\mu_q}=
-f_{\lambda_p}\cdot g_{\lambda_p}, 
\label{coup3}
\end{equation}
which is the same that the equality
\begin{equation}
f_{\lambda_p}g_{\lambda_p}
\cdot \Big(
\sum\limits_{1\leq q \leq n}
\frac{f_{\mu_q} g_{\mu_q}}{\lambda_p-\mu_q}\Big)   =
-f_{\lambda_p} g_{\lambda_p}\cdot
\label{coup4}
\end{equation}
Since
 \(f_{\lambda_p}g_{\lambda_p}= R_{\lambda_p},\,\,
f_{\mu_q}g_{\mu_q}= R_{\mu_q} \), the last equality takes the form
\begin{equation}
R_{\lambda_p}\cdot\bigg(I+
\sum\limits_{1\leq q \leq n}\frac{R_{\mu_q}}{\lambda_p-\mu_q}\bigg)=0.
\label{coup5}
\end{equation}
According to (\ref{addnonfact_b}),
\begin{displaymath}
I+\sum\limits_{1\leq q \leq n}\frac{R_{\mu_q}}{\lambda_p-\mu_q}=
R^{-1}({\lambda_p}).
\end{displaymath}
Hence, (\ref{coup5}) takes the form
  \(R_{\lambda_p}\cdot R^{-1}({\lambda_p})=0.\) According to (\ref{cons2a}), 
this equality is true.
\hfill\framebox[0.45em]{ }\\[0.0em]

 The proofs of theorems \ref{rightreprtheo} and  \ref{leftreprtheo} 
are analogous. We will prove only the first one of them.\\[10pt]
PROOF of Theorem \ref{rightreprtheo} .\ \
I. First of all, we obtain the representation (\ref{rsysrepr}).
The main idea of the proof is to expand the function \(K^r_R\) of two
variables into a double simple fraction series and then to interpret
this expansion as the matricial equality (\ref{rsysrepr}). Actually, we
derive not a {\it double} expansion of a function of two variables,
but an {\it iterated} one.  Let us fix a point
 \(y \in \overline{{\Bbb C}} \setminus {\cal N}(R)\).
For this fixed value \(y\), we consider the kernel
 \(K^r_{\scriptscriptstyle R}(x,y)\)
 as a function of the variable \(x\).
This function is rational with respect to \(x\), and
\(K^r_{\scriptscriptstyle R}(x,y)\rightarrow 0\)
 by \(x\rightarrow \infty\). It may have
 singularities only at the points of the set \({\cal P}(R)\) and at the point
\(y\) where denominator \(x-y\) vanishes. Actually this function is holomorphic
 at the point \(x=y\) because the numerator vanishes at the point \(x=y\) as
 well.
At each of the points \(\lambda_p \in {\cal P} (R)\) the function
 \(K^r_{\scriptscriptstyle R}(x,y)\)
(considered as a function of \(x\)) either is holomorphic or has a simple pole
with  residue \(K^r_{\lambda_p }(y)\)
of rank one:
\begin{equation}
K^r_{\lambda_p }(y)=\frac{R_{\lambda_p}\cdot R^{-1}(y)}{\lambda_p-y},
\quad \quad (1\leq p\leq n ).
\label{partres}
\end{equation}
Expanding the function \(K^r_{\scriptscriptstyle R}(x,y)\) into the simple
 fraction sum,
 we obtain
\begin{equation}
K^r_{\scriptscriptstyle R}(x,y)
=\sum\limits_{1\leq p \leq n}K^r_{\lambda_p }(y)\cdot
 \frac{1}{x-\lambda_p}.
\label{partexp}
\end{equation}
In its turn the residue \(K^r_{\lambda_p }(y)\), considered as a function of
 \(y\), is a rational function. It vanishes at the point \(\infty\).
 This function is also holomorphic at the point 
\(\lambda_p \) because the numerator 
\mbox{\(R_{\lambda_p}\cdot R^{-1}(y)\)} vanishes at the point \(\lambda_p\):
 the equality \(R_{\lambda_p}\cdot R^{-1}(\lambda_p)=0\) is the same as
the equality \(R_{\lambda}B_0=0\) in (\ref{cons2a}).
 Thus, the only possible singularities of the
function \(K^r_{\lambda_p }(y)\) are the points \(\mu_q\,\)
\mbox{\( (q=1,\,\dots\,, n ) \)}
 of the set \({\cal N} (R)\).
These singularities are simple poles, with the residues \(K^r_{p,q}\):
\begin{equation}
K^r_{p,q}=\frac{R_{\lambda_p}\cdot R_{\mu_q}}{\lambda_p-\mu_q},
\quad (1\leq p,q \leq n ).
\label{doubleres}
\end{equation}
Expanding the function \(K^r_{\lambda_p }\) into the simple fraction sum, we
 obtain
\begin{equation}
K^r_{\lambda_p }(y)=\sum\limits_{1\leq q \leq n}K^r_{p,q}\cdot
\frac{1}{y-\mu_q}.
\label{simplfracexp}
\end{equation}
Combining (\ref{partexp}) and (\ref{simplfracexp}) (and transforming
the iterated sum into the double sum), we obtain
the double expansion
\begin{equation}
K^r_{\scriptscriptstyle R}(x,y)
=\sum\limits_{1\leq p \leq n\atop 1\leq q \leq n}
\frac{1}{x-\lambda_p}\cdot K^r_{p,q} \cdot \frac{1}{y-\mu_q}.
\label{doubleexp}
\end{equation}
Substituting into (\ref{doubleres}) expressions (\ref{factres}), we get
\begin{equation}
K^r_{p,q}=f_{\lambda_p}\cdot h^r_{p,q}\cdot g_{\mu_q},
\label{doubleresoth}
\end{equation}
where \(h_{p,q}^r\) are defined by (\ref{rmatrentry}). Thus, the expansion
(\ref{doubleexp}) takes the form
\begin{equation}
 K^r_{\scriptscriptstyle R}(x,y)
=\sum\limits_{1\leq p \leq n\atop 1\leq q \leq n}
\frac{f_{\lambda_p}}{x-\lambda_p}\cdot h^r_{p,q}
 \cdot \frac{g_{\mu_q}}{y-\mu_q},
\label{modiddoubleexp}
\end{equation}
where \(f_{\lambda_p}\) and \(g_{\mu_q}\) are the left semi-residual vector
at the pole \(\lambda_p\) and the right semi-residual vector at the ``zero''
\(\mu_q\) respectively.

The representation (\ref{rsysrepr}) is simply the
 representation 
(\ref{modiddoubleexp})
in the matricial form.
{\sf The statement 1 of  Theorem \ref{rightreprtheo} is proved.}

From (\ref{doubleexp}) it follows, that
\begin{equation}
K^r_{p,q}=\lim_{x\rightarrow \lambda_p\atop y\rightarrow \mu_p}
(x-\lambda_p)(y-\mu_q)\cdot
 K^r_{\scriptscriptstyle R}(x,y).
\label{uniq}
\end{equation}
Thus, the values \(K_{p,q}\) are determined from  the kernel
 \(K^r_{\scriptscriptstyle R}(x,y)\)
uniquely. From (\ref{doubleresoth}) it follows, that (for given 
\(f_{\lambda_p},\,g_{\mu_q}\)) the values \(h^r_{p,q}\) are determined
uniquely. {\sf The statement 2 of Theorem \ref{rightreprtheo} is proved.}

The statement 3 of Theorem \ref{rightreprtheo} follows immediately from 
 Theorem \ref{invtheo} below,
 where we not only prove the invertibility of matrices \(H_r\) and
 \(H_r\), but also find their inverse matrices.
\hfill\framebox[0.45em]{ }\\[0.0em]     

\vspace*{0.0em}
\begin{theo}   
 \label{invtheo}
{\sl
Let \(R\) be a rational matrix function in general position, and
the matrices \(H^r,\,\, H^l \) are defined from it  according to
(\ref{rmatrentry}), (\ref{lmatrentry}), where
 \( \{\lambda_1,\, \dots \, \lambda_n\}={\cal P}(R)\), 
\(\{\mu_1,\, \dots \, \mu_n\}={\cal N}(R)\} \),
\( f_{\lambda},\,  g_{\lambda} ,\, f_{\mu},\,  g_{\mu} \) are appropriate
 semi-residual vectors.}
Then the equality holds:
\begin{equation}
H^r \cdot H^l = I, \qquad H^l \cdot H^r = I
\label{invrel}
\end{equation}
\end{theo}
\begin{rem}
If we already know from somewhere that the matrix \(H^l\) is invertible,
then we can easily deduce that \(\big(H^{l}\big)^{-1}=H^r\).
Indeed, multiplying the equation (\ref{LS2}) by the matrices 
\(\big(H^l\big)^{-1}\) from the both sides, from the right and from the left,
and taking into account that
\[\big(H^l\big)^{-1}G_{\cal N}=-G_{\cal P}, \qquad
F_{\cal P}\big(H^l\big)^{-1}=F_{\cal N }\] (these are equalities
(\ref{coup}.b) and (\ref{coup}.d)) , we obtain the equality
\begin{equation}
A_{\cal P}\big(H^l\big)^{-1}-\big(H^l\big)^{-1}A_{\cal N}=
G_{\cal P}F_{\cal N}.
\label{SL5}
\end{equation}
Thus, each of the matrices \( \big(H^l\big)^{-1} \) and \(H^r\) is the
 solution  of the same Sylvester-Lyapunov equation
\( A_{\cal P}X-XA_{\cal N}=G_{\cal P}F_{\cal N}.\)
The condition \( {\cal P}(R)\cap{\cal N}(R)=\emptyset \) means that
\({\sigma}_{A_{\cal P}}\cap {\sigma}_{A_{\cal N}}=\emptyset\).
Under this condition, the solution \(X\) of the Sylvester-Lyapunov equation
\(A_{\cal P}X-XA_{\cal N}=G_{\cal P}F_{\cal N}\)
 is unique. Hence, \( \big(H^l\big)^{-1}=H^r\).
\hfill\framebox[0.45em]{ }\\[0.0em]
\end{rem}

We give two proofs of Theorem \ref{invtheo}.\\[4pt]
PROOF I. We prove only the first equality in (\ref{invrel}).
Multiplying (\ref{LS1}) by the matrix \(H^l\) from the right and
(\ref{LS2}) by the matrix \(H^r\) from the left, we came to the equalities
\[
A_{\cal P}H^rH^l-H^rA_{\cal N}H^l=G_{\cal P}F_{\cal N}H^l
\]
and
\[
H^rA_{\cal N}H^l-H^rH^lA_{\cal P}=H^rG_{\cal N}F_{\cal P}.
\]
Taking into account the equalities (\ref{coup}.a) and  (\ref{coup}.d),
we obtain
\[
A_{\cal P}H^rH^l-H^rA_{\cal N}H^l=G_{\cal P}F_{\cal P}
\]
and
\[
H^rA_{\cal N}H^l-H^rH^lA_{\cal P}=-G_{\cal P}F_{\cal P}.
\]
Adding two last equalities, we see that  the matrices
 \(H^rH^l\) and \(A_{\cal P}\) commute:
\[
(H^rH^l)\,A_{\cal P}=A_{\cal P}\, H^rH^l.
\]
Hence,
\begin{equation}
(H^rH^l)\,\varphi(A_{\cal P})=\varphi(A_{\cal P})\,(H^rH^l).
\label{com}
\end{equation}
for every function \(\varphi\) which is holomorphic on the spectrum
 \(\sigma_{A_{\cal P}}\) of the matrix \(A_{\cal P}\).
>From (\ref{coup}.a) and (\ref{coup}.b) it follows that
\[(H^rH^l)\,G_{\cal P}=G_{\cal P}.\]
Multiplying this equality by \(\varphi(A_{\cal P}) \) from the left
and taking into account the commutational relation (\ref{com}), we obtain
that
\begin{equation}
(H^rH^l-I)\,\varphi(A_{\cal P})\,G_{\cal P}=0.
\label{min}
\end{equation}
Let us fix an index \(q \in [1,\,\dots \, ,\,n]\) and specify the function
\(\varphi :\,\, \varphi(\lambda_p)=\delta_{p,q},\,p=1,\,2,\,\dots ,\,n.\)
By such choice of \(\varphi\),
 \(\varphi_(A_{\cal P})=
\mbox{\rm diag}[\,\delta_{1,q},\,\delta_{2,q},\,\dots,\,\delta_{n,q}\,].\)
 Thus,
\[
\varphi(A_{\cal P})\,G_{\cal P}=
 \left[
\begin{array}{c}
\delta_{1,q}\, g_{{}_{\lambda_1}}
\\
 \vdots 
\\
\delta_{n,q}\, g_{{}_{\lambda_n}}
\end{array}\right].
\]
Therefore we obtain that \(m_{p,q}\,g_{{}_{\lambda_q}}=0\)\ \  for 
 \(1\leq p,q\leq n\),
where \(M\stackrel{\rm\tiny def}{=}H^rH^l-I,\quad
M=\|m_{p,q}\|_{1\leq p,q\leq n}.\) Since \(m_{p,q}\in {\Bbb C}\) and
\(g_{{}_{\lambda_q}}\) is a non-zero vector row, \( m_{p,q}=0\) for all
\(1\leq p,q\leq n\), i.e. \(M=0\). Hence,
\(H^rH^l-I=0.\)
\hfill\framebox[0.45em]{ }\\[4pt]
PROOF II. We prove only the first equality in (\ref{invrel}). This matrix
equality is equivalent to the system of \(n^2\) scalar equations
\begin{equation}
\sum\limits_{1\leq q \leq n} h_{p,q}^r \cdot h_{q,p}^l = 1, \qquad
p=1,\,2,\, \dots \,,\, n.
\label{inveq1}
\end{equation}
and
\begin{equation}
\sum\limits_{1\leq q \leq n} h_{p,q}^r \cdot h_{q,s}^l = 0, \qquad
p=1,\,2,\, \dots \,,\, n; \,\, s=1,\,2,\, \dots \,,\, n; \,\, p\not= s.
\label{inveq2}
\end{equation}
According to (\ref{rmatrentry}), (\ref{lmatrentry})
 equality (\ref{inveq1}) means that
\begin{equation}
\sum\limits_{1\leq q \leq n}
\frac{(g_{\lambda_p}f_{\mu_q})\cdot (g_{\mu_q}f_{\lambda_p})}
{(\lambda_p - \mu_q)^2}=-1.
\label{checkpp}
\end{equation}
Because \(f_{\lambda_p}\not=0,\,\, g_{\lambda_p}\not=0\),
the last equality is equivalent\,\footnote{\label{footadd}
\ If \(f\) is non-zero vector-columns, \(g\) is non-zero vector-row,
 then the equality \(c_1=c_2\), where \(c_1,\, c_2\)
are complex numbers, is equivalent 
 to the equality \(fc_1g=fc_2g\).
}
\  to the equality
\begin{equation}
f_{\lambda_p}\cdot
 \sum\limits_{1\leq q \leq n}
\frac{(g_{\lambda_p}f_{\mu_q})\cdot (g_{\mu_q}f_{\lambda_p})}
{(\lambda_p - \mu_q)^2}\cdot g_{\lambda_p}=
 - f_{\lambda_p}\cdot g_{\lambda_p}.
\end{equation}
or, what is the same\,\footnote{\label{footass}
\, Here we use the associativity of the matrix multiplication.}
, to the equality
\begin{equation}
f_{\lambda_p} g_{\lambda_p}\cdot
\sum\limits_{1\leq q \leq n}
\frac{f_{\mu_q} g_{\mu_q}}
{(\lambda_p - \mu_q)^2}
\cdot f_{\lambda_p}g_{\lambda_p}=
- f_{\lambda_p} g_{\lambda_p}.
\label{checkpp1}
\end{equation}
Taking into account the factorization (\ref{factres}), we see, that
the equality (\ref{checkpp1}) is equivalent to the equality
\begin{equation}
R_{\lambda_p} \cdot  \sum\limits_{1\leq q \leq n}
\frac{R_{\mu_q}}
{(\lambda_p - \mu_q)^2}\cdot R_{\lambda_p} =
 - R_{\lambda_p}.
\label{checkpp2}
\end{equation}
In view of (\ref{addnonfact_b}), 
\begin{equation}
\sum\limits_{1\leq q \leq n} \frac{R_{\mu_q}}{(\lambda_p - \mu_q )^2}=
-\big(R^{-1}\big)^{\prime}({\lambda_p}).
\label{verivpp}
\end{equation}
Thus, the equality (\ref{checkpp2}) takes the form
\begin{equation}
R_{\lambda_p}\cdot \big(R^{-1}\big)^{\prime}(\lambda_p)\cdot R_{\lambda_p}=
R_{\lambda_p}.
\label{finalverif}
\end{equation}
According to Remark \ref{usid}, the  equality (\ref{finalverif}) holds.
(See (\ref{rel1})).
Thus, the equalities (\ref{inveq1}) are established.

The equalities (\ref{inveq2}) can be established in the same way.
According to (\ref{rmatrentry}) and (\ref{lmatrentry}),
the  equality (\ref{inveq2}) means that
\begin{equation}
\sum\limits_{1 \leq q \leq n}
\frac{(g_{\lambda_p}f_{\mu_q})\cdot (g_{\mu_q}f_{\lambda_s})}
{(\lambda_p - \mu_q)\cdot (\lambda_s - \mu_q)}=0\, , \qquad p\not= s.
\label{checkpsint.tex}
\end{equation}
Because \(f_{\lambda_p}\not=0,\,\, g_{\lambda_s}\not=0\),
 the last equality is equivalent
(see the footnote \({}^{\ref{footadd}}\)) 
 to the equality
\begin{equation}
f_{\lambda_p}\cdot
\sum\limits_{1 \leq q \leq n}
\frac{(g_{\lambda_p}f_{\mu_q})\cdot (g_{\mu_q}f_{\lambda_s})}
{(\lambda_p - \mu_q)\cdot (\lambda_s - \mu_q)}
\cdot g_{\lambda_s}
=0,
\end{equation}
or (see the footnote \({}^{\ref{footass}}\))
, what is the same, to the equality
\begin{equation}
f_{\lambda_p}g_{\lambda_p}\cdot
\sum\limits_{1 \leq q \leq n} \bigg(
\frac{f_{\mu_q} g_{\mu_q}}{\lambda_p - \mu_q}
-\frac{f_{\mu_q} g_{\mu_q}}{\lambda_s - \mu_q}\bigg)
\cdot f_{\lambda_s}g_{\lambda_s}
=0, \qquad   p\not=s, 
\label{checkps1}
\end{equation}
Taking into account the factorization (\ref{factres}), we see, that
the equality (\ref{checkps1}) is equivalent to the equality
\begin{equation}
R_{\lambda_p}\cdot
\sum\limits_{1 \leq q \leq n} \bigg(
\frac{R_{\mu_q}}{\lambda_p - \mu_q}
-\frac{R_{\mu_q}}{\lambda_s - \mu_q}\bigg)
\cdot
R_{\lambda_s}=0.
\label{checkps}
\end{equation}
In view of (\ref{addnonfact_b}), 
\begin{equation}
\sum\limits_{1 \leq q \leq n}
\bigg( \frac{R_{\mu_q}}{\lambda_p - \mu_q}
-\frac{R_{\mu_q}}{\lambda_s - \mu_q} \bigg)
= R^{-1}(\lambda_p) - R^{-1}(\lambda_s), \qquad p\not= s.
\label{verivps}
\end{equation}
Thus, the equality (\ref{checkps}) takes the form
\begin{equation}
R_{\lambda_p}\cdot
\Big(R^{-1}(\lambda_p)- R^{-1}(\lambda_s)\Big)\cdot R_{\lambda_s}=0.
\label{finalverifs}
\end{equation}
In view of (\ref{cons2a}) and  (\ref{cons1a}), 
\begin{equation}
R_{\lambda_p}\cdot R^{-1}(\lambda_p)=0, \qquad
 R^{-1}(\lambda_s)\cdot R_{\lambda_s}=0.
\end{equation}
Thus, (\ref{checkps}) holds. The equalities (\ref{inveq2}) are established.
\hfill\framebox[0.45em]{ }\\[0.0em]

The representations (\ref{rpsjoint}) and (\ref{lpsjoint}) are almost
what we need.
However, there is an essential disadvantage in these representations:
each one of them contains explicitly  {\sf all the four} semi-residual
matrices. For example, the representation (\ref{rpsjoint}) contains 
explicitly not
only the semi-residual matrices \(F_{\cal P}\) and \(G_{\cal N}\)
(this is evident), but also  the matrices 
 \(F_{\cal N}\) and  \(G_{\cal P}\) (see the expression (\ref{rmatrentry})
for the  the right core matrix \(H^r\)). But the {\sf four}
 semi-residual matrices
(together with zero and pole locations) are {\sf over-determined data}:
 the matrix function \(R\) is completely determined by two of those
semi-residual  matrices only.
 For example, from the additive representation (\ref{matradddir}) it follows,
that (under the normalization \(R(\infty)=I\)) the zero and pole matrices 
 \(A_{\cal P},\,A_{\cal N}\) together with the two semi-residual matrices
\(F_{\cal P},\,G_{\cal P}\) determine completely
the matrix function \(R\), and hence the other two semi-residual matrices
\(F_{\cal N},\,G_{\cal N}\). Because of this, the semi-residual matrices
\(F_{\cal N},\,G_{\cal P}\) are (at least in principle) expressible in terms 
of the matrices \(A_{\cal P},\,A_{\cal N},\,F_{\cal P},\,G_{\cal N}\).
Hence, we can hope to express the right core matrix \(H^r\) in
terms of the matrices \(A_{\cal P},\,A_{\cal N},\,F_{\cal P},\,G_{\cal N}\).
Indeed, this can be done easily and explicitly: on the one hand, the
left core  matrix \(H^l\) is expressible in terms of the entries
of the matrices \(A_{\cal P},\,A_{\cal N},\,F_{\cal P},\,G_{\cal N}\) only
(see (\ref{lmatrentry}); on the other hand, \(H^r=(H^l)^{-1}\)
(see (\ref{invrel})).


This suggests us that it may be reasonable to use the inverse matrices
\(S^r \) and \(S^l \) instead of the matrices\footnote{
We remind, that, according to Theorem 2.3, the matrices \(H^r\) and
\(H^l\)  are mutually inverse, and hence, invertible.} \(H^r\) and \(H^l\):
\begin{equation} 
S^r\stackrel{\rm\tiny def}{=}\big(H^r\big)^{-1}
\quad  \mbox{\rm and} \quad
S^l\stackrel{\rm\tiny def}{=} \big(H^l\big)^{-1}.
\label{coupmatr}
\end{equation}

The equalities (\ref{coup}) and (\ref{LS1}),\,  (\ref{LS2}) can be rewritten
in term of the matrices \(S^r\) and \(S^l\).
Namely, equalities (\ref{coup}) take the form
\begin{equation}
\mbox{\rm a)}.\,\,     G_{\cal N}=-S^rG_{\cal P}; \quad
\mbox{\rm b)}.\,\,     G_{\cal P}=-S^lG_{\cal N}; \quad
\mbox{\rm c)}.\,\,     F_{\cal P}= F_{\cal N}S^r; \quad
\mbox{\rm d)}.\,\,     F_{\cal N}= F_{\cal P}S^l.
\label{recoup}
\end{equation}

Multiplying the equality (\ref{LS1}) by the matrix \((H^r)^{-1}\) from the
right and from the left and taking into account equalities (\ref{coup}.a)
and (\ref{coup}.c), we transform (\ref{LS1}) to the form
\begin{equation}
A_{\cal N}S^r - S^rA_{\cal P}=G_{\cal N}F_{\cal P}.
\label{Lyap1}
\end{equation}
Analogously, from (\ref{LS2}) and (\ref{coup}.b), (\ref{coup}.d)
we derive the equality
\begin{equation}
A_{\cal P}S^l - S^lA_{\cal N}=G_{\cal P}F_{\cal N}.
\label{Lyap2}
\end{equation}
Thus, the matrices \(S^r\) and \(S^l\) are solutions of the
Sylvester-Lyapunov equations
\begin{equation}
A_{\cal N}X-XA_{\cal P}=G_{\cal N}F_{\cal P}.
\label{SL4}
\end{equation}
\begin{equation}
A_{\cal P}X-XA_{\cal N}=G_{\cal P}F_{\cal N}.
\label{SL3}
\end{equation}
respectively.

Now we {\sf change} our point of view and {\sf define} the matrices
\(S^r\) and \(S^l\) as solutions of Sylvester-Lyapunov equations (but not as
the  matrices inverse to the core matrices \(H^r\) and  \(H^s\); see
(\ref{coupmatr})).

\vspace{1.0ex}
\noindent
\begin{defi}
{\rm
Let \(R\) be a rational matrix function in general position,
\(A_{\cal P}\) and \(A_{\cal N}\) be its pole and zero matrices,
\(F_{\cal P},\,G_{\cal P},\, F_{\cal N},\,G_{\cal N}\) be its appropriate
semi-residual matrices.
\begin{enumerate}
\setlength{\itemsep}{-0.5cm}

\vspace{-1.8ex}
\item
{\sl
The matrices \(S^r\) and \(S^l\) which are the solutions of the
Sylvester-Lyapunov equations (\ref{SL4}) and (\ref{SL3}), are said to be
{\sf the right zero-pole coupling matrix}
 and
 {\sf the left zero-pole coupling matrix} respectively.}\\
\item
{\sl
The relations {\rm (\ref{recoup})} are said to be
{\sf the zero-pole coupling relations}}.
\end{enumerate}
}
\end{defi}
\begin{rem}
Since the spectra of the matrices \(A_{\cal P}\) and \(A_{\cal N}\) do not
intersect, the Sylvester-Lyapunov equations (\ref{SL3}) and (\ref{SL4})
are uniquely solvable. (However, as the matrices 
 \(A_{\cal P}\) and \(A_{\cal N}\)
are diagonal, the solvability of these equations as well as the uniqueness is
obvious). Moreover, it is possible to obtain
the explicit expressions for the matrices \(S^r\) and 
\(S^r\) from (\ref{SL3}) and (\ref{SL4}):
\begin{equation}
{\rm r).\ }
S^r=\|s^r_{p,q}\|_{1\leq p,q\leq n}, \quad
s^r_{p,q}=\frac{g_{\mu _p}\cdot f_{\lambda _q}}{\mu _p -\lambda _q},
\qquad
{\rm l).\ }
S^l=\|s^l_{p,q}\|_{1\leq p,q\leq n}, \quad
s^l_{p,q}=\frac{g_{\lambda _p}\cdot f_{\mu _q}}{\lambda _p - \mu _q }.
\label{S^}
\end{equation}
(Actually, we {\it derived} the Sylvester-Lyapunov equations from
the explicit expressions for the matrices which we interpret now
as  solutions of these equations).
\hfill\framebox[0.45em]{ }\\[0.0em]
\end{rem}

\vspace{-0.2em}
According to (\ref{coupmatr}), the relations (\ref{invrel}) 
 can be rewritten in the form
\begin{equation}
S^r\cdot S^l = S^l\cdot S^r=I.
\label{invs}
\end{equation}
We may also  refer directly to the equalities
 (\ref{checkpp}) and (\ref{checkpsint.tex}):
 these equalities mean that the matrices (\ref{S^}.r) and (\ref{S^}.l)
are mutually inverse.

\vspace{0.7em}
\noindent
\begin{rem}
The representations (\ref{rsysrepr})  and (\ref{lsysrepr}) may be rewritten
 in terms of the 
matrices \(S^r\) and \(S^l\) (instead of the matrices \(H^r\) and \(H^l\)):
\begin{eqnarray}
 (R(x)(R(y))^{-1}&
=&I+(x-y)\,
F_{\cal P}\,\big(xI- A_{\cal P}\big)^{-1}\cdot
(S^r)^{-1}\cdot\big(yI- A_{\cal N}\big)^{-1}G_{\cal N},
 \label{rpsjointS}\\ 
 (R(x))^{-1}R(y)&=&I+(x-y\,)
F_{\cal N}\,\big(xI- A_{\cal N}\big)^{-1}\cdot
(S^l)^{-1}\cdot\big(yI- A_{\cal P}\big)^{-1}G_{\cal P}.
  \label{lpsjointS}
\end{eqnarray}
Under the normalizing condition (\ref{norm}),
the matrix functions \(R\),\,\(R^{-1}\) themselves admits the representations
\begin{eqnarray}
R(z)&=&I-F_{\cal P}\,(zI-A_{\cal P})^{-1}(S^{r})^{-1}G_{\cal N},
\label{forRs}
\\
R^{-1}(z)&=&I+F_{\cal P}\,(S^{r})^{-1}(zI-A_{\cal N})^{-1}G_{\cal N},
\label{forR^{-1}s}\\
R(z)&=&I+F_{\cal N} (S^{l})^{-1}  (zI-A_{\cal P})^{-1}G_{\cal P},
\label{forRR^{-1}s1}\\
R^{-1}(z)&=&I-F_{\cal N}(zI-A_{\cal N})^{-1}(S^{l})^{-1}G_{\cal P}.
\label{forRR^{-1}2s}
\end{eqnarray} 
These formulas may be obtained of from (\ref{rpsjointS}),\,
(\ref{lpsjointS}), letting \(x\) or \(y\) tend to \(\infty \) there,
or from (\ref{forR}) - (\ref{forRR^{-1}2}), rewriting them
 in terms of the matrices 
\(S^{r}\),\,\(S^{l}\) (instead of the matrices \(H^{r}\),\,\(H^{l}\)).
\noindent
\mbox{} \hfill\framebox[0.45em]{ }
 \end{rem}

\begin{rem}
Of course, we may obtain the zero-pole coupling relations
(\ref{recoup}) comparing the representations (\ref{matradddir}),\,
(\ref{matraddinv}) and 
(\ref{forR}) - (\ref{forRR^{-1}2}).
\noindent
\mbox{} \hfill\framebox[0.45em]{ }
 \end{rem}

\vspace{0.3em}
\noindent
\begin{rem}
The semi-residual vectors are defined not completely uniquely, but up to
transformations (\ref{polegauge}) and (\ref{zerogauge}) only.
If the semi-residual vectors are transformed according to
(\ref{polegauge}) and (\ref{zerogauge}), the right hand sides 
of the Sylvester-Lyapunov equations
 (\ref{SL4}) and (\ref{SL3}) are transformed as:
\begin{equation}
G_{\cal N}F_{\cal P}\rightarrow 
 (D_{\cal N})^{-1}\cdot G_{\cal N}F_{\cal P}\cdot D_{\cal P}; \qquad
G_{\cal P}F_{\cal N}\rightarrow 
 (D_{\cal P})^{-1}\cdot G_{\cal P}F_{\cal N}\cdot D_{\cal N}
\label{eqtransf}
\end{equation}
The solutions \(S^r\) and \(S^l\) of the Sylvester-Lyapunov equations
(\ref{SL4}) and (\ref{SL3}) are transformed as:
\begin{equation}
S^r \rightarrow (D_{\cal N})^{-1}\cdot S^r\cdot D_{\cal P}; \qquad
S^l \rightarrow (D_{\cal P})^{-1}\cdot S^l\cdot D_{\cal N}.
\label{transcoup}
\end{equation}
Of course, the expressions (\ref{rpsjointS}), (\ref{lpsjointS})
 (for \(R(x)(R(y))^{-1}\) and \((R(x))^{-1}R(y)\)) are invariant with respect
to the transformations
 (\ref{polegauge}), (\ref{zerogauge}), (\ref{transcoup})
(of the semi-residual and the zero-pole coupling matrices).
\hfill\framebox[0.45em]{ }\\[0.0ex]
\end{rem}

The representations (\ref{rpsjointS}) and (\ref{lpsjointS}) are exactly
what we need. Obtaining them is one of the main goal of this paper.
Therefore we choose a special name for this representation:

\begin{defi}\label{Jsr}
{\sl 
Let \(R\) be a rational matrix function in general position, \(A_{\cal P}\),
\(A_{\cal N}\) be its pole and zero matrices,
\(F_{\cal P}\),\,\(G_{\cal P}\),\,\(F_{\cal N}\),\,\(G_{\cal N}\) be
appropriate semi-residual matrices, \(S^r\), \(S^l\) be the solutions
of the Sylvester-Lyapunov equations (\ref{SL4}) and (\ref{SL3})
respectively. The formulas (\ref{rpsjointS}) and (\ref{lpsjointS})
are said to be 
{\sf the right joint system representation of the pair \(R,\,R^{-1}\)}
 and
{\sf the left joint system representation of the pair \(R,\,R^{-1}\)}
 respectively.

The formulas {\rm (\ref{forRs})}, {\rm (\ref{forR^{-1}s})},
{\rm (\ref{forRR^{-1}s1})}, {\rm (\ref{forRR^{-1}2s})} 
(which can be obtained from {\rm (\ref{rpsjointS})}, {\rm (\ref{lpsjointS})}
 by passage to the limit)
 are said to be  {\sf the right system representation of
the function \(R\),  the right system representation of
the function \(R^{-1}\),  the left system representation of
the function \(R\),  the left system representation of
the function \(R^{-1}\)} respectively.  
}
\end{defi}

\begin{rem}
The terminology is motivated by the so-called {\it system theory} or,
in more detail, by {\it the theory of linear time invariant dynamical system}.
In this theory, all the objects such as the zero and pole matrices,
the semi-residual matrices, the zero-pole coupling matrices are interpreted
from the point of view of dynamical systems. This interpretation does not
 play any role in our considerations. We need the joint system
 representations  as a tool to introduce a convenient coordinates 
in the set of all rational matrix functions (in general position).
\hfill\framebox[0.45em]{ }
\end{rem} 

\begin{rem}
In the realization theory one obtains formulas like
(\ref{forRs}) - (\ref{forRR^{-1}2s})
 for matrix functions \(R\) and \(R^{-1}\) considered
{\sf individually}. In the representations
 (\ref{rpsjointS}),\,(\ref{lpsjointS})
the matrix function \(R,\,R^{-1}\) are considered {\sf jointly}.
This is the reason for using the terminology {\sf joint} system
representation.
\hfill\framebox[0.45em]{ }
\end{rem}

Now we summarize the results of this section and formulate 

\begin{theo}\label{summarTheo}
{\sl
Let \(R\) be a rational matrix function in general position,
 \(A_{\cal P}\) and \(A_{\cal N}\) be its pole and zero matrices,
\(F_{\cal P},\,G_{\cal P}\) be its left and right pole semi-residual matrices,
\(F_{\cal N}, G_{\cal N}\) be its left and right zero semi-residual matrices.
Then:}

\vspace{-2.8ex}
 {\rm 
\begin{enumerate}
\setlength{\itemsep}{-1.5ex}
\item{\sl The matrices \(S^r\), which is a solution of
 the Sylvester-Lyapunov equation 
\(A_{\cal N}X-XA_{\cal P}=G_{\cal N}F_{\cal P}\),
\ \    and \(S^l\), which is a solution  of the Sylvester-Lyapunov equation
\(A_{\cal P}X-XA_{\cal N}=G_{\cal P}F_{\cal N}\):
\ \
$$
{\rm r).\ }\,
S^r=\|s^r_{p,q}\|_{1\leq p,q\leq n}, \,
s^r_{p,q}=\frac{g_{\mu_p}\cdot f_{\lambda_q}}{\mu _p -\lambda _q},
\quad
{\rm l).\ }\,
S^l=\|s^l_{p,q}\|_{1\leq p,q\leq n}, \,
s^l_{p,q}=\frac{g_{\lambda_p}\cdot f_{\mu_q}}{\lambda_p-\mu_q},
\eqno{(\ref{S^})}
$$
 are
mutually inverse, i.e. the equalities
$$
S^r\cdot  S^l =I,\qquad  S^l\cdot  S^r =I
\eqno{(\ref{invs})}
$$
 hold. In particular, the matrices
\(S^r\) and \(S^l\) are invertible.}\\[-0.9ex]
\item {\sl The matrix function \(R(x)(R(y))^{-1}\) can be recovered
from the data \(A_{\cal P},\,A_{\cal N},\,F_{\cal P},\,G_{\cal N}\)  
 by the formula (\ref{rpsjointS})
$$
 (R(x)(R(y))^{-1}
=I+(x-y)\,
F_{\cal P}\,\big(xI- A_{\cal P}\big)^{-1}\cdot
(S^r)^{-1}\cdot\big(yI- A_{\cal N}\big)^{-1}G_{\cal N},
\eqno{(\ref{rpsjointS})}
$$
 i.e.
 the right joint system representation  holds.}\\
\item {\sl  The matrix function \((R(x))^{-1}R(y)\) can be recovered
from the data \(A_{\cal P},\,A_{\cal N},\,F_{\cal N},\,G_{\cal P}\)  
 by the formula (\ref{lpsjointS}),
$$
 (R(x))^{-1}R(y)=I+(x-y\,)
F_{\cal N}\,\big(xI- A_{\cal N}\big)^{-1}\cdot
(S^l)^{-1}\cdot\big(yI- A_{\cal P}\big)^{-1}G_{\cal P}.
\eqno{(\ref{lpsjointS})}
$$
 i.e.
 the left joint system representation  holds.}\\
\item
{\sl Under the normalizing condition (\ref{norm}),
the matrix functions \(R\),\,\(R^{-1}\) themselves admit the representations}
\renewcommand{\theequation}{%
\mbox{\ref{forRs}}}%
\begin{equation}
R(z)=I-F_{\cal P}\,(zI-A_{\cal P})^{-1}(S^{r})^{-1}G_{\cal N},
\end{equation}
\renewcommand{\theequation}{%
\mbox{\ref{forR^{-1}s}}}%
\begin{equation}
R^{-1}(z)=I+F_{\cal P}\,(S^{r})^{-1}(zI-A_{\cal N})^{-1}G_{\cal N},
\end{equation}
\renewcommand{\theequation}{%
\mbox{\ref{forRR^{-1}s1}}}%
\begin{equation}
R(z)=I+F_{\cal N} (S^{l})^{-1}  (zI-A_{\cal P})^{-1}G_{\cal P};
\end{equation}
\renewcommand{\theequation}{%
\mbox{\ref{forRR^{-1}2s}}}%
\begin{equation}
R^{-1}(z)=I-F_{\cal N}(zI-A_{\cal N})^{-1}(S^{l})^{-1}G_{\cal P}.
\end{equation}
\addtocounter{equation}{-4}
 \renewcommand{\theequation}{%
\mbox{\arabic{section}.\arabic{equation}}}%
\item
{\sl The zero-pole coupling relations
 hold:}
\renewcommand{\theequation}{%
\mbox{\ref{recoup}}}%
\begin{equation}
\mbox{\rm a)}.\,\,     G_{\cal N}=-S^rG_{\cal P}; \quad
\mbox{\rm b)}.\,\,     G_{\cal P}=-S^lG_{\cal N}; \quad
\mbox{\rm c)}.\,\,     F_{\cal P}= F_{\cal N}S^r; \quad
\mbox{\rm d)}.\,\,     F_{\cal N}= F_{\cal P}S^l.
\end{equation}
\addtocounter{equation}{-1}
\renewcommand{\theequation}{%
\mbox{\arabic{section}.\arabic{equation}}}%
\item {\sl If the representations
 (\ref{rpsjoint}) and (\ref{lpsjoint}) hold with some
matrices \(H^r\) and \(H^s\), then \newline
\(H^r=(S^r)^{-1}\),
\(H^l=(S^l)^{-1}\) of necessity.
}
\end{enumerate}
}
\end{theo}

\begin{rem}\label{hybrid}
The representation (\ref{rpsjointS}) allows us to recover 
the matrix function
 \(R(x)(R(y))^{-1}\)
from
 {\sf the left pole- and the right zero-} semi-residual matrices 
\(F_{\cal P}\) and \(G_{\cal N}\). 
whereas  the representation (\ref{lpsjointS}) allows us to recover 
the matrix function
 \((R(x))^{-1}R(y)\)
from
 {\sf the right pole- and the left zero-} semi-residual matrices 
\(F_{\cal N}\) and \(G_{\cal P}\).
However, sometimes one have needs for some ``{\sf hybrid}\,'' formulas which
allow to recover 
the  matrix function
 \(R(x)(R(y))^{-1}\) from 
 {\sf the right pole- and the left zero-} semi-residual matrices 
\(F_{\cal N}\) and \(G_{\cal P}\)
and the matrix function
 \((R(x))^{-1}R(y)\)
from 
 {\sf the left pole- and the right zero-} semi-residual matrices 
\(F_{\cal P}\) and \(G_{\cal N}\).
Such formulas can be easily derived from
 the joint system realization formulas (\ref{rpsjointS}),
(\ref{lpsjointS}) combined with
 the  zero-pole coupling relations (\ref{recoup}).
These ``hybrid'' formulas are of the form:
\begin{eqnarray}
R(x)\cdot (R(y))^{-1}
&=&
I-(x-y)\, F_{\cal N}\, (S^{l})^{-1}\, (xI-A_{\cal P})^{-1}\, S^l \,
(yI-A_{\cal N})^{-1}\, (S^{l})^{-1}G_{P},\label{hybr1} \\
(R(x))^{-1}\cdot R(y)
&=&
I-(x-y)\, F_{\cal P}\, (S^{r})^{-1}\, (xI-A_{\cal N})^{-1}\, S^r \,
(yI-A_{\cal P})^{-1}\, (S^{r})^{-1}G_{N}\label{hybr2}.
\end{eqnarray}
 The matrix \(S^l\) can be calculated from the data: 
\(F_{\cal N},\,G_{P},\,A_{\cal P},\,A_{\cal N}  \);
the matrix \(S^r\) can be calculated from the data: 
\(F_{\cal P},\,G_{N},\,A_{\cal P},\,A_{\cal N}  \).
\hfill\framebox[0.45em]{ }\\[0.0em]
\end{rem}

\noindent
\vspace*{0.2cm} 
\begin{minipage}{15.0cm}
\section{\hspace{-0.4cm}.\hspace{0.19cm}  FROM THE CHAIN IDENTITY TO THE
 SYLVESTER - LYAPUNOV EQUATION AND BACK. }
\end{minipage}\\[-0.5cm]
\setcounter{equation}{0}

The consideration of this item are concentrated around of the so-called
{\it chain identity}. Let us give a number of  definitions.
 
Let \(T(\,.\,,\,.\,)\) be a  \(k\times k\) matrix function
 of two complex variables, with domain of definition 
\footnote{We recall that
\(\overline{{\Bbb C}}\stackrel{\scriptscriptstyle\rm def}{=}{\Bbb C}%
\cup\infty\)
is the extended complex plane.
 }
\({\mathcal D}_T\),
\({\mathcal D}_T\in\overline{{\Bbb C}}\times\overline{{\Bbb C}}\),\,\,
\(T: {\mathcal D}_T\rightarrow {\Mgot}_k\).

\begin{defi}\label{chainID}
{\sl
A function \(T(\,.\,,\,.\,)\) of two variables is said to satisfy
{\sf the chain identity} if 
\begin{equation}
T(x,y)\cdot T(y,z)= T(x,z)  
\label{chainid}
\end{equation}
for every \(x,y,z\) for which \((x,y)\in{\mathcal D}_T \) and
 \((y,z)\in{\mathcal D}_T\).
(In particular, \((x,z)\) must  belong to \({\mathcal D}_T\),
if \((x,y)\in{\mathcal D}_T\) and \((y,z)\in{\mathcal D}_T\)). 
}
\end{defi}

\begin{defi}
{\sl A function  \(T(\,.\,,\,.\,)\)  of two variables is said to satisfy
{\sf the diagonal unity identity} if
\begin{equation}
T(x,x)=I   
\label{unid}
\end{equation}
for every point \(x\) for which \( (x,x)\) belongs to the domain of definition 
\({\mathcal D}_T\)
of the function \(T\).
}
\end{defi}

A class of function \(T\) satisfying both the chain identity and the diagonal
unity identity can be constructed in the following way.

\begin{defi}\label{chainDefi}
{\sl 
Let \(\Phi\) and \(\Phi^{-1}\) be a \(k\times k\) matrix functions of
 one variables with domains of definition \({\mathcal D}_{\Phi}\) and
\({\mathcal D}_{{\Phi}^{-1}}\) respectively, 
\({\mathcal D}_{\Phi}\in{\overline{\Bbb C}}\),
\({\mathcal D}_{{\Phi}^{-1}}\in{\overline{\Bbb C}}\).
Let us define the matrix function \(T_{\Phi}\) of two variables by the equality
\begin{equation}
T_{\Phi}(x,y)\stackrel{\rm\tiny def}{=} {\Phi}(x)\cdot {\Phi}^{-1}(y)
\label{SpecChainDef}
\end{equation}
with domain of definition 
\begin{equation}
{\mathcal D}_{T}\stackrel{\rm\tiny def}{=}
 {\mathcal D}_{\Phi}\times {\mathcal D}_{{\Phi}^{-1}},\quad \mbox{\rm i.e.}
\quad
\Big((x,y)\in{\mathcal D}_{T}\Big)\Leftrightarrow
 \Big(( x\in{\mathcal D}_{\Phi})\,\&\,(y\in{\mathcal D}_{{\Phi}^{-1}})\Big).
\label{SpecChainDom}
\end{equation}
The function \(T_{\Phi}\) is said to be
 {\sf the chain function generated by the function \({\Phi}\).}
}
\end{defi}

\begin{lemma}\label{chainlemma}
{\sl
Let \(\Phi\) and \(\Phi^{-1}\) be a \(k\times k\) matrix functions of
 one variables with domains of definition \({\mathcal D}_{\Phi}\) and
\({\mathcal D}_{{\Phi}^{-1}}\) respectively, 
\({\mathcal D}_{\Phi}\in{\overline{\Bbb C}}\),
\({\mathcal D}_{{\Phi}^{-1}}\in{\overline{\Bbb C}}\). Let \(T_{\Phi}\) be
the chain function generated by the function \({\Phi}\).
 If the  matrix functions \(\Phi\) and \(\Phi^{-1}\) are
 {\it mutually inverse}, i.e. the identities
\({\Phi}(x)\cdot {\Phi}^{-1}(x)={\Phi}^{-1}(x)\cdot {\Phi}(x)=I\)
hold for all
 \(x\in {\mathcal D}_{\Phi}\cap{\mathcal D}_{{\Phi}^{-1}} \),
then for the matrix function \(T\) both chain identity and diagonal
unity identity hold.
}
\end{lemma}

\begin{rem}
Of course, Lemma \ref{chainlemma} is reach in content only under condition
\({\mathcal D}_{\Phi} \cap {\mathcal D}_{{\Phi}^{-1}}\not=\emptyset\). 
If this condition fails then the values \(T(x,x)\) and 
\(T(x,y)\cdot T(y,z)\) are defined for the empty set of arguments. 
\end{rem}

\noindent
PROOF of Lemma \ref{chainlemma}. The diagonal unity identity expressed that the functions \(\Phi\) and \({\Phi}^{-1}\) are mutually inverse. The chain identity for the function \(T_{\Phi}\) is the consequence of two facts: 1). The
function \(\Phi\) and \({\Phi}^{-1}\) are mutually inverse;
2). The matricial multiplication is associative.
\hfill\framebox[0.45em]{ }\\[0.0ex]

It turns out that each function \(T\) of two variables satisfying both
the chain identity and the diagonal unity identity is of the form \(T_{\Phi}\)
for some function \({\Phi}\) of one variable.

\begin{theo}\label{ChainGen}
{\sl
Let \(T\) be a \(k\times k\) matrix function of two variables, 
which domain of definition
 \({\mathcal D}_T\in\overline{{\Bbb C}}\times\overline{{\Bbb C}}\)
 is of the form
\({\mathcal D}_T={\mathcal D}_1\times{\mathcal D}_2\), where
\({\mathcal D}_1\in \overline{{\Bbb C}}\) and
\({\mathcal D}_2\in \overline{{\Bbb C}}\), with
 \({\mathcal D}_1\cap{\mathcal D}_2\not=\emptyset .\) If for the function \(T\)
both chain identity and diagonal unity identity are satisfied, then the
 function \(T\) is of the form \(T=T_{\Phi}\),
(i.e. \(T(x,y)={\Phi}(x)\cdot {\Phi}^{-1}(y)\)),
 where \(\Phi\) and
 \({\Phi}^{-1}\) are mutually inverse \(k\times k\) matrix
 functions of one variable, with 
\({\mathcal D}_{\Phi}={\mathcal D}_1\) and 
\({\mathcal D}_{{\Phi}^{-1}}={\mathcal D}_2\).
}
\end{theo}

\noindent
PROOF. Let us fix an arbitrary point \(a\) belonging to the set 
\({\mathcal D}_1\cap{\mathcal D}_2.\) (We will call this point \(a\)
{\sf the distinguished point}.) Let us define now
\begin{equation}
{\mathcal D}_{\Phi}\stackrel{\rm\tiny def}{=}{\mathcal D}_1,\,\,
 {\Phi}(x)\stackrel{\rm\tiny def}=T(x,a);\quad
{\mathcal D}_{\Phi}^{-1}\stackrel{\rm\tiny def}{=}{\mathcal D}_2,\,\,
 {\Phi}^{-1}(y)\stackrel{\rm\tiny def}=T(a,y).
\label{SpecialChainFunc}
\end{equation}
The functions \({\Phi}\) and \({\Phi}^{-1}\) are mutually inverse: this
follows from the chain and diagonal unity identities. The equality
\(T(x,y)={\Phi}(x)\cdot{\Phi}^{-1}(y)\) is the chain identity written down
for the triple of the points \(x,a,y\). In addition, we note that this function
\(\Phi\) satisfies {\it the normalizing condition} \({\Phi}(a)=I\).
\hfill\framebox[0.45em]{ }\\[0.0ex]

\begin{defi}\label{chainDefi1}
{\sl
Let \(R\) be a rational \(k\times k\) matrix function of one variables,
 \(\mbox{det}R\not\equiv 0\),
and \(R^{-1}\) is the inverse (in the commonly accepted sense) matrix
function; the domain of definition \({\mathcal D}_R\) is the set of
 holomorphicity
of the function \(R\); the domain of definition \({\mathcal D}_R^{-1}\) is
the set of holomorphicity of the function \(R_{-1}\). (In other words,
\footnote{We recall that \({\cal P}(R)\) is the pole set of the function \(R\),
  \({\cal N}(R)\)
is the zero set of the function \(R\), i.e. the pole  set of the function
\(R^{-1}\). }
\({\mathcal D}_R=\overline{\Bbb C}\setminus {\cal P}(R),\, 
{\mathcal D}_R^{-1}=\overline{\Bbb C}\setminus {\cal N}(R)\)).
We associate with the function \(R\) two functions of two variables,
\(T^{{}^{\,r}}_{{}_R}(x,y)\) and \(T^{{}^{\,l}}_{{}_R}(x,y)\):
\begin{equation}
\mbox{\rm a).} \,\,
T^{{}^{\,r}}_{{}_R}(x,y)   \stackrel{\mbox{\rm\tiny def}}{=}
R(x)\cdot R^{-1}(y);
\qquad
\mbox{\rm b).} \,\,
T^{{}^{\,l}}_{{}_R}(x,y)   \stackrel{\mbox{\rm\tiny def}}{=}
 R^{-1}(x)\cdot R(y).
\label{chaindef}
\end{equation}
The function \(T^{{}^{\,r}}_{{}_R}\) is said to be
{\sf
 the right chain function generated by \(R\).
}
The function \(T^{{}^{\,l}}_{{}_R}\) is said to be 
{\sf 
the left chain function generated
by \(R\).
} 
}
\end{defi}
\ \\[2.9ex] 
\begin{rem}\label{remT_R}
It is clear that that the function \(T^{{}^{\,r}}_{{}_R}\) is 
the function of the form 
\(T_{\Phi}\) (in the sense of the Definition \ref{chainDefi}) for
\({\Phi} =R\), and the function \(T^{{}^{\,l}}_{{}_R}\) is the function
\(T_{\Phi}\) for \({\Phi} =R^{-1}\). Thus,
 the right chain function generated by \(R\) is the left
chain function generated by \(R^{-1}\):
\begin{equation}
T^{{}^{\,r}}_{{}_R}(x,y)=
T^{{}^{\,l}}_{{{\scriptscriptstyle R}^{-1}}}(x,y).
\label{leftright}
\end{equation}
\end{rem}

>From Lemma \ref{chainlemma} and from Remark \ref{remT_R} it follows   

\begin{cor}\label{Rchain}
{\sl
Let \(R\) be a matrix function, \(\mbox{det}R\not\equiv 0\).
Then both matrix functions 
\(T^{{}^{\,r}}_{{}_R}(x,y)\) and \(T^{{}^{\,l}}_{{}_R}(x,y)\)
(see (\ref{chaindef}))
satisfy the chain identity and the diagonal unity identity:
\begin{eqnarray}
T^{{}^{\,r}}_{{}_R}(x,y)\cdot T^{{}^{\,r}}_{{}_R}(y,z)\equiv
T^{{}^{\,r}}_{{}_R}(x,z); \qquad \qquad T^{{}^{\,r}}_{{}_R}(x,x)\equiv I.
\label{spchainid1}
\\
T^{{}^{\,l}}_{{}_R}(z,y)\cdot T^{{}^{\,l}}_{{}_R}(y,x) \equiv
T^{{}^{\,l}}_{{}_R}(z,x); \qquad \qquad T^{{}^{\,l}}_{{}_R}(x,x)\equiv I.
\label{spchainid2}
\end{eqnarray}
}
\end{cor}\hfill\framebox[0.45em]{ }\\[0.0ex]

>From (\ref{rkern}) and (\ref{lkern}) it is clear that
 the right and left chain functions
 \(T^{{}^{\,r}}_{{}_R}\) and \(T^{{}^{\,l}}_{{{\scriptscriptstyle R}^{-1}}}\) 
can be expressed from the right and left kernels:
\begin{eqnarray}
T^{{}^{\,r}}_{{}_R}(x,y)=I+(x-y)\cdot K^{{}^{\,r}}_{\scriptscriptstyle R}(x,y);
\label{kernchain1}\\
T^{{}^{\,l}}_{{}_R}(x,y)=I+(x-y)\cdot K^{{}^{\,l}}_{\scriptscriptstyle R}(x,y).
\label{kernchain2}
\end{eqnarray}
Substituting in the chain identity (\ref{spchainid1})
 the expression (\ref{kernchain1}) for 
\(T^{{}^{\,r}}_{{}_R}(x,y)\), we obtain an identity for the right kernel
\(K^r_{\scriptscriptstyle R}\):
\begin{equation}
\Big(I+(x-y)K^r_{\scriptscriptstyle R}(x,y)\Big) \cdot 
\Big(I+(y-z)K^r_{\scriptscriptstyle R}(y,z)\Big)   \equiv
I+(x-z)K^r_{\scriptscriptstyle R}(x,z).
\label{kernchain3}
\end{equation}
Removing the parentheses, we obtain the identity
\begin{displaymath}
(x-y)K^r_{\scriptscriptstyle R}(x,y) +
(y-z)K^r_{\scriptscriptstyle R}(y,z)-
(x-z)K^r_{\scriptscriptstyle R}(x,z)
\equiv
-(x-y)(y-z)
K^r_{\scriptscriptstyle R}(x,y)
K^r_{\scriptscriptstyle R}(y,z).
\end{displaymath}
Dividing on \( (x-y)(y-z)\), we come to the identity
\begin{equation}
\frac{K^r_{\scriptscriptstyle R}(x,y)-K^r_{\scriptscriptstyle R}(x,z)}{y-z}
-
\frac{K^r_{\scriptscriptstyle R}(x,z)-K^r_{\scriptscriptstyle R}(y,z)}{x-y}
\equiv
-K(x,y)K(y,z).
\label{kernchain4}
\end{equation}
Assume now that the kernel \(K^r_{\scriptscriptstyle R}\) admits 
the representation (\ref{rsysrepr}), with some \(k \times n\) matrix
\(F_{\cal P}\) \(n \times k\) matrix \(G_{\cal N}\) and \(n \times n\) matrices
\(A_{\cal P}, A_{\cal N}, H^r\). Substituting the expressions (\ref{rsysrepr})
into (\ref{kernchain2}), we come to the identity
\begin{eqnarray}
F_{\cal P} \cdot
( xI-A_{\cal P} )^{-1}\cdot H^r\cdot 
\frac
{(yI-A_{\cal N}) ^{-1} - (zI-A_{\cal N}) ^{-1}} 
{y-z}\cdot G_{\cal N} \quad - \mbox{\hspace{55pt}}\nonumber \\
F_{\cal P}\cdot
\frac
{(xI-A_{\cal P})^{-1} - (yI-A_{\cal P})^{-1}} 
{x-y} 
\cdot H^r \cdot
(zI-A_{\cal N})^{-1} 
\cdot G_{\cal N}  \quad  \equiv \mbox{\hspace{55pt}}
 \label{kernchain5} \\
- \quad F_{\cal P}\cdot 
(xI-A_{\cal P})^{-1}
 \cdot H^r \cdot
 (yI-A_{\cal N})^{-1}\cdot 
 G_{\cal N} F_{\cal P}
 \cdot 
(yI-A_{\cal P})^{-1}
 \cdot H^r \cdot
 (zI-A_{\cal N})^{-1}\cdot  G_{\cal N}.
     \nonumber 
\end{eqnarray}
Using the Hilbert identity for resolvents, we come to the identity
\begin{eqnarray}
F_{\cal P} \cdot
( xI-A_{\cal P} ) ^{-1}\cdot  H^r\cdot
 (yI-A_{\cal N}) ^{-1}\cdot (zI-A_{\cal N}) ^{-1} \cdot G_{\cal N} 
\quad - 
\mbox{\hspace{65pt}}
 \nonumber \\
F_{\cal P} \cdot (xI-A_{\cal P} )^{-1} \cdot (yI-A_{\cal P} )^{-1} \cdot
H^r \cdot (zI-A_{\cal N})^{-1}\cdot  G_{\cal N} \quad \equiv 
\mbox{\hspace{65pt}}
 \label{kernchain6}\\
 F_{\cal P}\cdot 
(xI-A_{\cal P})^{-1}
 \cdot H^r \cdot
 (yI-A_{\cal N})^{-1}\cdot 
 G_{\cal N} F_{\cal P}
 \cdot 
(yI-A_{\cal P})^{-1}
 \cdot H^r \cdot
 (zI-A_{\cal N})^{-1}\cdot  G_{\cal N}.
     \nonumber
\end{eqnarray}
Putting the common factors outside the parentheses, we obtain
\begin{equation}
F_{\cal P} \cdot (xI-A_{\cal P} )^{-1}\cdot M 
\cdot
(zI-A_{\cal N})^{-1}\cdot  G_{\cal N} \quad \equiv 0,
\label{kernchain7}
\end{equation}
where
\begin{equation}
M=H^r \cdot (yI-A_{\cal N}) ^{-1} - 
(yI-A_{\cal P} )^{-1} \cdot H^r 
-H^r\cdot (yI-A_{\cal N})^{-1}\cdot G_{\cal N} F_{\cal P}
 \cdot (yI-A_{\cal P})^{-1}
 \cdot H^r
\label{kernchain8}
\end{equation}

Assume moreover that \(R\) is a rational matrix function in general position.
Let
\(A_{\cal P}\) and \(A_{\cal N}\) be its pole and zero matrices
(i.e. these matrices are of the form (\ref{polematr}), where all the numbers
\(\lambda_1,\,\dots\, 
,\,\lambda_n;\,\mu_1,\,\dots\, ,\,\mu_n\)
are pairwise different), and \(F_{\cal P}\) and \(G_{\cal N}\) be the pole
 left semi-residual and the zero right semi-residual matrices
(in particular, they have the form (\ref{leftsemires}), (\ref{rightsemires})
 where
 \(f_{\lambda_p}\) and \(g_{\mu_q}\) are non-zero k vector-columns).

>From (\ref{kernchain7}) it follows that 
\begin{equation}
F_{\cal P} \cdot \varphi (A_{\cal P})\cdot M \cdot \psi (A_{\cal N}) \cdot
G_{\cal N} \quad \equiv 0,
\label{kernchain9}
\end{equation}
where \(\varphi\) and \(\psi\) are arbitrary functions which are analytic
on the spectra of the  matrices \(A_{\cal P}\) and \(A_{\cal N}\) respectively.

Let us fix two indices,
 \(p\in [1,\, \dots\, , n]\) and \(q\in [1,\, \dots\, , n]\) and specify two
functions \(\varphi\) and \(\psi\):
\begin{displaymath}
\varphi ( {\lambda_j})=\delta_{j,p} \qquad
 \qquad \psi ({\mu_j})=\delta_{j,q} , \qquad (j=1,\,2,\, \cdot\, , \, n).
\end{displaymath}
Then
\begin{equation}
\varphi (A _{\cal P})=\mbox{\rm diag}\,\big(\delta _{1,p},\, \dots \,
 ,\delta_{n,p}\big )
 , \qquad
\psi (A _{\cal N})=\mbox{\rm diag}\,\big(\delta _{1,q},\, \dots \,
 ,\delta_{n,q}\big ).
\label{funct}
\end{equation}
For such functions \(\varphi\) and \(\psi\),
 (\ref{kernchain9}) becomes the form:
\begin{equation}
f_{\lambda_p}\cdot m_{p,q}\cdot g_{\mu_q}=0,
\label{kernchain10}
\end{equation}
where \(m_{p,q}\) is \(p,q\) entry of the matrix \(M\).
Since \(f_{\lambda_p}\) and \(g_{\mu_q}\) are non-zero vector-column and
vector-row,
\(m_{p,q}=0.\) Thus, \(M= 0\), or,
\begin{equation}
H^r \cdot (yI-A_{\cal N}) ^{-1} - 
(yI-A_{\cal P} )^{-1} \cdot H^r 
-H^r\cdot (yI-A_{\cal N})^{-1}\cdot G_{\cal N} F_{\cal P}
 \cdot (yI-A_{\cal P})^{-1}
 \cdot H^r \equiv 0.
\label{kernchain11}
\end{equation}
The ``left'' version of the equation (\ref{kernchain11}) has the form
\begin{equation}
H^l \cdot (yI-A_{\cal P}) ^{-1} - 
(yI-A_{\cal N} )^{-1} \cdot H^l 
-H^l\cdot (yI-A_{\cal P})^{-1}\cdot G_{\cal P} F_{\cal N}
 \cdot (yI-A_{\cal N})^{-1}
 \cdot H^l \equiv 0.
\label{kernchain12}
\end{equation}
To obtain (\ref{kernchain12}), we have to use the chain identity 
(\ref{spchainid2}), the expression (\ref{kernchain2}) for the chain function
 \(T^{l}_{R}\)
 in terms of its associated kernel
 \(K^{l}_{R}\),
 and the representation (\ref{lsysrepr}) of this kernel.
However, we may just replace in (\ref{kernchain11}):
\begin{displaymath}
H^r\rightarrow H^l;\qquad F_{\cal P}\rightarrow F_{\cal N};\qquad
G_{\cal N}\rightarrow G_{\cal P}:\qquad
 A_{\cal N}\rightarrow A_{\cal P};\qquad
A_{\cal P}\rightarrow A_{\cal N}.
\end{displaymath}
Let us examine the Laurent expansion (with respect to \(y\))
 of the function in the left hand side of (\ref{kernchain11}):
\begin{eqnarray}
H^r \cdot (yI-A_{\cal N}) ^{-1} - 
(yI-A_{\cal P} )^{-1} \cdot H^r 
-H^r\cdot (yI-A_{\cal N})^{-1}\cdot G_{\cal N} F_{\cal P}
 \cdot (yI-A_{\cal P})^{-1}
 \cdot H^r=\nonumber
\\
(H^r A_{\cal N}- A_{\cal P}H^r-H^r\,G_{\cal N}\cdot  F_{\cal P}
H^r)y^{-2}+\,\,O(y^{-3}) \qquad (y\rightarrow \infty)
\nonumber
\end{eqnarray}
In view of (\ref{kernchain11}), 
the Sylvester-Lyapunov equality  holds:
\begin{equation}
H^rA_{\cal N}-A_{\cal P}H^r =H^rG_{\cal N}\cdot F_{\cal P}H^r.
\label{LSH1}
\end{equation}
Analogously, from (\ref{kernchain12}) we can derive the equality
\begin{equation}
H^lA_{\cal P}-A_{\cal N}H^l =H^lG_{\cal P}\cdot F_{\cal N}H^l.
\label{LSH2}
\end{equation}
According to Theorem \ref{rightreprtheo}
 (actually, according to Theorem \ref{invtheo}),
the {\it core  matrix} \(H^r\) (see Definition \ref{defcoup}) is invertible.
(The first proof of the Theorem \ref{invtheo} is based on the
Sylvester-Lyapunov equalities (\ref{LS1}),\,(\ref{LS2}),
 but the second one is independent from  them).
Multiplying (\ref{kernchain11}) by the matrix 
\(\big(H^r\big)^{-1}\) from the right and by the matrix
\(\big(H^r\big)^{-1}\) from the left
and denoting, as before (see (\ref{coupmatr})),
\(S^r\stackrel{\rm\tiny def}{=}\big(H^r\big)^{-1}\),
 we come to the equality (\ref{Lyap1}).
>From (\ref{kernchain12}) we can derive the equality (\ref{Lyap2})
 in the same way as we already derived the equality (\ref{Lyap1})
from (\ref{kernchain11}).

 Thus, we obtained the equalities (\ref{Lyap1}) and
 (\ref{Lyap2}) in two different ways.
The first one is based on the explicit expressions
 (\ref{rmatrentry}) and (\ref{lmatrentry})
for the core matrices \(H^r\) and  \(H^s\).  This method uses essentially
the specific character of a rational matrix function
{\it in general position}.  The second method works for  much more  
broad classes of rational matrix functions.
Up to certain point, the method works for 
{\it  arbitrary} rational matrix function.
First of all, this method uses the chain identities (\ref{spchainid1}) and
(\ref{spchainid2}). These identities  are evidently true for arbitrary
 matrix functions \(R\) which are non-degenerate
 (i.e. \(\mbox{\rm det}R \not\equiv 0\)).
 Then we use the representations
(\ref{rsysrepr}), (\ref{lsysrepr}) for the kernels, associated with \(R\).
However, we nowhere use that the pole and zero matrices
\(A_{\cal P}\) and \(A_{\cal N}\) are diagonal,
with disjoint simple spectra. 
Actually, we have obtained the equality (\ref{kernchain9})
 for {\it any}
rational matrix function \(R\) which associated kernel \(K_R^r\) admits 
the representation (\ref{rkern}) with {\it arbitrary}
\(A_{\cal P}\),\,\(A_{\cal N}\),\,\(F_{\cal P}\),\,\(G_{\cal N}\,\,
\mbox{\rm and}\,\, H^r\)
(or, what is the same, for {\it any} rational matrix function \(R\),
such that the function \(R(x)\cdot ( R(y))^{-1}\) admits the representation
(\ref{rpsjoint}) with {\it arbitrary}
\(A_{\cal P}\),\,\(A_{\cal N}\),\,\(F_{\cal P}\),\,\(G_{\cal N}\,\,
\mbox{\rm and}\,\, H^r\)).
Then we have to conclude from (\ref{kernchain9}) (under  the assumption
that this equality holds for arbitrary functions \(\varphi\), \(\psi\)
 which are holomorphic on the spectra of \(A_{\cal P}\) and \(A_{\cal P}\)
 respectively),
that \((\ref{kernchain11})\) holds.

\begin{defi}\label{obstroll}(The row-version):
{\sl
Let \(\Gamma \) be a \(k\times n\) matrix (\(k\) rows, \(n\) columns), and
\(A\) be \(n\times n\) matrix. The pair \( (\Gamma ,A) \) is said to be
 {\sf obstrollable}\footnote{
The word {\bf obstrollable} is a mixture of the words {\bf obs}ervable
and con{\bf trollable.}
}
\  if the linear span of the set of \(n\)-vector-rows
\(\{v\,\Gamma (\lambda I - A )^{-1} \}\),
where \(v\) runs over the space \({\Bbb C}^k\) of all \(k\)-vector-rows
 and \(\lambda\)
runs over \({\Bbb C}\setminus{\sigma}_A \)
 ( \({\sigma}_A\) is the spectrum of \(A\)),
 coincides with the whole
 space \({\Bbb C}^n\) (of all \(n\)-vector-rows). \\
{\sf The equivalent definitions:}\\
{\sf I}. The pair \( (\Gamma,A) \) is said to be {\sf obstrollable}, if
 the linear span of the set of vectors
\(\{v\,\Gamma A^m \}\),
where \(v\) runs over \({\Bbb C}^k\)  of all \(k\) vector-rows  and \(m\) runs
over the set \({\Bbb N}\) of all natural numbers,
 coincides with the whole  space \({\Bbb C}^n\) (of all \(n\)-vector-rows).\\
{\sf II}. The pair \( (\Gamma,A) \) is said to be {\sf obstrollable}, if the
linear span of the set of vectors
\(\{v\,\Gamma\varphi (A)\}\),
 where \(v\) runs over \({\Bbb C}^k\)  of all \(k\) vector-rows  and 
\(\varphi\) runs over the set of all functions holomorphic on \({\sigma}_A\),
coincides with the whole space \({\Bbb C}^n\) (of all \(n\)-vector-rows).
\hfill\framebox[0.45em]{ }\\[0.0em]
}
\end{defi}

\noindent
DEFINITION \ref{obstroll} (The column-version):
{\sl
Let \(\Gamma \) be a \(n\times k\) matrix (\(n\) rows,\,\(k\) columns), and
\(b\) be \(n\times n\) matrix. the pair \( (B,\,\Gamma) \) is said to be
 {\sf obstrollable},
 if the linear span of the set of \(n\)-vector-columns
\(\{ (\mu I - B )^{-1}\Gamma v \}\),
where \(v\) runs over the space \({\Bbb C}^k\) of all \(k\)-vector-columns
 and \(\mu\)
runs over \({\Bbb C}\setminus{\sigma}_B \)
 ( \({\sigma}_B\) is the spectrum of \(B\)),
 coincides with the whole
 space \({\Bbb C}^n\) (of all \(n\)-vector-columns). \\
{\sf The equivalent definitions:}\\
{\sf I}. The pair \( (B,\,\Gamma) \) is said to be {\sf obstrollable}, if
 the linear span of the set of vectors
\(\{\Gamma B^m\,v \}\),
where \(v\) runs over \({\Bbb C}^k\)  of all \(k\) vector-columns  and \(m\)
runs over the set \({\Bbb N}\) of all natural numbers,
 coincides with the whole  space \({\Bbb C}^n\) (of all \(n\)-vector-columns).\\
{\sf II}. The pair \( (B,\,\Gamma) \) is said to be {\sf obstrollable}, if the
linear span of the set of vectors
\(\{\varphi (B)\,\Gamma\,v\}\),
 where \(v\) runs over \({\Bbb C}^k\)  of all \(k\) vector-columns  and 
\(\varphi\) runs over the set of all functions holomorphic on \({\sigma}_B\),
coincides with the whole space \({\Bbb C}^n\) (of all \(n\)-vector-columns).
\hfill\framebox[0.45em]{ }\\[0.0em]
}

\noindent
COMMENT TO TERMINOLOGY: This terminology is motivated by the system theory. 
(In more detail, by the theory of linear time invariant dynamical systems).
 In this theory structures like \(F(\lambda I-A)^{-1}\) and
 \((\mu I-B)^{-1}G\) appear, where 
\(F\) and \(G\) are \(k\times n\) and \(n\times k\) matrices and \(A,B\) are
\(n\times n\) matrix,
and usually \(n\) is much bigger then \(k\) (\(n\gg k)\).
If \(F\) is interpreted as the input operator of the system and \(G\) is 
interpreted as its output operator, then the notions of
controllability of the pair \((F,A)\) and the notion of observability 
of the pair \((B,\,G)\) are introduced. 
If \(G\) is interpreted as the input operator of the system and \(F\) is 
interpreted as its output operator, then the notions of
controllability of the pair \((B,\,G)\) and the notion of observability 
of the pair \((F,A)\) are introduced. 
Structures analogous to the structure 
\(F(\lambda I-A)^{-1},\,\,(\mu I-B)^{-1}G\)
appear in the system representation problems as well.
We need to formulate the property which is controllability if \(\Gamma\) is
interpreted as the input operator, and is observability, if
 \(\Gamma\) is interpreted as the output operator. However, we would not like
to give the preference to one of two: in- or out- interpretations of the
matrix \(\Gamma\). Because this, we choose the ``neutral'' term
{\it obstrollability}.
\hfill\framebox[0.45em]{ }\\[0.0em]

\begin{lemma}
{\sl 
Let \(F,G\) be \(k\times n\) and \(n\times k\) matrices, and \(A,B,M\) are
 \(n\times n\)
matrices. Assume that {\rm
\begin{enumerate}
\setlength{\topsep}{0.1cm}
\setlength{\itemsep}{-0.1cm}
\item{\sl
\(F(\lambda I-A)^{-1}M(\mu I-B)^{-1}G\equiv 0\)
\qquad \big(\(\forall \lambda \in {\Bbb C}\setminus\sigma_{A},\quad
\forall\mu \in {\Bbb C}\setminus\sigma_{B}\)\big).
}
\item{\sl
The pairs \((F,A)\) and \((B,G)\) are obstrollable.
}
\end{enumerate}

\vspace{-2.0ex}
\noindent
}
Then \(M=0\).
}
 \end{lemma}
 
\noindent
PROOF. The proof follows immediately from the definition of obstrollability.
\hfill\framebox[0.45em]{ }\\[0.0em]

Thus, the foregoing reasonings (\,the reasoning of this section and the
 reasoning used for the first proof of Theorem \ref{invtheo}\,)
 prove the following

\begin{theo}\label{robstLS1}
{\sl
{\sf I).}
Let \(R\) be a rational function such that\,\footnote{\label{noassump1}
We  assume neither the matrices
 \(A_{\cal P},\,A_{\cal N}\)
 are diagonal, nor
their spectra are simple, nor their spectra are  non-intersecting.
We also don't assume a'priori that the matrices \(H^r\),\,\(H^l\)
are invertible.
}
 the chain matrix function 
\(R(x)(R(y))^{-1}\) admits the representation of the form
{\rm  (\ref{rpsjoint})}  some \(k\times n\)  matrix \(F_{\cal P}\),\, 
\(n\times k\) matrix \(G_{\cal N}\) and \(n\times n\) matrices
 \(A_{\cal P}\),\,\(A_{\cal N}\) and \(H^r\).\\
\hspace*{5pt}  If the pairs \((F_{\cal P},\,A_{\cal P})\) and
 \((A_{\cal N},\,G_{\cal N})\)
are obstrollable, then the the matrix \(H^r\) satisfy the equality
(\ref{LSH1});

{\sf II).}
Let \(R\) be a rational function such that\,\({}^{\ref{noassump1}}\)
 the chain matrix function 
\((R(x))^{-1}R(y)\) admits the representation of the form
{\rm  (\ref{lpsjoint})}  some \(k\times n\)  matrix \(F_{\cal N}\),\, 
\(n\times k\) matrix \(G_{\cal P}\) and \(n\times n\) matrices
 \(A_{\cal P}\),\,\(A_{\cal N}\) and \(H^l\).\\
\hspace*{5pt} If the pairs \((F_{\cal N},\,A_{\cal N})\) and \((A_{\cal P},
\,G_{\cal P})\)
are obstrollable, then the the matrix \(H^l\) satisfy the equality
(\ref{LSH2});

{\sf III).}
Let \(R\) be a rational function such that the chain  matrix  functions 
\(R(x)(R(y))^{-1}\) and \((R(x))^{-1}R(y)\)
 admit the representations of the form
{\rm(\ref{rpsjoint})} and {\rm(\ref{lpsjoint})} respectively, with some
\(k\times n\) matrices \(F_{\cal P},\,F_{\cal N}\), some \(n\times k\)
 matrices \(G_{\cal N},\,G_{\cal P}\) and some \(n\times n\) matrices
 \(A_{\cal P}\),\,\(A_{\cal N}\),\,\(H^r\),\,\(H^l\).\\
\hspace*{5pt} If all four pairs
 \((F_{\cal P},\,A_{\cal P})\),\,
 \((A_{\cal N},\,G_{\cal N})\),\,
\((F_{\cal N},\,A_{\cal N})\) and
\((A_{\cal P},\,G_{\cal P})\)
are obstrollable, and if moreover the coupling relations (\ref{coup}) hold,
then the matrices \(H^r\),\,\(H^l\) are mutually inverse:
\(H^r\cdot H^l=I\),\,\(H^l\cdot H^r=I\),\, and for their inverse matrices
\(S^r=(H^r)^{-1}\) and \(S^l=(H^l)^{-1}\)
the  equalities 
{\rm (\ref{Lyap1})}, {\rm (\ref{Lyap2})} holds, i.e.
 the matrices \(S^r\) and  \(S^l\) are solutions of the
 Sylvester-Lyapunov equations {\rm  (\ref{SL4})} and {\rm  (\ref{SL3})}.
}
\end{theo}

\vspace{5pt}
\begin{theo}\label{robstLS2}
{\sl
Let \(R\) be a rational function such that the chain function 
\(R(x)(R(y))^{-1}\) admits the representation of the form
{\rm  (\ref{rpsjointS})}
with
\footnote{\label{noassump2}We  assume neither the matrices
 \(A_{\cal P},\,A_{\cal N}\)
 are diagonal, nor
their spectra are simple, not their spectra are  non-intersecting. Of course,
the invertibility of the matrix \(S^r\) imposes implicitly some
restrictions on the matrices \(A_{\cal P},\,A_{\cal N},\,F_{\cal P},\,
G_{\cal N}\).
}
 some \(k\times n\) and \(n\times k\) matrices \(F_{\cal P},\,G_{\cal N}\), 
\(n\times n\) matrices \(A_{\cal P}\) and \(A_{\cal N}\)
 and some {\sf invertible} \(n\times n\) matrix \(S^r\).
If the pairs \((F_{\cal P},\,A_{\cal P})\) and
 \((A_{\cal N},\,G_{\cal N})\) are obstrollable, then the equality 
{\rm (\ref{Lyap1})} holds, i.e. the matrix \(S^r\) is a solution of the
 Sylvester-Lyapunov equation {\rm  (\ref{SL4})}.
}
\end{theo}

The ``left'' version of this theorem holds as well.

\begin{theo}\label{lobstLS}
{\sl
Let \(R\) be a rational function such that the chain function 
\((R(x))^{-1}R(y)\) admits the representation of the form
{\rm  (\ref{lpsjointS})}
with\({}^{\ref{noassump2}}\)
 some \(k\times n\) and \(n\times k\) matrices \(F_{\cal N},\,G_{\cal P}\), 
\(n\times n\) matrices \(A_{\cal P}\) and \(A_{\cal N}\)
 and some {\it invertible} \(n\times n\) matrix \(S^l\).
If the pairs \((F_{\cal N},\,A_{\cal N})\) and
 \((A_{\cal P},\,G_{\cal P})\) are obstrollable, then the equality 
{\rm (\ref{Lyap2})}
 holds, i.e. the matrix \(S^l\) is a solution of the
 Sylvester-Lyapunov equation {\rm (\ref{SL3})}.
}
\end{theo}

\begin{lemma}
{\sl
Let \(F=(f_1, \,\dots,\, f_n)\) be a \(k\times n\) matrix (i.e.
 \(f_1\), ... , \(f_n\)
  are \(k\) vector-columns), and
 \(A=\mbox{\rm diag}(\alpha_1,\,\dots,\,\alpha_n)\) be a diagonal matrix
with simple spectrum (i.e. all diagonal entries
 \(\alpha_1,\,\dots\,\alpha_n\) are pairwise different). If no
column \(f_1\), ... , \(f_n\) of the matrix \(F\) is equal to zero,
 then the pair (\(F,A\)) is obstrollable.
}
 \end{lemma}

\noindent
PROOF. Actually, the proof of the statement was already done
 (See how we obtained the equality (\ref{kernchain10})).
\hfill\framebox[0.45em]{ }

This lemma (together with Theorem \ref{rightreprtheo})
 shows that Theorem \ref{robstLS1} and Theorem \ref{robstLS2}
 are applicable to rational matrix
functions in general position.

Now we ``inverse'' our reasonings leading from a chain identity to
a Sylvester-Lyapunov equation. Let
 \(F,\,G\) be \(k\times n\) and \(n\times k\) matrices,
 \(A,\,B\) be \(n\times n\) matrices, with spectra \({\sigma}_A\) and 
\({\sigma}_B\),
and
\(S\) be an {\it invertible} \(n\times n\)
matrix for which
 the equality
\begin{equation}
BS-SA=GF
\label{LyGen}
\end{equation}
 holds. Evidently, this equality is equivalent
 to the identity (with respect to \(y\in {\Bbb C}\)):
\begin{displaymath}
S(yI-A)-(yI-B)S\equiv GF.
\end{displaymath}
Multiplying the last identity by the matrix \(S^{-1}(yI-B)^{-1}\)
from the left and by the matrix \((yI-A)^{-1}S^{-1}\)
from the right, we come to the identity 
\begin{displaymath}
S^{-1}(yI-B)^{-1}-(yI-A)^{-1}S^{-1}-
S^{-1}(yI-B)^{-1}G\cdot
F(yI-A)^{-1}S^{-1}\equiv 0.
\end{displaymath}
(This is nothing more than the equality \(M=0\), where \(M\) is defined by
(\ref{kernchain8})). Multiplying the last identity by the matrix
\(F(xI-A)^{-1}\) from the left and by the matrix 
\((zI-B)^{-1}G\)
 from the right, we come to the identity (with respect to
\(x\in {\Bbb C},\,y\in {\Bbb C},\,z\in {\Bbb C} \)):
\begin{eqnarray}
F(xI-A)^{-1}S^{-1}\cdot (yI-B)^{-1}
(zI-B)^{-1}\cdot  G- &\nonumber \\
-F\cdot (xI-A)^{-1}(yI-A)^{-1}\cdot 
S^{-1}(zI-B)^{-1}G-& \nonumber\\
F(xI-A)^{-1}
S^{-1}(yI-B)^{-1}G\times & \nonumber \\
\times F(yI-A)^{-1}S^{-1}
(zI-B)^{-1}G&\equiv0.
\nonumber
\end{eqnarray}
Using Hilbert identities
\begin{eqnarray}
(xI-A)^{-1}(yI-A)^{-1}\equiv -
\frac{(xI-A)^{-1}-(yI-A)^{-1}}{x-y},\nonumber\\
(yI-B)^{-1}(zI-B)^{-1}\equiv -
\frac{(yI-B)^{-1}-(zI-B)^{-1}}{y-z},\nonumber
\end{eqnarray}
we reduce the last identity to the form
\begin{displaymath}
\frac{K(x,z)-K(y,z)}{x-y}-\frac{K(x,y)-K(x,z)}{y-z}\equiv
K(x,y)\cdot K(y,z),
\end{displaymath}
or, what is the same, to the identity
\begin{equation}
(x-y)\,K(x,y)+(y-z)\,K(y,z)-(x-z)\,K(x,z)\equiv -(x-y)(y-z)\,K(x,z).
\label{hyddenchain}
\end{equation}
where {\sf now} the function \(K(\xi,\eta)\) of two variables
is {\sf defined} by the formula
\begin{equation}
K(x,y)\stackrel{\rm\tiny def}{=}
F(x I-A)^{-1}\cdot
S^{-1}\cdot
(y I-B)^{-1}G.
\end{equation}
The  identity (\ref{hyddenchain}) may be rewritten in the form a chain identity
(\ref{chainID}):
\(
T(x,y)\cdot T(y,z)\equiv T(x,z),
\)
where the function \(T(.,.)\) of two variables is {\sf defined} as
\(
T(x,y)\stackrel{\rm\tiny def}{=}I+(x - y )K(x,y),
\)\\
or

\noindent
\(
T(x,y)\stackrel{\rm\tiny def}{=}
I+(x -y )F(x I-A)^{-1}S^{-1}
(y I-B)^{-1}G,\qquad
 {\mathcal D}_T\stackrel{\rm\tiny def}{=}
\Big(\overline{{\Bbb C}}\setminus{\sigma}_A\Big)\times
\Big (\overline{{\Bbb C}}\setminus{\sigma}_B\Big).
\) 

Thus, we proved the following\\[2.0ex]
\begin{theo}\label{chaingener}  
{\sl  Let  \(F,\,G\) be \(k\times n\) and \(n\times k\) matrices, 
 \(A,\,B\) be \(n\times n\) matrices with spectra \({\sigma}_A\)
and \({\sigma}_B\),
and
\(S\) be an {\sf invertible} \(n\times n\)
matrix for which  the equality \(BS-SA=GF\)
  holds.

 Then the matrix function \(T(.,.)\), which is {\sf defined} by
\begin{equation}
T(x,y)\stackrel{\rm\tiny def}{=}
I+(x -y )F(x I-A)^{-1}S^{-1}
(y I-B)^{-1}G,\quad
 {\mathcal D}_T\stackrel{\rm\tiny def}{=}
\Big(\overline{{\Bbb C}}\setminus{\sigma}_A\Big)\times
\Big (\overline{{\Bbb C}}\setminus{\sigma}_B\Big),
\label{defchain}
\end{equation}
satisfies the chain identity {\rm (\ref{chainid})}: \,\,
\[
T(x,\, y)\cdot T(y,\, z)\equiv T(x ,\, z)
\]
  and the diagonal unity identity {\rm (\ref{unid}):} \,\,
\[
T(x , x)\equiv I,
\]
 and hence\,\footnote{According to Theorem \ref{ChainGen}:
the set
 \((\overline{{\Bbb C}}\setminus{\sigma}_A)\cap
 (\overline{{\Bbb C}}\setminus{\sigma}_B)=
 \overline{{\Bbb C}}\setminus ({\sigma}_A\cup{\sigma}_B)\) is not only 
nonempty but also very rich. So, we have many possibilities for choice of 
a distinguished point.( See the proof of the Theorem \ref{ChainGen}).
 In particular, we can choose the point \(\infty\) as a distinguished point.
},
 is of form
\[T(x,y)=R(x)R^{-1}(y),\] where
 \(R(x)\stackrel{\rm\tiny def}{=}T(x,\infty),\, 
(R(y))^{-1}\stackrel{\rm\tiny def}{=}T(\infty, y)\):
\begin{equation}
R(x)\stackrel{\rm\tiny }{=}I-F(xI-A)^{-1}S^{-1}G, \qquad 
R^{-1}(y)\stackrel{\rm\tiny }{=}I+FS^{-1}(yI-B)^{-1}G.
\label{explR}
\end{equation}
are mutually inverse
{\rm (}i.e. \(R(x)R^{-1}(x)\equiv R^{-1}(x)R(x)\equiv I\){\rm )}
rational  matrix functions.
}
\end{theo}

\noindent
PROOF. The chain identity (\ref{chainid}) for the function \(T\),
defined by (\ref{defchain}), was proved immediately before. The diagonal
unity identity (\ref{unid}) evidently follows from the expression
 (\ref{defchain}).
 The equality \(T(x,y)\equiv T(x,\infty)\cdot T(\infty , y))\)
is the special case of the chain identity  (\ref{chainid}) (written for
the triple of points \(x,\infty, y)\)). See Theorem \ref{ChainGen}.
\voffset = - 2.3 truecm
\voffset = - 2.3 truecm
\voffset = - 2.3 truecm

Letting \(y\) tend to \(\infty\), we obtain the expression (\ref{explR})
for the matrix function \(R(x)\stackrel{\rm\tiny def}{=}T(x,\infty)\).
Letting \(x\) tend to \(\infty \), we obtain the expression (\ref{explR})
for the matrix function \(R^{-1}\stackrel{\rm\tiny def}{=}T(\infty,y).\)

That the functions \(R\) and \(R_{-1}\), defined by (\ref{explR})),
are mutually inverse follows from the chain identity written for the triples
\(x,\infty,x\) and \(\infty ,x,\infty \). That the function \(R\) is
rational is evident.
\hfill\framebox[0.45em]{ }

However, Theorem \ref{chaingener} says nothing about the nature
 of the  rational function \(R\). Imposing restrictions on the data
\(A,\,B,\, F,\, G\), we can say more about the matrix function function \(R\).

\begin{theo}\label{chainspecial}
{\sl  Let  \(F\) be \(k\times n\) matrix and \(G\) be \(n\times k\) matrix
  with non-zero columns and non-zero rows respectively, i.e.
\[F=[f_1\,f_2\, \dots\,f_n] ,\qquad G=
\left[
\begin{array}{c}
 g_1
\\
g_2
\\
 \vdots 
\\
 g_n
\end{array}\right]
 ,\]
where no column
 \(f_1,\,f_2,\, \dots\, ,f_n \) and no row \( g_1,\,g_2,\, \dots\, ,g_n\)
 are zero , and let
 \(A,\,B\) be \(n\times n\) be diagonal matrices with simple disjoint spectra,
 i.e.

\vspace{-13pt}
\[A=\mbox{\rm diag}(\lambda_1,\,\dots \, , \lambda_n), \qquad
  B=\mbox{\rm diag}(\mu_1,\,\dots \, , \mu_n),\]

\vspace{-13pt}
\noindent
where \(\lambda_1,\,\dots \, , \lambda_n;\,\,\mu_1,\,\dots \, , \mu_n\) are
pairwise different complex numbers.\\
Assume that the  \(n\times n\) matrix \(S\),
\begin{equation}S=\|s_{p,q}\|_{1\leq p,q\leq n},\qquad 
 s_{p,q}= \frac{g_p\,f_q}{\mu_p-\lambda_q}
\label{s1}
\end{equation}
(which can be obtained from the data \(F,\,G,\,A,\,B\)
as the unique solution of the Sylvester-Lyapunov equation
\(BX-XA=GF\))
 is  {\it invertible}. \\[2pt]
Then:}

\vspace{-5pt}
\begin{enumerate}
\item
{\sl The matrix function \(T(\,.\,,\,.\,)\) of two variables, 
which is {\sf defined} by
$$
T(x,y)\stackrel{\rm\tiny def}{=}
I+(x -y )F(x I-A)^{-1}S^{-1}
(y I-B)^{-1}G,\quad
 {\mathcal D}_T\stackrel{\rm\tiny def}{=}
\Big(\overline{{\Bbb C}}\setminus{\sigma}_A\Big)\times
\Big (\overline{{\Bbb C}}\setminus{\sigma}_B\Big),
\eqno{(\ref{defchain})}
$$
satisfies the chain identity {\rm (\ref{chainid})}:\
\[T(x,y)\cdot T(y,z)\equiv T(x,z),\]
and the diagonal unity identity {\rm (\ref{unid})}:\ \[T(x,x)\equiv I.\]
}
\item
{\sl 
The matrix function \(T(x,y)\) is of the form 
\[T(x,y)=R(x)R^{-1}(y),\] where the matrix functions \(R,\, R^{-1}\) are
{\sf defined} by the formulas
 \(R(x)\stackrel{\rm\tiny def}{=}T(x,\infty),\,
R^{-1}(y)\stackrel{\rm\tiny def}{=}T(\infty,y)\):
$$
R(x)\stackrel{\rm \tiny }{=}I-F(xI-A)^{-1}S^{-1}G, \qquad 
R^{-1}(y)\stackrel{\rm \tiny }{=}I+FS^{-1}(yI-B)^{-1}G.
\eqno{(\ref{explR})}
$$
and are mutually inverse,
  {\rm (}i.e. \(R(x)R^{-1}(x)\equiv R^{-1}(x)R(x)\equiv I\){\rm )}.
}
\item
{\sl
The matrix functions \(R\) and \(R^{-1}\) are 
rational matrix functions in general position.
}
\item
{\sl
The pole set
\({\cal P}(R)\) of the function \(R\) coincides with the set
 \(\{{\lambda}_1, \, \dots \,{\lambda}_n\}\); the zero set
\({\cal N}(R)\) of the function \(R\) coincides with the set
 \(\{{\mu}_1, \, \dots \,{\mu}_n\}\), i.e.
\begin{equation}
A_{\cal P}=A, \quad A_{\cal N}=B,
\end{equation}
where \(A_{\cal P}\) and \(A_{\cal }\) are the pole and zero matrices of
the matrix function \(R\).
}
\item
{\sl 
The semi-residual matrices\,\footnote{\label{f}More precisely, one of the 
representatives of the equivalence class of the set of semi-residual
 matrices of the matrix function \(R\). See Remark \ref{gauge}.}
\ 
\(F_{\cal P},\,G_{\cal P},\,F_{\cal N},\,G_{\cal N}\) of the matrix
function \(R\) can be expressed in terms of the data \(F,\,G\) and 
of the matrix \(S\) (which in its turn is expressible from the data
\(F,\,G,\,A,\,B\)):
\begin{equation}
F_{\cal P}=F;\quad G_{\cal N}=G; \quad
\quad\quad F_{\cal N}=F\cdot S^{-1}; 
\quad G_{\cal P}=-S^{-1}\cdot G; 
\label{finalexpr1}
\end{equation}
}
\item
{\sl
The right zero-pole coupling matrix \(S^r\)  and the  left zero-pole
coupling matrix \(S^l\) for the matrix function \(R\) can be expressed 
in terms of the  matrix \(S\):
\begin{equation}
S^r=S;\qquad S^l=S^{-1}.
\label{finalexpr2}
\end{equation}
}
\end{enumerate}
\end{theo}

\noindent
PROOF.\\
\(\bullet\) Items 1 and 2 of Theorem \ref{chainspecial} are already proved.
(See Theorem \ref{chaingener}).\\
\(\bullet\) Let's prove that the matrix functions
\(R\) and \(R^{-1}\) are in general positions and investigate their 
singularities.
The expression (\ref{explR}) for \(R\) may be written in the form
\[R(x)=I + \sum\limits_{1\leq j \leq n}\frac{R_{\lambda_j}}{x-\lambda_j},\]
with the matrix \(R_{\lambda_j}\) is of the form
\[R_{\lambda_j}=f_j\cdot v_j,\]
where \(v_j\) is \(k\)-th row of the \(n\times k \) matrix 
\(V\stackrel{\rm\tiny def}{=}-S^{-1}G.\)
 From this expression  it follows that \(R\) is
holomorphic outside the points \(\{\lambda_1,\,\dots \, ,\lambda_n\}\)
\voffset = - 2.3 truecm
and its inverse \(R^{-1}\) is holomorphic outside the points 
\(\{\mu_1,\,\dots \, ,\mu_n\}\). 
 Let's focus on  the point \(\lambda_j.\) There are two
possibilities: or \(v_j=0\), or \(v_j\not=0.\) If \(v_j=0\) then 
 \(R_{\lambda_j}=0\) and hence, the function \(R\) is holomorphic at 
the point \(\lambda_j\). If \(v_j\not=0.\) than the matrix
\(R_{\lambda_j}\) is non-zero, and has rank one. (We recall that, according to
the assumptions of Theorem \ref{chainspecial}, \(f_j\not=0\).) We
show now that the equality \(v_j=0\) is impossible. This equality
may be written in the form \(E_jS^{-1}G=0\), where 
\(E_j=\mbox{\rm diag}(\delta_{1j},\,\delta_{2j},\,\dots \, , \delta_{nj})\),
( \( \delta\) is the Kronecker symbol).
Multiplying the identity\,\footnote{Which, in fact, serves as
{\sf the definition} of the matrix \(S\).} 
\(BS-SA=GF\) by the matrix \(E_jS^{-1}\) from the left and by the
matrix \(S^{-1}\) from the right and taking into account that the matrices
\(e_j\) and \(A\) are permutable (both  are diagonal), we come to the
equality \((E_jS^{-1})B-A(E_jS^{-1})=0\). Because the spectra of \(A\) and
\(B\) are disjoint, we obtain that \(E_jS^{-1}=0\), and hence \(E_j=0\).
The contradiction shows that the equality \(v_j=0\) is impossible.
Thus, each point \(\lambda_j\), \(j=1,\,2,\,\dots\,,n\), is simple
pole of the matrix function \(R\), with residue matrices  of rank one.\\
\(\bullet\)
Analogously, we can show that the matrix function \(R^{-1}\) is
holomorphic outside the points \(\{\mu_1,\,\dots \, ,\mu_n\}\) and that
each point \(\mu_j\), \(j=1,\,2,\,\dots\,,n\), is a simple pole of 
the matrix function \(R^{-1}\), with residue matrices  of rank one.
Thus, items 3 and 4 of the claim of Theorem \ref{chainspecial} are proved.\\
\(\bullet\)
Item 5 of the claim follows from the representations (\ref{explR}).
(Compare (\ref{explR}) with (\ref{matradddir})-(\ref{matraddinv})).\\
\(\bullet\)
Now that we have established the relation (\ref{finalexpr1}) we may
rewrite the equality \mbox{\(BS-SA=GF\)} in the form
\(A_{\cal N}S-SA_{\cal P}=G_{\cal N}F_{\cal P}\). Comparing the
last equality  with equality (\ref{Lyap1}), we conclude that \(S^r=S\).
>From (\ref{invs}) it follows now that \(S^l=S^{-1}\).
\hfill\framebox[0.45em]{ }

\vspace{1ex}
The ``hybrid'' version of this theorem
 (see Remark \ref{hybrid} and formulas (\ref{hybr1}),\,(\ref{hybr2}))
 can be
formulated as well. This is the form which is convenient for applications
 in study of the Schlesinger system.

\begin{theo}\label{chainhybr}
{\sl
 Let  \(F\) be \(k\times n\) matrix and \(G\) be \(n\times k\) matrix
  with non-zero columns and non-zero rows respectively, i.e.
\[F=[f_1\,f_2\, \dots\,f_n] ,\qquad G=
\left[
\begin{array}{c}
 g_1
\\
g_2
\\
 \vdots 
\\
 g_n
\end{array}\right]
 ,\]
where no column
 \(f_1,\,f_2,\, \dots\, ,f_n \) and no row \( g_1,\,g_2,\, \dots\, ,g_n\)
 are zero , and let
 \(A,\,B\) be \(n\times n\) be diagonal matrices with simple disjoint spectra,
i.e
\[A=\mbox{\rm diag}(\lambda_1,\,\dots \, , \lambda_n), \qquad
  B=\mbox{\rm diag}(\mu_1,\,\dots \, , \mu_n),\]
where \(\lambda_1,\,\dots \, , \lambda_n;\,\,\mu_1,\,\dots \, , \mu_n\) are
pairwise different complex numbers.\\
Assume that the  \(n\times n\) matrix \(S\),
\begin{equation}S=\|s_{p,q}\|_{1\leq p,q\leq n},\qquad 
 s_{p,q}= \frac{g_p\,f_q}{\lambda_p-\mu_q}
\label{s2}
\end{equation}
{\rm (}which   can be obtained from the data \(F,\,G,\,A,\,B\)
as the  unique solution of the Sylvester-Lyapunov equation
\(AX-XB=GF\){\rm )}\,
 is  {\it invertible}. \\[2pt]
Then:}
\begin{enumerate}
\item
{\sl
The matrix function \(T\) of two variables, which is {\sf defined} by
the formula
\begin{equation}
T(x,y)\stackrel{\rm\tiny def}{=}
I-(x-y)FS^{-1}(xI-A)^{-1}S(yI-B)^{-1}S^{-1}G,
\label{hybrexpr}
\end{equation}
satisfies the chain identity {\rm (\ref{chainid})}:\
\[T(x,y)\cdot T(y,z)\equiv T(x,z),\]
and the diagonal unity identity {\rm (\ref{unid})}:\ \[T(x,x)\equiv I.\]
}
\item
{\sl 
The matrix function \(T(x,y)\) is of the form 
\[T(x,y)=R(x)R^{-1}(y),\] where the matrix functions \(R,\, R^{-1}\) are
{\sf defined} by the formulas
 \(R(x)\stackrel{\rm\tiny def}{=}T(x,\infty),\,
R^{-1}(y)\stackrel{\rm\tiny def}{=}T(\infty,y)\):
$$
R(x)\stackrel{\rm \tiny }{=}I+FS^{-1}(xI-A)^{-1}G, \qquad 
R^{-1}(y)\stackrel{\rm \tiny }{=}I-F(yI-B)^{-1}S^{-1}G.
\eqno{(\ref{explR})}
$$
and are mutually inverse,
  {\rm (}i.e. \(R(x)R^{-1}(x)\equiv R^{-1}(x)R(x)\equiv I\){\rm )}.
}
\item
{\sl
The matrix functions \(R\) and \(R^{-1}\) are 
rational matrix functions in general position.
}
\item
{\sl
The pole set
\({\cal P}(R)\) of the function \(R\) coincides with the set
 \(\{{\lambda}_1, \, \dots \,{\lambda}_n\}\); the zero set
\({\cal N}(R)\) of the function \(R\) coincides with the set
 \(\{{\mu}_1, \, \dots \,{\mu}_n\}\), i.e.
\begin{equation}
A_{\cal P}=A, \quad A_{\cal N}=B,
\end{equation}
where \(A_{\cal P}\) and \(A_{\cal N}\) are the pole and zero matrices of
the matrix function \(R\).
}
\item
{\sl
The semi-residual matrices\,\footnote{See the footnote 
\({}^{{\rm\scriptscriptstyle\ref{f}}}\).}
\ 
\(F_{\cal P},\,G_{\cal P},\,F_{\cal N},\,G_{\cal N}\) of the matrix
function \(R\) can be expressed in terms of the data \(F,\,G\) and of
the matrix \(S\) (which in its turn is expressible from the data
\(F,\,G\,A\,B\)):
\begin{equation}
F_{\cal N}=F;\quad G_{\cal P}=G; \quad
\quad\quad F_{\cal P}=F\cdot S^{-1}; 
\quad G_{\cal N}=-S^{-1}\cdot G; 
\label{hybrexpr1}
\end{equation}
}
\item
{\sl
The right zero-pole coupling matrix \(S^r\)  and the  left zero-pole
coupling matrix \(S^l\) for the matrix function \(R\) can be expressed 
in terms of the matrix \(S\):
\begin{equation}
S^l=S;\qquad S^r=S^{-1}.
\label{hybrexpr2}
\end{equation}
}
\end{enumerate}
\end{theo}

\noindent
PROOF. Theorem \ref{chainhybr} is nothing more then Theorem \ref{chainspecial}
 in ``other variables''.  Let \(A,\, B,\, F,\, G\) be the data of 
Theorem \ref{chainhybr}, and \(S\) be the matrix (\ref{s2}) generated
by this data. Let's introduce the matrices 
\[\tilde{F}=FS^{-1},\quad
 {\tilde{G}}=-S^{-1}G,\quad \tilde{S}=S^{-1}.
\]
Equality \(AS-SB=GF\) rewritten in terms of
 \(A,B,\tilde{F},\tilde{G},\tilde{S}\) becomes the form
\(B\tilde{S}-\tilde{S}A=\tilde{G}\tilde{F}.\)
>From the last equality it is easy to see that no  column of the matrix
\(\tilde{F}\) and no row of the matrix \({\tilde{G}}\) equals zero:
the equality \(E_j{\tilde{G}}=0\) or \(\tilde{F}E_j=0\),
where
 \(E_j=\mbox{\rm diag}(\delta_{1j},\,\delta_{2j},\,\dots \, , \delta_{nj})\),
( \( \delta\) is the Kronecker symbol), leads to the equality
\(E_j\tilde{S}=0\) or \(\tilde{S}E_j=0\), what contradicts to the
invertibility of \(\tilde{S}\). Now Theorem \ref{chainspecial}, applied
to the matrix function
\[T(x,y)
=I+(x-y)\tilde{F}(xI-A)^{-1}(\tilde{S})^{-1}(yI-B)^{-1}{\tilde{G}},\]
gives the chain and diagonal unity identities for this \(T\)
and the factorization \(T(x,y)=R(x)R^{-1}(y)\), as well as 
the expressions for the semi-residual matrices 
\(F_{\cal P}, G_{\cal P}, F_{\cal N}, G_{\cal N}\)
 of the matrix function \(R\):
\[F_{\cal P}=\tilde{F},\quad 
 {G_{\cal N}}={\tilde{G}}, \quad
F_{\cal N}=\tilde{F}\cdot {\tilde{S}}^{-1}, \quad
{G_{\cal P}}=-{\tilde{S}}^{-1}{\tilde{G}}.
\]
Rewritten in terms of \(F, G, S\), these relations becomes the form
(\ref{hybrexpr1}).
\hfill\framebox[0.45em]{ }

\vspace*{0.6cm}

\noindent
\vspace*{0.2cm} 
\begin{minipage}{15.0cm}
\section{\hspace{-0.4cm}.\hspace{0.19cm} THE SYSTEM REPRESENTATION AS A
TOOL FOR THE SPECTRAL (WIENER-HOPF)
 FACTORIZATION OF MATRIX FUNCTIONS.}
\end{minipage}\\[-0.5cm]
\setcounter{equation}{0}

In this section we show that the system representation may be used as
 an efficient tool for the so called {\sf spectral factorization}
(or the  {\sf Wiener-Hopf factorization}) of
a matrix function. The problem of the spectral factorization can
be formulated in the following way.\\[2.0ex]
GEOMETRIC CONFIGURATION. {\sl
 In the extended complex plane \(\overline{\Bbb C}\)
 a simple closed contour \(\Gamma\) is given. This contour separates
\(\overline{\Bbb C}\) inter two regions, \(G_{+}\) and \(G_{-}\). 
These regions \(G_{+}\) and \(G_{-}\) are connected open sets. We assume that
 the point
\(\{\infty\}\) does not belong to the contour \(\Gamma\), thus
 one of the components, say \(G_{-}\), contains the point \(\{\infty\}\).}

\begin{defi}\label{spectrFact}
{\sl
Given a \(k\times k\) matrix
function \(\Phi\) on the contour \(\Gamma\), the factorization  of the form
\begin{equation}
\Phi(\zeta)={\Phi}_{+}(\zeta)\cdot {\Phi}_{-}(\zeta)\quad (\zeta \in \Gamma),
\label{WHF}
\end{equation}
where \({\Phi}_{+}\) and \({\Phi}_{-}\) are \(k\times k\) matrix functions,
the matrix function \({\Phi}_{+}\) and its inverse \(({\Phi}_{+})^{-1}\) are
holomorphic on  \(G_{-}\cup\Gamma\), and
the matrix function \({\Phi}_{-}\) and its inverse \(({\Phi}_{-})^{-1}\) are
holomorphic%
\footnote{In particular,  the functions 
\({\Phi}_{+}\) and \({\Phi}_{-}\) are  holomorphic on the common
 boundary \(\Gamma\) of the domains \(G_{+}\) and \(G_{-}\), so
the relation (\ref{WHF}) makes sense.}%
\ on \(G_{+}\cup\Gamma\), is said to be {\sf the spectral factorization
 ({\sl  or} the Wiener-Hopf factorization) of the matrix-function
\(\Phi\) with respect to \(\Gamma\).
\hfill\framebox[0.45em]{ }
}
}
\end{defi}\\[1.6ex]
We impose the normalizing condition
\begin{equation}
\Phi_{+}(\infty)=I.
\label{normPhi+}
\end{equation}
on the factor  \(\Phi_{+}\).
(The function \(\Phi_{+}\) is holomorphic
and invertible at the point \(\infty\), so the condition (\ref{normPhi+})
 makes sense).

 {\sf Under the normalizing condition (\ref{normPhi+}),
 the spectral factorization (\ref{WHF}) is unique}.
             
  Even in the scalar case
\(k=1\) (i.e. \(\Phi\) is a complex valued function) the factorization problem
(\ref{WHF}) is not always solvable: there is a topological obstacle for
the solvability. For a smooth nonvanishing complex valued function \(\Phi\) on
\(\Gamma\), the factorization problem (\ref{WHF}) solvable if and only if 
there exists an univalued continuous branch of the function
 \(\ln\Phi (\zeta )\)
on \(\Gamma\). If this condition is fulfilled,
the solution of the factorization problem may be expressed in terms of the data
\(\Phi\) by the formula
\begin{equation}
\Phi_{\pm}(\zeta)=\mbox{\rm exp}\bigg\{\pm \frac{1}{2\pi i}\int\limits_{\Gamma}
\frac{\ln \Phi (t)}{t-z}\,dt\bigg\} \qquad (z\in G_{\pm}).
\label{SoPl}
\end{equation}
The proof of the fact, that the formula (\ref{SoPl}) gives the solution of the
factorization problem (\ref{WHF}) in the scalar case, is based essentially on
the Sokhotski\u{\i}-Plemelj formulas on the boundary behavior of the
Cauchy integral. Actually, in the scalar case we solve the additive problem
\[
\Psi_{+}(\zeta)+\Psi_{-}(\zeta)=\ln\Phi(\zeta) \qquad (\zeta\in \Gamma),
\]
and then we exponentiate. In the matricial case we still can solve the
appropriate additive problem using the Cauchy integral, but exponentiating
 does not lead to the desirable result:
In view of noncommutativity of the matricial multiplication,
\(\mbox{\rm exp}\{A+B\}\not=\mbox{\rm exp}\{A\}\cdot \mbox{\rm exp}\{B\}\) for
 matrices \(A\) and \(B\) in general.
In the matricial case, the situation with the factorization problem
(\ref{WHF}) is much more complicated than in the scalar case.
There are not only topological obstacles to the solvability of this problem.
The factorization problem (\ref{WHF}) is equivalent to some system of
singular integral equations on \(\Gamma\), and to analyze this system is
approximately so hard as to investigate the original factorization problem 
(\ref{WHF}). The factorization problem (\ref{WHF}) appeared firstly in the
 context
 of Hilbert's twenty-first problem: to construct a Fuchsian linear differential
 system with the prescribed monodromy group. See \cite{Pl2}, \cite{Bo}
and \cite{Gah} for details and historical references. The factorization of the
 type (\ref{WHF}) is used also in the solving systems of singular integrals
 equations with Cauchy kernel on the contour as well as for solving 
of systems of integral equations  which kernel depends
on the differences of the arguments on the half axis. See
\cite{Vek}, \cite{GoKr} and \cite{ClGo} on this subject. It should be
 mentioned
that the pioneer papers \cite{Bir1} and \cite{Bir2} by G.\,Birkgoff
had a profound impact on the further investigations on matrix factorization.

We confine ourself to the case the function \(\Phi\) is a rational matrix
 function
(or, more precisely, the restriction on \(\Gamma\) of a rational matrix
 function)
such that the functions \(\Phi\) and \(\Phi^{-1}\) are holomorphic on the
 contour
\(\Gamma\). In this case the factorization (\ref{WHF}) is {\sf global}, i.e.
the  matrix functions \(\Phi_{+}\) and \(\Phi_{-}\) are rational, and
the equality
\begin{equation}
\Phi(z)=\Phi_{+}(z)\cdot \Phi_{-}(z) \qquad
 (\forall z \in \overline{\Bbb C })
\label{globFactoriz}
\end{equation}
holds. Indeed, in this case the function \({\Phi_{+}}^{-1}\Phi\) is
 holomorphic 
within \(G_{-}\) except finite many poles located on the set
 \({\cal P}(\Phi)\cap G_{-}\). In view view of (\ref{WHF}), this function
 continues analytically into
the function \(\Phi_{-}\) which is holomorphic on \(G_{+}\cup \Gamma\).
Thus, the matrix-function \(\Phi_{-}\) has no other singularities in 
\(\overline{\Bbb C}\) than finite many poles and hence is rational. 
For the same reasoning, the matrix function \(\Phi_{+}\) is rational.

Thus, in the case that the initial matrix function \(\Phi\) is rational,
the problem of the spectral factorization may be reformulated
in the following manner:\\[2.5ex]
\renewcommand{\thedefi}{%
\mbox{\ref{spectrFact}\(\,{}^{\prime}\)}}%
\begin{defi}
{\sl 
Given a \(k\times k\) rational matrix function \(\Phi\),
\(\mbox{\rm det}\Phi\not\equiv 0\), its factorization of the form
{\rm (\ref{globFactoriz})}, where \(\Phi_{+},\,\Phi_{-}\) are rational matrix
functions with zero and pole location
\begin{equation}
 {\cal P}(\Phi_{+})\subset G_{+}, \quad
 {\cal N}(\Phi_{+})\subset G_{+},   \quad
 {\cal P}(\Phi_{-})\subset G_{-},   \quad
 {\cal N}(\Phi_{-})\subset G_{-},
\label{LocSingFac}
\end{equation}
is said to be {\sf the spectral factorization
 ({\sl  or} the Wiener-Hopf factorization) of the matrix-function
\(\Phi\) with respect to \(\Gamma\).
\hfill\framebox[0.45em]{ }
}
}
\end{defi}
\addtocounter{defi}{-1}
\renewcommand{\thedefi}{%
\mbox{\arabic{section}.\arabic{defi}}}%

\vspace{0.5em}
\noindent
We consider even the more special case: {\sf the function \(\Phi\) is a
 rational
function in general position}. In this case the calculation of the factors
\(\Phi_{+}\) and \(\Phi_{-}\) can be performed completely by hand, in terms of
 poles 
and ``zeros'' of the matrix function \(\Phi\) and its semiresidual vectors.

So, let \(\Phi\) be a rational matrix function in general position, normalized
by the condition
\begin{equation}
\Phi(\infty)=I.
\label{normPhi}
\end{equation}
Let 
\({\cal P}(\Phi)\) and \({\cal N}(\Phi)\) be its pole and zero sets,
\(A_{{\cal P}}(\Phi)\) and \(A_{{\cal N}}(\Phi)\) be its pole and zero
 matrices,
\(F_{\cal P}(\Phi)\), \(G_{\cal P}(\Phi)\), \(F_{\cal N}(\Phi)\),
 \(G_{\cal N}(\Phi)\)
be the appropriate semiresidual matrices. 
According to the Theorem \ref{summarTheo}, {\sf  the zero-pole coupling
 matrices
\(S^r(\Phi)\) and \(S^l(\Phi)\) are invertible}, and the matrix-function 
\(\Phi\) 
admits the representations of the form  (\ref{forRs}) and (\ref{forRR^{-1}s1}):
\begin{equation}
\Phi(z)=
I-F_{\cal P}(\Phi)(zI-A_{{\cal P}(\Phi)})^{-1}{S^r(\Phi)}^{-1}G_{\cal N}(\Phi)
\label{RSystForPhi}
\end{equation}

\vspace*{-3.4ex}
\noindent
and 
\begin{equation}
\Phi(z)=I+
F_{\cal N}(\Phi){S^l(\Phi)}^{-1}(zI-A_{{\cal P}(\Phi)})^{-1}G_{\cal P}(\Phi).
\label{LSystForPhi}
\end{equation}
The Sylvester-Lyapunov equations for the matrices \(S^r(\Phi)\) and
 \(S^l(\Phi)\) (which actual-y are
the {\sf definitions} of these matrices) are of the form:
\begin{equation}
A_{\cal N}(\Phi)S^r(\Phi)-S^r(\Phi)A_{\cal P}(\Phi)=
G_{\cal N}(\Phi)F_{\cal P}(\Phi),
\label{SLRforPhi}
\end{equation}
and
\begin{equation}
A_{\cal P}(\Phi)S^l(\Phi)-S^l(\Phi)A_{\cal N}(\Phi)=
G_{\cal P}(\Phi)F_{\cal N}(\Phi).
\label{SLLforPhi}
\end{equation}
Moreover, the matrices \(S^r(\Phi)\) and \(S^l(\Phi)\) satisfy the equality
\begin{equation}
 S^l(\Phi)\cdot S^r(\Phi)=I.
\label{invsForPhi}
\end{equation}
(This is (\ref{invs}) for the matrix function \(\Phi\)).
The zero-pole coupling relation hold:
\begin{equation}
\begin{array}{ll}
\mbox{\rm a)}.\,\,     G_{\cal N}(\Phi)=-S^r(\Phi)G_{\cal P}(\Phi); \quad&
\mbox{\rm b)}.\,\,     G_{\cal P}(\Phi)=-S^l(\Phi)G_{\cal N}(\Phi);\quad
\\[1.5ex]
\mbox{\rm c)}.\,\,     F_{\cal P}(\Phi)= F_{\cal N}(\Phi)S^r(\Phi); \quad&
\mbox{\rm d)}.\,\,     F_{\cal N}(\Phi)= F_{\cal P}(\Phi)S^l(\Phi).
\end{array}
\label{z-prelforPhi}
\end{equation}
(This is (\ref{recoup}) for the matrix-function \(\Phi\)).

 According to the assumptions,
\begin{equation}
{\rm a).}\ \ {\cal P}(\Phi)\cap \Gamma =\emptyset;\qquad\qquad\qquad
{\rm b).}\ \  {\cal N}(\Phi)\cap \Gamma =\emptyset.
\label{Count_assump}
\end{equation}
Assume that the factorization (\ref{globFactoriz})
holds, and that the normalizing condition (\ref{normPhi+}) is satisfied.

 Since the functions \(\Phi_{-}\), \(\Phi^{-1}_{-}\) are
 holomorphic in \(G_{+}\) and the functions \(\Phi\), \(\Phi^{-1}\) have only
simple poles in \(G_{+}\), from the relations
\begin{equation}
{\rm a).}\ \ \Phi_{+}=\Phi \cdot \Phi^{-1}_{-},\qquad \qquad \qquad {\rm b).}\
 \ \Phi_{+}^{-1}=\Phi_{-}\cdot \Phi^{-1}
\label{divRel}
\end{equation}
it follows that the functions \(\Phi_{+}\), \(\Phi_{+}^{-1}\) have only simple
 poles
 in \(G_{+}\). In \(G_{-}\) and on \(\Gamma\) the functions \(\Phi_{+}\) and
\(\Phi_{+}^{-1}\) don't have singularities at all. Thus, 
\begin{equation}
{\cal P}(\Phi_{+})={\cal P}(\Phi)\cap G_{+}, \qquad
{\cal N}(\Phi_{+})={\cal N}(\Phi)\cap G_{+}.
\label{locSing+}
\end{equation}

Let \(\lambda\in{\cal P}(\Phi_{+})\). From (\ref{divRel}.a) it follows that
 the residues
\(R_{\lambda}(\Phi_{+})\) and \(R_{\lambda}(\Phi)\) of the matrix functions
 \(\Phi_{+}\) and
\(\Phi\) at the point \(\lambda\) are related by the equality
\begin{equation}
R_{\lambda}(\Phi_{+})=R_{\lambda}(\Phi)\cdot ({\Phi_{-}}(\lambda))^{-1} 
\qquad\qquad (\ \forall \lambda\in {\cal P}(\Phi_{+})\ ).
\label{relRes+}
\end{equation}
Since \(\Phi\) is a matrix function in general position, the rank of the
 residue matrix 
 \(R_{\lambda}(\Phi)\) is equal to one. According to the assumptions, the
 matrix 
\({\Phi_{-}}(\lambda)\) is invertible for \(\lambda\in G_{+}\).
>From (\ref{relRes+}) it follows now that the rank of the residue matrix
\( R_{\lambda}(\Phi_{+})\) is equal to one as well. 
Let now \(\mu\in{\cal N}(\Phi_{+})\).
>From (\ref{divRel}.b) it follows that the residues
\(R_{\mu}(\Phi_{+})\) and \(R_{\mu}(\Phi)\) of the matrix functions 
\(\Phi_{+}^{-1}\) and
\(\Phi^{-1}\) at the point \(\mu\) are related by the equality
\begin{equation}
R_{\mu}(\Phi_{+})= {\Phi_{-}}(\mu)\cdot R_{\mu}(\Phi)
\qquad\qquad (\ \forall \mu\in {\cal N}(\Phi_{+})\ ). 
\label{relResInv+}
\end{equation}
The rank of the matrix \(R_{\mu}(\Phi)\) is equal to one (\(\Phi\) is a matrix
 function
 in general
position);
 the matrix \({\Phi_{-}}(\mu)\)
is invertible (according to the assumptions, the matrix \(\Phi_{-}(\mu)\) is
 invertible
 for
\(\mu\in G_{+}\)). From (\ref{relResInv+}) it follows now, that
 the the rank of the residue matrix \(R_{\mu}(\Phi_{+})\) is equal to one as
 well.
>From (\ref{locSing+}) it follows that
\[
{\cal P}(\Phi_{+})\cap {\cal N}(\Phi_{+})=\emptyset; \qquad
\{\infty\}\not\in {\cal P}(\Phi_{+}), \quad \{\infty\}\not\in
 {\cal N}(\Phi_{+}).
\]
{\sf Thus, \(\Phi_{+}\) is a rational matrix function in general position}.

In the same way we obtain that
{\sf  \(\Phi_{-}\) is a rational matrix function in general position}, and
\begin{equation}
{\cal P}(\Phi_{-})={\cal P}(\Phi)\cap G_{-}, \qquad
{\cal N}(\Phi_{-})={\cal N}(\Phi)\cap G_{-}.
\label{locSing-}
\end{equation}
Moreover, the residues \(R_{\lambda}(\Phi_{-})\) and \(R_{\mu}(\Phi_{-})\) of
the matrix functions \(\Phi_{-}\) and \(\Phi_{-}^{-1}\) at the poles
 \(\lambda\in{\cal P}(\Phi_{-})\)
and \(\mu\in{\cal N}(\Phi_{-})\,(\,={\cal P}(\Phi_{-}^{-1})\, )\) are related
to the residues of the matrix functions
\(\Phi\) and \(\Phi^{-1}\) at the same points by the equalities
\begin{equation}
R_{\lambda}(\Phi_{-})= ({\Phi_{+}}(\lambda))^{-1}\cdot R_{\lambda}(\Phi)
\qquad\qquad (\ \forall \lambda\in {\cal P}(\Phi_{-})\ )
\label{relRes-}
\end{equation}
and
\begin{equation}
R_{\mu}(\Phi_{-})= R_{\mu}(\Phi)\cdot {\Phi_{+}}(\mu)
\qquad\qquad (\ \forall \mu\in {\cal N}(\Phi_{+}). 
\label{relResInv-}
\end{equation}
>From (\ref{relRes+}) and (\ref{relResInv+}) it follows that the left
 semiresidual vectors of the
 matrix functions
\(\Phi_{+}\) and \(\Phi\) at the poles \(\lambda\in {\cal P}(\Phi_{+})\)
 coincide:
\begin{equation}
f_{\lambda}(\Phi_{+})=f_{\lambda}(\Phi) \qquad\qquad (\,\forall\lambda\in
 {\cal P}(\Phi_{+})\,). 
\label{ResRel+}
\end{equation}
and the right semiresidual vectors of the matrix functions \(\Phi_{+}^{-1}\)
 and \(\Phi^{-1}\)
at the poles \(\mu\in {\cal N}(\Phi_{+})\) coincide:
\begin{equation}
g_{\mu}(\Phi_{+}^{-1})=g_{\mu}(\Phi^{-1})
 \qquad\qquad (\,\forall\mu\in {\cal N}(\Phi_{+})\,). 
\label{ResInvRel+}
\end{equation}
In the same way we can obtain that the right semiresidual vectors of the
 matrix functions
\(\Phi_{-}\) and \(\Phi\) coincide:
\begin{equation}
g_{\lambda}(\Phi_{-})=g_{\lambda}(\Phi) \qquad\qquad (\,\forall\lambda\in
 {\cal P}(\Phi_{-})\,). 
\label{ResRel-}
\end{equation}
and the left semiresidual vectors of the matrix functions \(\Phi_{-}^{-1}\)
 and \(\Phi^{-1}\)
coincide:
\begin{equation}
f_{\mu}(\Phi_{-}^{-1})=f_{\mu}(\Phi^{-1}) \qquad\qquad (\,\forall\mu\in
 {\cal N}(\Phi_{-})\,). 
\label{ResInvRel-}
\end{equation}
{\sf The eqialities (\ref{ResRel+}) -- (\ref{ResInvRel-}) are crucial
for solving of the considered factorization problem.}

According to Lemma \ref{isocard}, the equality
\[
\#{\cal P}(\Phi_{+})=\#{\cal N}(\Phi_{+})
\]
holds for the rational matrix function in general position \(\Phi_{+}\).
 Taking into account
the equality (\ref{locSing+}), we obtain the following equality
\begin{equation}
\#({\cal P}(\Phi)\cap G_{+})=\#({\cal N}(\Phi)\cap G_{+})
\quad\stackrel{\rm\tiny def}{=}n_{+}.
\label{cardequal+}
\end{equation}
Of course, the equality
\begin{equation}
\#({\cal P}(\Phi)\cap G_{-})=\#({\cal N}(\Phi)\cap G_{-})
\quad\stackrel{\rm\tiny def}{=}n_{-}.
\label{cardequal-}
\end{equation}
holds as well.

To simplify notations, we denote
\begin{equation}
{\cal P}_{+}\stackrel{\rm\tiny def}{=}{\cal P}(\Phi)\cap  G_{+};\quad
{\cal P}_{-}\stackrel{\rm\tiny def}{=}{\cal P}(\Phi)\cap  G_{-};\quad
{\cal N}_{+}\stackrel{\rm\tiny def}{=}{\cal N}(\Phi)\cap  G_{+};\quad
{\cal N}_{-}\stackrel{\rm\tiny def}{=}{\cal N}(\Phi)\cap  G_{-}.
\label{PartSets}
\end{equation}
To the decompositions
\renewcommand{\arraystretch}{1.3}
\[
\begin{array}{cc}
{\cal P}(\Phi)={\cal P}_{+}\cup{\cal P}_{-},&
{\cal P}_{+}\cap{\cal P}_{-}=\emptyset ,\\
{\cal N}(\Phi)={\cal N}_{+}\cup{\cal N}_{-},&
{\cal N}_{+}\cap{\cal N}_{-}=\emptyset 
\end{array}
\]
of the pole and zero sets \({\cal P}(\Phi)\), \({\cal N}(\Phi)\)
 of the matrix function \(\Phi\) there
correspond natural block-decompositions of the matrices  which appear
in the system representations of the matrix functions \(\Phi\) and 
\({\Phi}^{-1}\): 
 pole and zero matrices
\(A_{\cal P}(\Phi)\), \(A_{\cal N}(\Phi)\), the semiresidual matrices
\(F_{\cal P}(\Phi), F_{\cal N}(\Phi),
 G_{\cal P}(\Phi), G_{\cal N}(\Phi)\) as well as the
zero-pole coupling matrices \(S^r(\Phi)\) and \(S^l(\Phi)\).

Namely, the decompositions of the pole and zero matrices
\(A_{\cal P}(\Phi)={\rm diag}(\lambda_l)_{\lambda_l\in{\cal P}(\Phi)}\)
 and \(A_{\cal N}(\Phi)={\rm diag}(\mu_l)_{\mu_l\in{\cal N}(\Phi)}\)
(of the dimension \(n\times n\)) are of the form: 
\begin{equation} 
A_{\cal P}(\Phi)=
\left[
\begin{array}{cc}
A_{\cal P}(\Phi)_1 & 0 \\
0 & A_{\cal P}(\Phi)_2
\end{array}
\right],
\label{PolMatrBlockDec}
\end{equation}
\begin{equation} 
A_{\cal N}(\Phi)=
\left[
\begin{array}{cc}
A_{\cal N}(\Phi)_1 & 0 \\
0 & A_{\cal N}(\Phi)_2
\end{array}
\right],
\label{ZerMatrBlockDec}
\end{equation}
where
\[
A_{\cal P}(\Phi)_1={\rm diag}(\lambda_l)_{\lambda_l\in{\cal P}_{+}};\qquad
A_{\cal P}(\Phi)_2={\rm diag}(\lambda_l)_{\lambda_l\in{\cal P}_{-}};
\]
\[
A_{\cal N}(\Phi)_1={\rm diag}(\mu_l)_{\mu_l\in{\cal N}_{+}};\qquad
A_{\cal N}(\Phi)_2={\rm diag}(\mu_l)_{\mu_l\in{\cal N}_{-}};
\]
\(A_{\cal P}(\Phi)_1\) and \(A_{\cal N}(\Phi)_1\) are diagonal 
\(n_{+}\times n_{+}\) matrices, 
\(A_{\cal P}(\Phi)_2\) and \(A_{\cal N}(\Phi)_2\) are diagonal 
\(n_{-}\times n_{-}\) matrices;
(\(n_{+}\) and \(n_{-}\) are defined in (\ref{cardequal+}) and
(\ref{cardequal-}):
 \(n_{+}=\#{\cal P}_{+}=\#{\cal N}_{+};\ \ 
n_{-}=\#{\cal P}_{-}=\#{\cal N}_{-}).\)
 
The block-decompositions of the  semiresidual matrices
\(F_{\cal P}(\Phi), F_{\cal N}(\Phi), G_{\cal P}(\Phi), G_{\cal N}(\Phi)\): 
\[
\begin{array}{ccc}
F_{\cal P}(\Phi)={\rm row}(f_{\lambda_l})_{\lambda_l\in{\cal P}(\Phi)},&
F_{\cal N}(\Phi)={\rm row}(f_{\mu_l})_{\mu_l\in{\cal N}(\Phi)}&
{\rm are}\quad k\times n \quad {\rm matrices},\\
 G_{\cal P}(\Phi)={\rm col}(g_{\lambda_l})_{\lambda_l\in{\cal P}(\Phi)},&
 G_{\cal N}(\Phi)={\rm col}(g_{\mu_l})_{\mu_l\in{\cal N}(\Phi)}&%
{\rm are}\quad n\times k\quad {\rm matrices},
\end{array}\]
related to the matrix-function \(\Phi\),
are of the form:
\begin{equation}
{\rm p).}\quad F_{\cal P}(\Phi)= 
\left[F_{\cal P}(\Phi)_1 \, \, \, F_{\cal P}(\Phi)_2\right];
\qquad\qquad
{\rm n).}\quad F_{\cal N}(\Phi)= 
\left[F_{\cal N}(\Phi)_1 \, \, \, F_{\cal N}(\Phi)_2\right].
\label{BlockResMatrL}
\end{equation}
and
\begin{equation}
{\rm p).}\quad G_{\cal P}(\Phi)= 
\left[
\begin{array}{cc}
G_{\cal P}(\Phi)_{\,1}\\
\vspace*{-2.0ex}
\\ 
 G_{\cal P}(\Phi)_{\,2}
\end{array}
\right];
\qquad\qquad
{\rm n).}\quad G_{\cal N}(\Phi)= 
\left[
\begin{array}{cc}
G_{\cal N}(\Phi)_{\,1}\\
\vspace*{-2.0ex}
\\ 
 G_{\cal N}(\Phi)_{\,2}
\end{array}
\right],
\label{BlockResMatrR}
\end{equation}
where
\renewcommand{\arraystretch}{1.3}
\[
\begin{array}{ccccl}
F_{\cal P}(\Phi)_1={\rm row}(f_{\lambda_l})_{\lambda_l\in{\cal P}_{+}},&
F_{\cal N}(\Phi)_1={\rm row}(f_{\mu_l})_{\mu_l\in{\cal N}_{+}}&
\mbox{\rm  are}&k \times  n_{+} & \mbox{\rm matrices};\\
F_{\cal P}(\Phi)_2={\rm row}(f_{\lambda_l})_{\lambda_l\in{\cal P}_{-}},&
F_{\cal N}(\Phi)_2={\rm row}(f_{\mu_l})_{\mu_l\in{\cal N}_{-}}&
\mbox {\rm are}&
 k\times n_{-}& \mbox{\rm matrices};\\
G_{\cal P}(\Phi)_1={\rm col}(g_{\lambda_l})_{\lambda_l\in{\cal P}_{+}},&
G_{\cal N}(\Phi)_1={\rm col}(g_{\mu_l})_{\mu_l\in{\cal N}_{+}}&
\mbox {\rm are}&\ n_{+} \times  k&  \mbox{\rm matrices};\\
G_{\cal P}(\Phi)_2={\rm col}(g_{\lambda_l})_{\lambda_l\in{\cal P}_{-}},&
G_{\cal N}(\Phi)_2={\rm col}(g_{\mu_l})_{\mu_l\in{\cal N}_{-}}&
\mbox {\rm are}& n_{-} \times k&  \mbox{\rm matrices}.
\end{array}
\]
>From (\ref{locSing+}) and  (\ref{locSing-}) it follows that
\begin{equation}
{\rm p).}\  A_{\cal P}(\Phi_{+})=A_{\cal P}(\Phi)_{\,1} \qquad
{\rm n).}\  A_{\cal N}(\Phi_{+})=A_{\cal N}(\Phi)_{\,1}
\label{RightPoleMatr+}
\end{equation}

\vspace*{-4.0ex}
\noindent
and
\begin{equation}
{\rm p).}\ A_{\cal P}(\Phi_{-})=A_{\cal P}(\Phi)_{\,2}\qquad
{\rm n).}\  A_{\cal N}(\Phi_{-})=A_{\cal N}(\Phi)_{\,2}.
\label{RightPoleMatr-}
\end{equation}
In view of (\ref{ResRel+}) and (\ref{ResInvRel+}),
\begin{equation}
{\rm p).}\quad F_{\cal P}(\Phi_{+})=F_{\cal P}(\Phi)_{\,1},\qquad
{\rm n).}\quad G_{\cal N}(\Phi_{+})=G_{\cal N}(\Phi)_{\,1}.
\label{PlusSemiRes}
\end{equation}
In view of (\ref{ResRel-}) and (\ref{ResInvRel-}),
\begin{equation}
{\rm p).}\quad F_{\cal N}(\Phi_{-})=F_{\cal N}(\Phi)_{\,2},\qquad
{\rm n).}\quad G_{\cal P}(\Phi_{-})=G_{\cal P}(\Phi)_{\,2}.
\label{MinusSemiRes}
\end{equation}

Thus, we have expressed the pole matrix \(A_{\cal P}(\Phi_{+})\)
 and  semiresidual matrices \(F_{\cal P}(\Phi_{+})\), \(G_{\cal N}(\Phi_{+})\)
for the left factor \(\Phi_{+}\) in terms of blocks of the appropriate
block-decompositions of the pole matrix \(A_{\cal P}(\Phi)\) and
and semiresidual matrices  \(F_{\cal P}(\Phi)\), \(G_{\cal N}(\Phi)\) for
the factorized matrix function \(\Phi\). We also have expressed
the zero matrix \(A_{\cal N}(\Phi_{-})\)
 and semiresidual matrices \(F_{\cal N}(\Phi_{-})\), \(G_{\cal P}(\Phi_{-})\)
for the right factor \(\Phi_{-}\) in terms of blocks  of the matrices
\(A_{\cal N}(\Phi)\), \(F_{\cal N}(\Phi)\), \(G_{\cal P}(\Phi)\). 
In principle, these data are sufficient to recover the factors
\(\Phi_{+}\) and \(\Phi_{-}\) (\,from the appropriate blocks of the pole
and semiresidual matrices for the factored matrix function \(\Phi\)\,).
To carry out the recovering, we have to solve the Sylvester-Lyapunov
equations (\ref{Lyap1}) and (\ref{Lyap2}) for the matrix functions
\(\Phi_{+}\) and \(\Phi_{-}\) respectively to find from these equations
the zero-pole coupling matrices \(S^r(\Phi_{+})\) and \(S^l(\Phi_{-})\):
\begin{eqnarray}
A_{\cal N}(\Phi_{+})S^r(\Phi_{+})-S^r(\Phi_{+})A_{\cal P}(\Phi_{+})&=&
G_{\cal N}(\Phi_{+})F_{\cal P}(\Phi_{+}),\label{LSfact+}
\\
A_{\cal P}(\Phi_{-})S^l(\Phi_{-})-S^l(\Phi_{-})A_{\cal P}(\Phi_{-})&=&
G_{\cal P}(\Phi_{-})F_{\cal N}(\Phi_{-}).
\label{LSfact-}
\end{eqnarray}
According to the assertion 1 of Theorem \ref{summarTheo} (actually, according
to Theorem \ref{invtheo}: see (\ref{invrel})), {\sf their solutions
\(S^r(\Phi_{+})\)  and \(S^r(\Phi_{-})\) are invertible matrices.}
Then we construct the factors \(\Phi_{+}\) and \(\Phi_{-}\) according
to the formulas (\ref{forRs}) and  (\ref{forRR^{-1}s1}):
\begin{equation}
\Phi_{+}(z)=I-F_{\cal P}(\Phi_{+})\big(zI-A_{\cal P}(\Phi_{+})\big)^{-1}
(S^r(\Phi_{+}))^{-1}G_{\cal N}(\Phi_{+})
\label{forPhiR}
\end{equation}

\vspace*{-2.8ex}
\noindent
and
\begin{equation}
\Phi_{-}(z)=I+F_{\cal N}(\Phi_{-})(S^l(\Phi_{-}))^{-1}
\big(zI-A_{\cal P}(\Phi_{-})\big)^{-1}G_{\cal P}(\Phi_{-}).
\label{forPhiL}
\end{equation}
We express now the matrices \(S^r(\Phi_{+})\) and \(S^l(\Phi_{-})\) in
terms of blocks of the matrices \(S^r(\Phi)\) and \(S^l(\Phi)\).
The block-decomposition of the  zero-pole coupling matrices 
\(S^r(\Phi)\) and \(S^l(\Phi)\) are of the form:
\begin{equation}
S^r(\Phi)=
\left[
\begin{array}{cc}
S^r(\Phi)_{\,11}& S^r(\Phi)_{\,12}\\
S^r(\Phi)_{\,21}& S^r(\Phi)_{\,22}
\end{array}
\right],
\label{blocDecS^r}
\end{equation}
\begin{equation}
S^l(\Phi)=
\left[
\begin{array}{cc}
S^l(\Phi)_{\,11}& S^l(\Phi)_{\,12}\\
S^l(\Phi)_{\,21}& S^l(\Phi)_{\,22}
\end{array}
\right],
\label{blocDecS^l}
\end{equation}
where \(11\)- block-entries \(S^r(\Phi)_{\,11}\),\ \(S^l(\Phi)_{\,11}\)
  are \(n_{+}\times n_{+}\)
 matrices, and
\(22\)- block-entries \(S^r(\Phi)_{\,22}\),  \(S^l(\Phi)_{\,22}\)
  are \(n_{-}\times n_{-}\)  matrices
(\(n_{+}\),\, \(n_{-}\) are defined in
 (\ref{cardequal+}) and (\ref{cardequal-})).
The block-decompositions (\ref{blocDecS^r}) and (\ref{blocDecS^l})
are consistent with the block-decompositions
 (\ref{PolMatrBlockDec})--(\ref{BlockResMatrR}) of the pole
and semiresidual matrices. The Sylvester-Lyapunov equation (\ref{SLRforPhi})
 for the matrix 
\(S^r(\Phi)\),  written in the
block-matricial form%
\[
\hspace*{-1.5ex}
\begin{array}{r}
\left[
\begin{array}{cc}
A_{\cal N}(\Phi)_1 & 0 \\
0 & A_{\cal N}(\Phi)_2
\end{array}
\right]
\left[
\begin{array}{cc}
S^r(\Phi)_{\,11}& S^r(\Phi)_{\,12}\\
S^r(\Phi)_{\,21}& S^r(\Phi)_{\,22}
\end{array}
\right]
-
\left[
\begin{array}{cc}
S^r(\Phi)_{\,11}& S^r(\Phi)_{\,12}\\
S^r(\Phi)_{\,21}& S^r(\Phi)_{\,22}
\end{array}
\right]
\left[
\begin{array}{cc}
A_{\cal P}(\Phi)_1 & 0 \\
0 & A_{\cal P}(\Phi)_2
\end{array}
\right]\\

\vspace*{-3ex}
\\
=
\left[
\begin{array}{cc}
G_{\cal N}(\Phi)_{\,1}\\
\vspace*{-3.5ex}
\\ 
 G_{\cal N}(\Phi)_{\,2}
\end{array}
\right]
\left[F_{\cal P}(\Phi)_1 \, \, \, F_{\cal P}(\Phi)_2\right],
\end{array}
\]
may be considered as a system of matricial equations
for the block-entries of the matrix \(S^r(\Phi)\).
This system is decomposed into four equations for block-entries of this
 matrix.
 
In particular, the equation for the entry \(S^r(\Phi)_{11}\) is of the form
\begin{equation}
A_{\cal N}(\Phi)_{\,1}S^r(\Phi)_{11}-S^r(\Phi)_{11}A_{\cal P}(\Phi)_{\,1}=
G_{\cal N}(\Phi)_1F_{\cal P}(\Phi)_1.
\label{forS^r_{11}}
\end{equation}
The Sylvester-Lyapunov equation (\ref{SLLforPhi})
 for the matrix \(S^l(\Phi)\), written in the
block-matricial form, may be considered as a system of matricial equations
for the block-entries of the matrix \(S^l(\Phi)\).
In particular, the equation for the entry \(S^l(\Phi)_{22}\) is of the form
\begin{equation}
A_{\cal P}(\Phi)_{\,2}S^l(\Phi)_{22}-S^l(\Phi)_{22}A_{\cal N}(\Phi)_{\,2}=
G_{\cal P}(\Phi)_2F_{\cal N}(\Phi)_2.
\label{forS^l_{22}}
\end{equation}
We show that
\begin{equation}
S^r(\Phi_{+})=S^r(\Phi)_{\,11}
\label{coins+}
\end{equation}

\vspace*{-3.5ex}
\noindent
and

\vspace*{-3.5ex}
\begin{equation}
S^l(\Phi_{-})=S^l(\Phi)_{\,22}.
\label{coins-}
\end{equation}
The easiest way to do this is to use the explicit formulas (\ref{S^})
for the matrices \(S^r\), \(S^l\) in terms of the poles, zeros and 
semiresidual matrices of the matrix function. For the function
\(\Phi \) the formula (\ref{S^}.r) takes the form
\begin{equation}
S^r(\Phi)=\|s^r_{p,q}(\Phi)\|
_{\mu_p\in{\cal N}(\Phi), \lambda_q\in{\cal P}(\Phi)}, \quad
s^r_{p,q}(\Phi)
=\frac{g_{\mu _p}(\Phi)\cdot f_{\lambda _q}(\Phi)}{\mu _p -\lambda _q}.
\label{explS^r}
\end{equation}
In particular,
\begin{equation}
S^r(\Phi)_{11}=\|s^r_{p,q}(\Phi)\|
_{\mu_p\in{\cal N}(\Phi)\cap G_{+},\, \lambda_q\in{\cal P}(\Phi)\cap G_{+}}, 
\quad
s^r_{p,q}(\Phi)
=\frac{g_{\mu _p}(\Phi)\cdot f_{\lambda _q}(\Phi)}{\mu _p -\lambda _q}.
\label{explS^r_11}
\end{equation}
 For the function \(\Phi_{+}\) the formula (\ref{S^}.r) takes the form
\begin{equation}
S^r(\Phi_{+})=\|s^r_{p,q}(\Phi_{+})\|
_{\mu_p\in{\cal N}(\Phi_{+}), \lambda_q\in{\cal P}(\Phi_{+})}, \quad
s^r_{p,q}(\Phi_{+})
=\frac{g_{\mu _p}(\Phi_{+})\cdot f_{\lambda _q}(\Phi_{+})}{\mu _p -\lambda _q}
\label{explS^r+}
\end{equation}
Comparing two last formulas and taking into account (\ref{locSing+}),
(\ref{ResRel+}) and  (\ref{ResInvRel+}), we conclude that
(\ref{coins+}) holds. In the same way, comparing the formulas
\begin{equation}
S^l(\Phi)=\|s^l_{p,q}(\Phi)\|
_{\lambda_p\in{\cal P}(\Phi), \mu_q\in{\cal N}(\Phi)}, \quad
s^l_{p,q}(\Phi)
=\frac{g_{\lambda_p}(\Phi)\cdot f_{\mu _q}(\Phi)}{\lambda _p -\mu _q}
\label{explS^l}
\end{equation}
and
\begin{equation}
S^l(\Phi_{-})=\|s^l_{p,q}(\Phi_{-})\|
_{\lambda_p\in{\cal P}(\Phi_{-}), \mu_q\in{\cal N}(\Phi_{-})}, \quad
s^l_{p,q}(\Phi_{-})
=\frac{g_{\lambda_p}(\Phi_{-})\cdot f_{\mu _q}(\Phi_{-})}{\lambda_p -\mu _q}
\label{explS^l-}
\end{equation}
and taking into account (\ref{locSing-}),
(\ref{ResRel-}) and  (\ref{ResInvRel-}), we conclude that
(\ref{coins-}) holds.

The explicit expressions (\ref{explS^r+}), (\ref{explS^l-})
 for the zero-pole coupling matrices 
are consequence of the Sylvester-Lyapunov equalities
 (\ref{LSfact+}), (\ref{LSfact-}). It is also possible to derive the equations
(\ref{coins+}) and (\ref{coins-}) directly from the equalities
(\ref{LSfact+}), (\ref{LSfact-}), bypassing the explicit
expressions (\ref{explS^r+}), (\ref{explS^l-}).
 (\ref{explS^r}), (\ref{explS^l}).
 The latter way is better, because it is applicable
{\sf not only to matrix functions in general position.}

Taking into account the equalities (\ref{RightPoleMatr+}),
and (\ref{PlusSemiRes}), we came from (\ref{LSfact+}) 
 to the equation
\begin{equation}
A_{\cal N}(\Phi)_{1}S^r(\Phi_{+})-S^r(\Phi_{+})A_{\cal P}(\Phi)_{1}=
G_{\cal N}(\Phi)_{1}F_{\cal P}(\Phi)_{1}.
\label{forS+}
\end{equation}
Comparing (\ref{forS^r_{11}}) and (\ref{forS+}) and taking into account the
uniqueness of the solution of the Sylvester-Lyapunov equation, we conclude
 that the equality (\ref{coins+}) holds.

In the same way, we can establish the equality (\ref{coins-}).
Taking into account the equalities (\ref{RightPoleMatr-}),
 (\ref{MinusSemiRes}), we came from (\ref{LSfact-}) to the
equation
\begin{equation}
A_{\cal P}(\Phi)_{2}S^l(\Phi_{-})-S^l(\Phi_{-})A_{\cal N}(\Phi)_{2}=
G_{\cal P}(\Phi)_2F_{\cal N}(\Phi)_2.
\label{forS-}
\end{equation}
Comparing (\ref{forS^l_{22}}) and (\ref{forS-}) and taking into account the
uniqueness of the solution of the Sylvester-Lyapunov equation, we conclude
 that the equality (\ref{coins-}) holds.

According to Theorem \ref{summarTheo} (applied to the rational matrix functions
in general position \(\Phi_{+}\) and \(\Phi_{-}\)),
 the matrices \(S^r(\Phi_{+})\) and \(S^l(\Phi_{-})\)
are invertible. In view of (\ref{coins+}) and (\ref{coins-}),
{\sf  the block-entries 
\boldmath\((S^r(\Phi))_{11}\)\unboldmath \,
 and \,\boldmath\((S^l(\Phi))_{22}\)}\unboldmath
\  (\,of the matrices \(S^r(\Phi)\) and
\(S^l(\Phi)\) respectively\,) {\sf are invertible.}

Now we may rewrite the formulas (\ref{forPhiR}), (\ref{forPhiL})
 for the factors
\(\Phi_{+}\) and \(\Phi_{-}\) in terms of block-entries of the
pole-, zero-, semiresidual- and zero-pole coupling matrices for the 
factorized matrix function \(\Phi\). Substituting the expressions
(\ref{RightPoleMatr+}.p), (\ref{PlusSemiRes}) and (\ref{coins+}) into
 (\ref{forPhiR}),
we obtain
\begin{equation}
\Phi_{+}(z)=I-{F_{\cal P}(\Phi)}_1(zI-{A_{\cal P}(\Phi)}_1)^{-1} 
({S^r(\Phi)}_{\,11})^{-1}{G_{\cal N}(\Phi)}_1. 
\label{forPhi+}
\end{equation} 
Analogously, substituting the expressions
(\ref{RightPoleMatr+}.n), (\ref{MinusSemiRes}) and (\ref{coins-}) into
 (\ref{forPhiL}),
we obtain
\begin{equation}
\Phi_{-}(z)=I+{F_{\cal N}(\Phi)}_2({S^l(\Phi)}_{\,22})^{-1}
(zI-{A_{\cal P}(\Phi)}_2)^{-1} 
{G_{\cal P}(\Phi)}_2. 
\label{forPhi-}
\end{equation} 
The formula (\ref{forPhi+}) expresses the factor \(\Phi_{+}\) in terms of the
 values
 \(F_{\cal P}(\Phi)\), \(G_{\cal N}(\Phi)\) and \(S^r(\Phi)\),
 which appear in the representation
(\ref{RSystForPhi}), whereas the formula (\ref{forPhi-}) expresses the factor
 \(\Phi_{-}\)
 in terms of
 the values \(F_{\cal N}(\Phi)\),
\(G_{\cal P}(\Phi)\) and \(S^l(\Phi)\), which appear in the {\sf other}
 representation
(\ref{LSystForPhi}). This disagreement is inconvenient for some calculations.
Therefore we also give the formula which expresses the factor \(\Phi_{-}\) in
 terms of
the the same values \(F_{\cal P}(\Phi)\), \(G_{\cal N}(\Phi)\) and
 \(S^r(\Phi)\),
 which appear in the expression (\ref{forPhi+})
 for the factor  \(\Phi_{+}\). To do it, we have to use the zero-pole
coupling relations (\ref{z-prelforPhi}.a), (\ref{z-prelforPhi}.c),
Together with (\ref{MinusSemiRes}.p) and (\ref{PlusSemiRes}.n), these
relations mean:
\begin{equation}
\begin{array}{cc}
F_{\cal N}(\Phi_{-})=\left(F_{\cal P}(\Phi)S^r(\Phi\right)^{-1})_2;&
G_{\cal P}(\Phi_{-})=-((S^r(\Phi))^{-1}G_{\cal N}(\Phi))_2.
\end{array}
\label{MinussemiresAltern}
\end{equation}
 \(\Phi\)):
The equality (\ref{invsForPhi}) together with (\ref{coins-}) means:
\begin{equation}
S^l(\Phi_{-})=\left(S^r(\Phi)^{-1}\right)_{22}.
\label{coins-altern}
\end{equation}

Substituting the expressions (\ref{MinussemiresAltern}) and
 (\ref{coins-altern})  into
(\ref{forPhiL}), we obtain
\begin{equation}
\Phi_{-}(z)=I-\Big({F_{\cal P}(\Phi)}\,S^r(\Phi)^{-1}\Big)_{2}
\Big(\Big({S^r(\Phi)^{-1}}\Big)_{\,22}\Big)^{-1}
\Big(zI-{A_{\cal P}(\Phi)}_2\Big)^{-1} 
\Big(S^r(\Phi)^{-1}\,G_{\cal N}(\Phi)\Big)_2. 
\label{uniforPhi-}
\end{equation}
The last formula already expresses the factor \(\Phi_{-}\) in terms of 
the values 
 \(F_{\cal P}(\Phi)\), \(G_{\cal N}(\Phi)\) and \(S^r(\Phi)\),
i.e. in terms of the same values which appear in the representation
 (\ref{forPhi+}) of the factor \(\Phi_{+}\). 

However, for the further considerations it will be useful to transform this
 formula,
substituting into it the expression for the inverse matrix \(S^r(\Phi)^{-1}\)
in terms of the block-entries of the matrix \(S^r(\Phi)\) itself.
First of all we recall a formula for the inversion of a \(2\times 2\) block
 matrix
with square diagonal block-entries. Let
\begin{equation}
M=\left[
\begin{array}{cc}
m_{11}&m_{12}\\
m_{21}&m_{22}
\end{array}
\right]
\label{blockM}
\end{equation}
be a square \(n\times  n\)
block-matrix matrix, which  block-entries \(m_{11}\) and \(m_{22}\) 
be square \(n_1\times  n_1\) and \(n_2\times  n_2\) matrices respectively
\((n=n_1+n_2)\).
 We assume that the matrix \(m_{11}\) is invertible. To inverse
the matrix \(M\), we first of all factorize it:
\begin{equation}
\left[
\begin{array}{cc}
m_{11}&m_{12}\\
m_{21}&m_{22}
\end{array}
\right]
=
\left[
\begin{array}{cc}
I&0\\
m_{21}\,m_{11}^{-1}&I
\end{array}
\right]
\cdot
\left[
\begin{array}{cc}
m_{11}      &        0\\
0           &        m_{22}-m_{21}\,m_{11}^{-1}\,m_{12}
\end{array}
\right]
\cdot
\left[
\begin{array}{cc}
I       &         m_{11}^{-1}m_{12}\\
0       &         I
\end{array}
\right].
\label{fact}
\end{equation}
The matrices \
 \(\left[
\begin{array}{cc}
I&0\\
m_{21}\,m_{11}^{-1}&I
\end{array}
\right]\) \
and \
\(
 \left[
\begin{array}{cc}
I       &         m_{11}^{-1}m_{12}\\
0       &         I
\end{array}
\right]
\) \ are invertible, and
\[
\left[
\begin{array}{cc}
I&0\\
m_{21}\,m_{11}^{-1}&I
\end{array}
\right]^{-1}=\left[
\begin{array}{cc}
I&0\\
-m_{21}\,m_{11}^{-1}&I
\end{array}
\right] ;\quad\quad
\left[
\begin{array}{cc}
I       &         m_{11}^{-1}m_{12}\\
0       &         I
\end{array}
\right]^{-1}=
\left[
\begin{array}{cc}
I       &        - m_{11}^{-1}m_{12^{}}\\
0       &         I
\end{array}
\right]
\]
 Therefore the matrix \(M\) is invertible if and only if
the matrix \ \(\left[
\begin{array}{cc}
m_{11}      &        0\\
0           &        m_{22}-m_{21}\,m_{11}^{-1}\,m_{12}
\end{array}
\right]\) \ is invertible. The latter is invertible if and only if the matrix
 \
\((m_{22}-m_{21}\,m_{11}^{-1}\,m_{12})\) \ is invertible. If the matrix
\((m_{22}-m_{21}\,m_{11}^{-1}\,m_{12})\) is invertible
 \footnote{The invertibility of
the matrix \ \(m_{11}\)\ was assumed from the very beginning}, then
\[
\left[
\begin{array}{cc}
m_{11}&m_{12}\\
m_{21}&m_{22}
\end{array}
\right]^{-1}=
\left[
\begin{array}{cc}
I       &        - m_{11}^{-1}m_{12}\\
0       &         I
\end{array}
\right]\cdot
\left[
\begin{array}{cc}
m_{11}^{-1}      &        0\\
0           &       ( m_{22}-m_{21}\,m_{11}^{-1}\,m_{12})^{-1}
\end{array}
\right]
\cdot
\left[
\begin{array}{cc}
I&0\\
-m_{21}\,m_{11}^{-1}&I
\end{array}
\right]\cdot
\]
Finally, for the matrix \(M\) with the decomposition (\ref{blockM}),
\begin{equation}
M^{-1}
=\quad 
\left[
\begin{array}{cc}
m_{11}^{-1}   &       0\\
0             &       0
\end{array}
\right]\
+\
\left[
\begin{array}{c}
- m_{11}^{-1}m_{12}   \\
I
\end{array}
\right]
\cdot
\Big[m_{22}-m_{21}\,m_{11}^{-1}\,m_{12}\Big]^{-1}
\cdot
\left[
\begin{array}{cc}
 -m_{21}\,m_{11}^{-1} \  &\  I
\end{array}
\right]\, \cdot
\label{blockinv}
\end{equation}
In particular,
\begin{equation}
(M^{-1})_{\,22}=(m_{22}-m_{21}\,m_{11}^{-1}\,m_{12})^{-1},
\label{m22}
\end{equation}
If \(F\) is \(k\times n\) matrix, with the block-decomposition
 \(F=[F_1\ \ F_2]\),
where \(F_1\), \(F_2\)  are  \(k\times n_1\) and \(k\times n_2\) matrices
 respectively,
then \(2\)-entry of the matrix \(F\,M^{-1}\) is of the form
\begin{equation}
(F\,M^{-1})_{\,2}=
\left[ F_1\ \ \ F_2\right]
\cdot
\left[
\begin{array}{c}
 -m_{11}^{-1}m_{12}\\
I
\end{array}
\right]
\cdot
\left(m_{22}-m_{21}\,m_{11}^{-1}\,m_{12}\right)^{-1}.
\label{2entryF}
\end{equation}
If \(G\) is \(n\times k\) matrix, with the block-decomposition
 \(
G=
\left[
\begin{array}{c}
G_1\\
G_2
\end{array}
\right]
\),
where \(G_1\), \(G_2\)  are  \(n_1\times k\) and \(n_2\times k\) matrices
 respectively,
then \(2\)-entry of the matrix \(M^{-1}\,G\) is of the form
\begin{equation}
\left(M^{-1}G\right)_{\,2}=
\left(m_{22}-m_{21}\,m_{11}^{-1}\,m_{12}\right)^{-1}\cdot
\left[
-m_{21}m_{11}^{-1}\ \ \ I
\right]
\cdot
\left[
\begin{array}{c}
G_1\\
G_2
\end{array}
\right] \cdot
\label{2entryf}
\end{equation}
Let us take the matrix \(S^r(\Phi)\) with the decomposition
 (\ref{blocDecS^r})  as the
matrix \(M\) \((n_1=n_{+}, n_2=n_{-})\) as the matrix \(M\). Both the
 matrix \(S^r(\Phi)\)
and its block-entry the matrix \(S^r(\Phi)_{11}\) are invertible. Therefore
the entry \((S^r(\Phi)^{-1})_{22}\)
of the inverse matrix \((S^r(\Phi))^{-1}=S^l(\Phi)\) is invertible as well,
 and,
according to (\ref{m22}),
\begin{equation}
((S^r(\Phi )^{-1})_{22})^{-1}=
S^r(\Phi)_{22}- S^r(\Phi)_{21}\,S^r(\Phi)_{11}^{-1}\,S^r(\Phi)_{12}.
\label{S^l_{22}}
\end{equation}
In particular, {\sf the matrix \    \boldmath
\( S^r(\Phi)_{22} - S^r(\Phi)_{21}\,S^r(\Phi)_{21}^{-1}\,S^r(\Phi)_{12}\)
\unboldmath  \  
is invertible}.\\
In view of {\rm (\ref{invsForPhi})},
the equality (\ref{coins-}) can be presented in the form
\begin{equation}
S^l(\Phi_{-})_{22}=
\left(S^r(\Phi)_{22} -
 S^r(\Phi)_{21}\,S^r(\Phi)_{21}^{-1}\,S^r(\Phi)_{12}\right)^{-1}.
\label{forS^l_Altern}
\end{equation}

Taking the matrix \(F_{\cal P}(\Phi)\)
with the decomposition (\ref{BlockResMatrL}.p) as the matrix \(F\) and
 the matrix \(G_{\cal N}(\Phi)\) with the decomposition
 (\ref{BlockResMatrR}.n)
as the matrix \(G\). we reduce the formula (\ref{uniforPhi-}) to the form
\begin{equation}
\begin{array}{c}
\Phi_{-}(z)=I- 
[F_{\cal N}(\Phi)_1 \ \ F_{\cal N}(\Phi)_2]\cdot
\left[
\begin{array}{c}
\ -S^r(\Phi)_{11}^{-1}\, S^r(\Phi)_{12}\\
I
\end{array}
\right]
\cdot
 \Big(zI-A_{\cal P}(\Phi)_2\Big)^{-1}\cdot \\
\cdot    
 \big(s^r(\Phi)_{22} -
 S^r(\Phi)_{21}\,S^r(\Phi)_{11}^{-1}\,S^r(\Phi)_{12}\Big)^{-1}
\cdot
\left[-S^r(\Phi)_{21}\,S^r(\Phi)_{11}^{-1} \ \ \  I \right]
\left[
\begin{array}{c}
 G_{\cal N}(\Phi)_1\\
G_{\cal N}(\Phi)_2
\end{array}
\right]\, \cdot
\label{finalForPhi-}
\end{array}
\end{equation}

Thus, we proved  the following

\vspace*{2ex}
\begin{theo}
\label{WNFdirTheo}
{\sl  Let \(\Phi \) be a rational \(k \times k\) matrix function in general
position, satisfying the normalizing condition {\rm(\ref{normPhi})}:
 \(\Phi(\infty)=I\),
with pole and zero sets \({\cal P}(\Phi)\) and \({\cal N}(\Phi)\),
pole and zero matrices \(A_{\cal P}(\Phi)\) and \(A_{\cal N}(\Phi)\),
semiresidual matrices 
\(F_{\cal P}(\Phi)\), \(F_{\cal N}(\Phi)\),
\(G_{\cal P}(\Phi)\), \(G_{\cal N}(\Phi)\),
and zero-pole coupling matrices \(S^r(\Phi)\) and \(S^l(\Phi)\), which are
 decomposed into blocks
as described above. \\[0.5em]
Assume that the function \(\Phi \) admits the spectral factorization
 with respect to \(\Gamma\), i.e.
\renewcommand{\theequation}{%
\mbox{\ref{globFactoriz}}}%
\begin{equation}
\Phi(z)={\Phi}_{+}(z)\cdot {\Phi}_{-}(z)\quad (z \in \overline{\Bbb C}),
\end{equation}
\renewcommand{\theequation}{%
\mbox{\arabic{section}.\arabic{equation}}}%
\addtocounter{equation}{-1}%
where \(\Phi_{+}\), \(\Phi_{-}\) are rational matrix functions,
with the pole- and zero-location: 
\renewcommand{\theequation}{%
\mbox{\ref{LocSingFac}}}%
\begin{equation}
 {\cal P}(\Phi_{+})\subset G_{+}, \quad
 {\cal N}(\Phi_{+})\subset G_{+},   \quad
{\cal P}(\Phi_{-})\subset G_{-},   \quad
{\cal N}(\Phi_{-})\subset G_{-}.
\end{equation}
\renewcommand{\theequation}{%
\mbox{\arabic{section}.\arabic{equation}}}%
\addtocounter{equation}{-1}%
Let  the normalizing condition
 {\rm (\ref{normPhi+})}: \(\Phi_{+}(\infty)=I\)
hold.\\[1.5ex]
Then:}
\begin{enumerate}
\item
{\sl For the matrix function \(\Phi\),
the number of poles and the number of ``zeros'',
located in \(G_{+}\), are equal:
\(\#({\cal P}(\Phi)\cap G_{+})=\#({\cal N}(\Phi)\cap G_{+})\);  
the number of poles and the number of ``zeros'',
 located in \(G_{-}\), are also equal:
\(\#({\cal P}(\Phi)\cap G_{-})=\#({\cal N}(\Phi)\cap G_{-})\).
}\\[-3.0ex]
\item
{\sl
The block-entries \(S^r(\Phi)_{11}\) and \(S^l(\Phi)_{22}\)\
{\rm (}in the above described block decompositions
{\rm (\ref{blocDecS^r}), (\ref{blocDecS^l})}
 of the  matrices
\(S^r(\Phi)\) and \(S^l(\Phi)\) respectively{\rm)},
as well as the matrix \
\mbox{\(S^r(\Phi)_{22}-S^r(\Phi)_{21}S^r(\Phi)_{11}^{-1}S^r(\Phi)_{12}\) }\
 are square invertible matrices.
}\\[-3ex]
\item
{\sl
The factors \(\Phi_{+}\) and \(\Phi_{-}\) are rational matrix functions in
general position, which are representable in the form
{\rm\((\ref{forPhi+})\)} and {\rm(\ref{forPhi-})}  respectively.
}\\[-3ex]
\item
{\sl
The representation {\rm\((\ref{forPhi+})\)} is the right system representation
of the factor \(\Phi_{+}\); the entries of this representation is expressible
in terms of the entries of the right system representation for \(\Phi\):
{\rm (\ref{locSing+})}, {\rm (\ref{RightPoleMatr+})},
{\rm (\ref{PlusSemiRes})} and {\rm(\ref{coins+})} hold.
The representation {\rm\((\ref{forPhi-})\)} is the left system representation
of the factor \(\Phi_{-}\); the entries of this representation is expressible
in terms of the entries of the right system representation for \(\Phi\):
{\rm (\ref{locSing-})}, {\rm (\ref{RightPoleMatr-})},
{\rm (\ref{MinusSemiRes})} and {\rm(\ref{coins-})} hold. The representation
{\rm\((\ref{forPhi-})\)} may be rewritten in terms of the entries
of the right system representation for \(\Phi_{+}\):
{\rm(\ref{finalForPhi-})} holds.
}
\end{enumerate}
\end{theo}

The converse statement is true as well.\\[1.5ex]
\begin{theo}
\label{invFactor}
{\sl  Let \(\Phi \) be a rational \(k \times k\) matrix function in general
position, satisfying the normalizing condition {\rm(\ref{normPhi})}:
 \(\Phi(\infty)=I\),
with pole and zero sets \({\cal P}(\Phi)\) and \({\cal N}(\Phi)\),
pole and zero matrices \(A_{\cal P}(\Phi)\) and \(A_{\cal N}(\Phi)\),
semiresidual matrices 
\(F_{\cal P}(\Phi)\), \(F_{\cal N}(\Phi)\),
\(G_{\cal P}(\Phi)\), \(G_{\cal N}(\Phi)\),
and zero-pole coupling matrices \(S^r(\Phi)\) and \(S^l(\Phi)\), 
which are decomposed into blocks as described above. 
 
\noindent
Assume that the following conditions are satisfied: 
}\\[-3.5ex]
\begin{enumerate}
\item
{\sl
None of the poles of \(\Phi\) as well as none of its ``zeros''
belongs to \(\Gamma\):
\[{\cal P}(\Phi)\cap\Gamma=\emptyset,\qquad
{\cal P}(\Phi^{-1})\cap\Gamma=\emptyset.\]
The number of poles and the number of zeros of the matrix-function \(\Phi\),
located in \(G_{+}\), are equal:
\[\#({\cal P}(\Phi)\cap G_{+})=\#({\cal N}(\Phi)\cap G_{+}),\]
or, what is the same,
the number of poles and the number of zeros of the matrix-function \(\Phi\),
located in \(G_{-}\), are equal:
\[\#({\cal P}(\Phi)\cap G_{-})=\#({\cal N}(\Phi)\cap G_{-}).\]
}\\[-3ex]

\item
{\sl
The block-entry \(S^r(\Phi)_{11}\) of the matrix \(S^r(\Phi)\)
{\rm (}in the above described block decompositions
{\rm (\ref{blocDecS^r})}{\rm)}  is  invertible matrix.
}
\end{enumerate}

\vspace*{-2.0ex}
\noindent
{\sl Then:}\\[-4.0ex]
\begin{enumerate}
\item
{\sl
The matrix functions \(\Phi_{+}\) and \(\Phi_{-}\), which are
{\sf  defined}  by
 {\rm (\ref{forPhi+})} and\,%
\footnote%
{
If the block-entry \(S^r(\Phi)_{11}\) is invertible,
the block-entry \(S^l(\Phi)_{22}\) is invertible as well.
}
\ {\rm (\ref{forPhi-})}, are
rational matrix functions in general position.
For these functions \(\Phi_{+}\) and \(\Phi_{-}\),
the conditions {\rm (\ref{LocSingFac})} are satisfied.
}\\[-3.5ex]
\item
{\sl
The matrix function \(\Phi\) admits the spectral factorization
{\rm (\ref{globFactoriz})} with respect to \(\Gamma\),
 with these matrix-functions 
\(\Phi_{+}\) and  \(\Phi_{-}\) as the factors.
}
\end{enumerate}
\end{theo}

\noindent
PROOF. To investigate the properties of the matrix functions \(\Phi_{+}\),
\(\Phi_{-}\), we will use Theorem \ref{chainspecial}.
To this theorem be applicable to the functions \(\Phi_{+}\), \(\Phi_{-}\),
we have to check whether its assumptions follow from the assumptions of
Theorem \ref{invFactor}. First of all, we have to be sure that none of
the columns of the matrices
\(F_{\cal P}(\Phi)_{1}\), \(F_{\cal N}(\Phi)_{2}\) is a zero-column
and none of the rows of the matrices
\(G_{\cal N}(\Phi)_{1}\), \(G_{\cal P}(\Phi)_{2}\)
is a zero-row. This property holds indeed because it holds for
the including matrices \(F_{\cal P}(\Phi)\), \(F_{\cal N}(\Phi)\),
\(G_{\cal P}(\Phi)\), \(G_{\cal N}(\Phi)\) (as for semiresidual
matrices for the rational matrix function \(\Phi\) in general
position). Then we have to check that the matrices 
\(S^r(\Phi)_{11}\), which appears at the ``core'' of the representation
(\ref{forPhi+}), satisfy some Sylvester-Lyapunov equality of the form
\(B\,S^r(\Phi)_{11}-S^r(\Phi)_{11\,}A=
G_{\cal N}(\Phi)_{1}\cdot F_{\cal P}(\Phi)_{1}\),
where \(A=A_{\cal P}(\Phi)_{1}\)
and \(B\) is a diagonal matrix. In the considered case, the identity
(\ref{forS^r_{11}}) plays the role of such Sylvester-Lyapunov equality, with
\(B=A_{\cal N}(\Phi)_{1}\). As we are already established, the equality
(\ref{forS^r_{11}}) is the consequence of the equality (\ref{SLRforPhi}).
The latter holds as the Sylvester-Lyapunov identity for the matrix function
\(\Phi\). Thus, the representation
 (\ref{forPhi+}) of the function
\(\Phi_{+}\) is its {\sf system} representation. According to 
Theorem \ref{chainspecial}, the matrix-functions \(\Phi_{+}\) is
a  rational matrix-function in general position, which
pole- and zero-sets are totalities of the diagonal entries
of the matrices \(A_{\cal P}(\Phi)_{1}\) and \(A_{\cal N}(\Phi)_{1}\)
 respectively:
\(
{\cal P}(\Phi_{+})=\Big\{\mbox{\rm diag}\,A_{\cal P}(\Phi)_{1}\Big\},
\) 
and
\(
{\cal N}(\Phi_{+})=\Big\{\mbox{\rm diag}\,A_{\cal N}(\Phi)_{1}\Big\}.\)
Hence (see the definition  (\ref{PolMatrBlockDec}), (\ref{ZerMatrBlockDec})
 of the matrices
\(A({\cal P})_{1}\), \(A({\cal N})_{1}\)\,), the conditions
\({\cal P}(\Phi_{+})={\cal P}(\Phi)\cap G_{+}\)
 \({\cal N}(\Phi_{+})={\cal N}(\Phi)\cap G_{+}\)
hold.

The matrix \(S^r(\Phi)\) is invertible as the Sylvester-Lyapunov matrix
related to the matrix-function \(\Phi\) (Theorem \ref{invtheo}.
 See  (\ref{invsForPhi})).
The corner block-entry \(S^r(\Phi)_{11}\) is invertible by the
assumptions of Theorem. Hence (see the identity 
{\rm (\ref{blockinv})} with \(S^r(\Phi)\) as \(M\)),
the matrix 
\( S^r(\Phi)_{22} - S^r(\Phi)_{21}\,S^r(\Phi)_{21}^{-1}\,S^r(\Phi)_{12}\)
is invertible as well, and  (\ref{S^l_{22}}) holds.
In view of  (\ref{invsForPhi}),  (\ref{forS^l_Altern}) holds.
In particular, the block-entry \(S^l(\Phi)_{22}\) is invertible as well.
Thus, the matrix-function \(\Phi_{-}\) is well defined by the formula
{\rm (\ref{forPhi-})}. The Sylvester-Lyapunov equality (\ref{forS^l_{22}})
for the matrix \(S^l(\Phi)_{22}\) follows from the equality
(\ref{SLLforPhi}). (The latter holds as the 'left''
 Sylvester-Lyapunov equality for \(\Phi\)). 
According to Theorem \ref{chainspecial}
(to be more precise, according to the ``left'' version of this Theorem),
\(\Phi_{-}\) is a rational matrix function in general position,
 which
pole- and zero-sets are totalities of the diagonal entries
of the matrices \(A_{\cal P}(\Phi)_{2}\)
 and \(A_{\cal N}(\Phi)_{2}\) respectively:
\(
{\cal P}(\Phi_{-})=\Big\{\mbox{\rm diag}\,A_{\cal P}(\Phi)_{2}\Big\},
\) 
and
\(
{\cal N}(\Phi_{-})=\Big\{\mbox{\rm diag}\,A_{\cal N}(\Phi)_{2}\Big\}.\)
Hence, the conditions
\({\cal P}(\Phi_{-})={\cal P}(\Phi)\cap G_{-}\)
 \({\cal N}(\Phi_{-})={\cal N}(\Phi)\cap G_{-}\)
hold.  All the more, the conditions (\ref{LocSingFac}) hold.

Thus, if the equality  (\ref{globFactoriz}) holds
 for \(\Phi_{+}\), \(\Phi_{-}\)
defined by  (\ref{forPhi+}),  (\ref{forS^l_Altern}), 
it gives the spectral factorization of \(\Phi\).

One remains only to verify the equality  (\ref{globFactoriz}).
The representation (\ref{forPhi-}) is convenient to investigate
the properties of \(\Phi_{-}\). However, to   verify   the equality
(\ref{globFactoriz}), it is more convenient to use
the representation (\ref{finalForPhi-}). To derive (\ref{finalForPhi-})
from  (\ref{forPhi-}), we have to use the zero-pole coupling relations
(\ref{z-prelforPhi}),
the formula (\ref{blockinv}) for inversion of \(2\times 2\) block- matrix
(applied to the matrix \(S^r(\Phi)\) with the block-decomposition
(\ref{blocDecS^r})) and, in particular, the equality (\ref{forS^l_Altern}).
 
To abbreviate the notation, we omit some notation entries, like indices etc.
So, for example, we write \(S\) instead \(S^r(\Phi\),
\(S_{11}\) instead \(S^r(\Phi)_{11}\),
\(F= [F_1 \ F_2]\) instead (\ref{BlockResMatrL}.p),
\(G=
\left[
\begin{array}{c}
G_1\\
G_2
\end{array}
\right]
{\rm instaed} \ 
{\rm (\ref{BlockResMatrL}.n)},
A=\left[
\begin{array}{cc}
A_1 & 0 \\
0   & A_2
\end{array}
\right]
\ {\rm instead} \
{\rm (\ref{PolMatrBlockDec})}, \
B=\left[
\begin{array}{cc}
B_1 & 0 \\
0   & B_2
\end{array}
\right]
\ {\rm instead} \
{\rm (\ref{ZerMatrBlockDec})}.
\)

Using the \(2\times 2\)-block-matrix inversion rule (\ref{blockinv}),
we present the representation (\ref{RSystForPhi}) for \(\Phi\) in the form
\begin{equation}
\Phi (z)=I- F\ (xI-A)^{-1}
\left(
\left[
\begin{array}{cc}
S_{11}^{-1}  & 0\\
0            & 0
\end{array}
\right]\ 
-
\left[
\begin{array}{c}
-\,S_{11}^{-1}\,S_{12}^{\ }\\
I
\end{array}
\right]
\Delta^{-1}
\left[
\begin{array}{cc}
-S_{21}S_{11}^{-1}& I
\end{array}
\right] \
\right)
G,
\label{PhiReprDetail}
\end{equation}
where
\[ \Delta =S_{22}-S_{21}S_{11}^{-1}S_{12}.\]
The representation (\ref{forPhi+}) for the function \(\Phi_{+}\)
we present in the form
\begin{equation}
\Phi_{+}(z)=I-F(zI-A)^{-1}
\left[
\begin{array}{cc}
S_{11}^{-1} & 0\\
0           & 0
\end{array}
\right]
G.
\label{forPhi_{+}verif}
\end{equation}
The representation (\ref{forPhi-}) for \(\Phi_{-}\),
rewritten in the form (\ref{finalForPhi-}), is:
\begin{equation}
\Phi_{-}(z)=I-F
\left[
\begin{array}{c}
-S_{11}^{-1}S_{12}\\
I
\end{array}
\right]
(zI-A_{2})^{-1}\Delta^{-1}
\left[
\begin{array}{cc}
-S_{21}S_{11}^{-1}  &   I 
\end{array}
\right]
G.
\label{forPhi_{-}verif}
\end{equation}
Multiplying the expressions in the right hand sides of 
(\ref{forPhi_{+}verif}) and (\ref{forPhi_{-}verif}) term by term, we obtain:
\begin{equation}
\hspace*{-2.0ex}
\begin{array}{lr}
\Phi_{+}(z)\Phi_{-}(z)=I- & \\[0.3ex]
-F(zI-A)^{-1}
\left[
\begin{array}{cc}
S_{11}^{-1} & 0\\
0           & 0
\end{array}
\right]
G -         &    {\rm (II)}   \\[2.3ex]
- F
\left[
\begin{array}{c}
-S_{11}^{-1}S_{12}\\
I
\end{array}
\right]
(zI-A_{2})^{-1}\Delta^{-1}
\left[
\begin{array}{cc}
-S_{21}S_{11}^{-1}  &   I 
\end{array}
\right]
G +           &   {\rm (III)}    \\[2.8ex]
+F(zI-A)^{-1}
\left[
\begin{array}{cc}
S_{11}^{-1} & 0\\
0           & 0
\end{array}
\right]
G
\cdot
F
\left[
\begin{array}{c}
-S_{11}^{-1}S_{12}\\
I
\end{array}
\right]
(zI-A_{2})^{-1}\Delta^{-1}
\left[
\begin{array}{cc}
-S_{21}S_{11}^{-1}  &   I 
\end{array}
\right].&   {\rm (IV)}
\end{array}
\label{termByTerm}
\end{equation}
Substituting the expression for \(GF\) from the
Sylvester-Lyapunov identity \(GF=BS-SA\), or
 \[
GF=S(zI-A)-(zI-B)S,
\]
into the expression (IV), we obtain:
\begin{equation}
\begin{array}{l}
 {\rm (IV)}=  \\
F(zI-A)^{-1}
\left[
\begin{array}{cc}
S_{11}^{-1} & 0\\
0           & 0
\end{array}
\right]
\cdot
S(zI-A)
\cdot
\left[
\begin{array}{c}
-S_{11}^{-1}S_{12}\\
I
\end{array}
\right]
(zI-A_{2})^{-1}\Delta^{-1}
\left[
\begin{array}{cc}
-S_{21}S_{11}^{-1}  &   I 
\end{array}
\right] -                          \\
-F(zI-A)^{-1}
\left[
\begin{array}{cc}
S_{11}^{-1} & 0\\
0           & 0
\end{array}
\right]
\cdot
(zI-B)S
\cdot
\left[
\begin{array}{c}
-S_{11}^{-1}S_{12}\\
I
\end{array}
\right]
(zI-A_{2})^{-1}\Delta^{-1}
\left[
\begin{array}{cc}
-S_{21}S_{11}^{-1}  &   I 
\end{array}
\right]
\end{array}
\label{forFour}
\end{equation}
It is clear that
\[
\left[
\begin{array}{cc}
S_{11}^{-1} & 0 \\
0 & 0
\end{array}
\right]
(zI-B)^{-1}=
\left[
\begin{array}{cc}
\ast  & 0 \\
0 & 0
\end{array}
\right],                 \quad
S
\left[
\begin{array}{c}
-S_{11}^{-1}S_{12}\\
I
\end{array}
\right]=
\left[
\begin{array}{c}
0 \\
\ast
\end{array}
\right].
\]
Thus,
\[
\left[
\begin{array}{cc}
S_{11}^{-1} & 0 \\
0 & 0
\end{array}
\right]
(zI-B)^{-1}
S
\left[
\begin{array}{c}
-S_{11}^{-1}S_{12}\\
I
\end{array}
\right]=
\left[
\begin{array}{c}
0 \\
0
\end{array}
\right],
\]
the second summand in the expression (\ref{forFour})
 for (IV) vanishes, and
\begin{equation}
 {\rm (IV)}=  \\
F(zI-A)^{-1}
\left[
\begin{array}{cc}
S_{11}^{-1} & 0\\
0           & 0
\end{array}
\right]
\cdot
S(zI-A)
\cdot
\left[
\begin{array}{c}
-S_{11}^{-1}S_{12}\\
I
\end{array}
\right]
(zI-A_{2})^{-1}\Delta^{-1}
\left[
\begin{array}{cc}
-S_{21}S_{11}^{-1}  &   I 
\end{array}
\right] 
G.
\end{equation}
As
\[
\left[
\begin{array}{cc}
S_{11}^{-1}&  0\\
0          &  0
\end{array}
\right]
\cdot S=
I-
\left[
\begin{array}{cc}
0          &  -S_{11}^{-1}S_{12}\\
0          &  I
\end{array}
\right]
,
\]
\[
\begin{array}{l}
{\rm (IV)}=
F
\left[
\begin{array}{c}
-S_{11}^{-1}S_{12}\\
I         
\end{array}
\right]
(zI-A_2)^{-1}{\Delta}^{-1}
\left[
\begin{array}{cc}
-S_{21}S_{11}^{-1}& I
\end{array}
\right]G
-                                       \\[2.8ex]
-
F(zI-A)^{-1}
\left[
\begin{array}{cc}
0      & -S_{11}^{-1}S_{12}\\
0      &  I
\end{array}
\right]
(zI-A)
\left[
\begin{array}{c}
-S_{11}^{-1}S_{12}\\
I         
\end{array}
\right]
(zI-A_2)^{-1}{\Delta}^{-1}
\left[
\begin{array}{cc}
-S_{21}S_{11}^{-1}& I
\end{array}
\right]G.
\end{array}
\]
It is clear that
\[
(zI-A)
\left[
\begin{array}{c}
-S_{11}^{-1}S_{12}\\
I         
\end{array}
\right]
(zI-A_2)^{-1}=
\left[
\begin{array}{c}
\ast\\
I         
\end{array}
\right],
\]
hence
\[
\left[
\begin{array}{cc}
0&-S_{11}^{-1}S_{12}\\
0&I         
\end{array}
\right]
(zI-A)
\left[
\begin{array}{c}
-S_{11}^{-1}S_{12}\\
I         
\end{array}
\right]
(zI-A_2)^{-1}=
\left[
\begin{array}{c}
-S_{11}^{-1}S_{12}\\
I         
\end{array}
\right].
\]
Thus,
\begin{equation}
\begin{array}{l}
{\rm (IV)}=
F
\left[
\begin{array}{c}
-S_{11}^{-1}S_{12}\\
I         
\end{array}
\right]
(zI-A_2)^{-1}{\Delta}^{-1}
\left[
\begin{array}{cc}
-S_{21}S_{11}^{-1}& I
\end{array}
\right]G
-                                       \\[2.8ex]
-
F(zI-A)^{-1}
\left[
\begin{array}{c}
-S_{11}^{-1}S_{12}\\
I         
\end{array}
\right]
{\Delta}^{-1}
\left[
\begin{array}{cc}
-S_{21}S_{11}^{-1}& I
\end{array}
\right]G.
\end{array}
\label{verifSpectFactId}
\end{equation}
>From (\ref{termByTerm}) and (\ref{verifSpectFactId}) it follows that
\[
\Phi_{+}(z)\,\Phi_{-}(z)=I-
F(zI-A)^{-1}
\left(
\left[
\begin{array}{cc}
S_{11}^{-1} & 0\\
0           & 0
\end{array}
\right]
+\left[
\begin{array}{c}
-S_{11}^{-1}S_{12}\\
I         
\end{array}
\right]
{\Delta}^{-1}
\left[
\begin{array}{cc}
-S_{21}S_{11}^{-1}& I
\end{array}
\right]
\right)
G
.
\]

Comparing the last expression to (\ref{PhiReprDetail}), we conclude that
the factorization equality (\ref{globFactoriz}) holds.
Theorem \ref{invFactor} is proved. \hfill\framebox[0.45em]{ }

\begin{rem}
\label{possibDegen}
Let \(\Phi\) be a rational matrix function in general position,
satisfying the normalizing condition
{\rm (\ref{normPhi})}: \(\Phi(\infty)=I\).
Assume that none of the poles of \(\Phi\) and  none of the ``zeros''
 of \(\Phi\) belongs to the
contour \(\Gamma\), and the conditions {\rm (\ref{cardequal+})}
:\ \ \(\#({\cal P}(\Phi)\cap G_{+})=\#({\cal N}(\Phi)\cap G_{+})\)
is satisfied;
(or, what is the same, the condition {\rm (\ref{cardequal-})}:\ \
\( \#({\cal P}(\Phi)\cap G_{-})=\#({\cal N}(\Phi)\cap G_{-})\)
is satisfied). This means that the matrices \(S^r(\Phi)_{11}\) and 
\(S^l(\Phi)_{22}\) are square one. According to Theorems
{\rm \ref{WNFdirTheo} and \ref{invFactor}},
 {\sf  the matrix function \(\Phi\) admits
the spectral factorization with respect to \(\Gamma\) if and only if
the matrix \(S^r(\Phi)_{11}\) is invertible, or, what is the same,
the matrix \(S^l(\Phi)_{22}\) is invertible.}
What can we say if the condition of the invertibility of the matrix
\(S^r(\Phi)_{11}\) is violated? According to the matrix factorization theory,
the matrix function \(\Phi\) admits the factorization of the form
\begin{equation}
\Phi(z)=\Phi_{+}(z)\,D(z)\,\Phi_{-}(z),
\label{genFact}
\end{equation}
where the matrix functions \(\Phi_{+}\) and \(\Phi_{+}^{-1}\) are holomorphic
on the set \(G_{-}\cup\Gamma\), 
the matrix functions \(\Phi_{-}\) and \(\Phi_{-}^{-1}\) are holomorphic
on the set \(G_{+}\cup\Gamma\) and \(D(z)\) is the matrix function of
the form
\begin{equation}
D(z)= {\rm diag}\,\Big( (z-z_{0})^{\kappa_1},\, (z-z_{0})^{\kappa_2}, \,
 \dots,\,
(z-z_{0})^{\kappa_k}\Big),
\label{partIndMatr}
\end{equation}
 \(\kappa_1,\,\kappa_2,\,\dots,\,\kappa_k\)  are integer numbers.
The numbers \(\kappa_1,\,\kappa_2,\,\dots,\,\kappa_k\) are said to be
{\sf the partial indices of the matrix function \(\Phi\) with respect
to the contour \(\Gamma\).}

The point \(z_{0}\) is an arbitrary chosen fixed point from \(G_{+}\).
Of course, the factors \(\Phi_{+}\),\, \(\Phi_{-}\) and \(D\) from
 the factorization
{\rm (\ref{genFact})} depend of the choice of the point \(z_{0}\).
However, {\sf the set} \(\{\kappa_1,\,\kappa_2,\,\dots,\,\kappa_k\}\)
(and under the normalizing condition
\begin{equation}
\kappa_1\leq \kappa_2 \leq\, \dots\, \kappa_k
\label{partIndNormCond}
\end{equation}
the matrix \(D\)) is determined uniquely by the function
\(\Phi\). In particular, the partial indices do not depend on the choice
of the distinguished  point \(z_{0}\). In contrast to this, even for the given
distinguished  point \(z_{0}\) the matrices \(\Phi_{+}\) and \(\Phi_{-}\)
are determined non-uniquely. However, this non-uniqueneess can be easily
described. (See Theorem 7.1 from \cite{GoKr}).

The partial indices play a fundamental role; in the homogeneous Hilbert
problem they were first introduced by N.I.Muskhelishvili and
N.P.Vekua \cite{MuVe}. (The factorization of another form:
\(\Phi(z)=\Phi_{+}(z)\,\Phi_{-}(z)\,D(z)\) was considered by G.Birkgoff
much earlier, in 1913 (see \cite{Bir1}). The matrices \(\Phi_{+}\),
\(\Phi_{-}\) and \(D\) have the same properties that in above described
factorization {\rm (\ref{genFact})} (in particular, the matrix \(D\) is of the
 form
{\rm (\ref{partIndMatr})}),  but the matrices themselves are different. In
particular, Birkgoff's partial indices may be different. The relation between
these kinds of factorizations is studed in \cite{FM}.

The natural question arise {\sf how to compute the factors
\(\Phi_{+}\), \(\Phi_{-}\), \(D\) in terms of the pole, zero and semiresidual
matrices for the given rational matrix function
\(\Phi\) in general position?
}
Of course, there are general methods for performing such a factorization
(see, for example, {\rm \cite{ClGo}}). However, we believe that for a rational
matrix function {\sf in general position} the factorization may be done
much more explicitly. 
\hfill\framebox[0.45em]{ }

\end{rem}
\begin{rem}
\label{bifurkFact}
Let as assume now that the rational matrix function \(\Phi\) in general
position depends holomorphically on some parameter \(\alpha\in {\Bbb C}\)
and satisfy the normalizing condition {\rm (\ref{normPhi})}:
\(\Phi(\infty, \alpha)\equiv I\).
The factors \(\Phi_{+}\), \(\Phi_{-}\), \(D\) in the factorization
{\rm (\ref{genFact})}:
\renewcommand{\theequation}{%
\mbox{\({\ref{genFact}}_{\,\,\alpha}\)}}%
\begin{equation}
\Phi(z,\alpha)=\Phi_{+}(z,\alpha)\,D(z,\alpha)\,\Phi_{-}(z,\alpha),
\end{equation}
\renewcommand{\theequation}{%
\mbox{\arabic{section}.\arabic{equation}}}%
\addtocounter{equation}{-1}%
\renewcommand{\theequation}{%
\mbox{\({\ref{partIndMatr}}_{\,\,\alpha}\)}}%
\begin{equation}
D(z,(\alpha))=
 {\rm diag}\,\Big( (z-z_{0})^{\kappa_{1}(\alpha)},\,
 (z-z_{0})^{\kappa_2(\alpha)}, \,
 \dots,\,
(z-z_{0})^{\kappa_k(\alpha)}\Big)
\end{equation}
\renewcommand{\theequation}{%
\mbox{\arabic{section}.\arabic{equation}}}%
\addtocounter{equation}{-1}%
depend now on \(\alpha\), and under the normalizing conditions
 {\rm (\ref{normPhi+})}: \(\Phi_{+}(\infty, \alpha)\equiv I\)
and {\rm (\ref{partIndNormCond})}:
\(\kappa_1(\alpha)\leq \kappa_2(\alpha) \leq\, \dots\, \kappa_k(\alpha)\)
are determined uniquely.

Assume now  that for all \(\alpha\) from some neighborhood of some point 
\(\alpha_{\scriptscriptstyle 0}\in{\Bbb C}\) the condition
\(\#({\cal P}(\Phi(\alpha))\cap G_{+})=\#({\cal N}(\Phi(\alpha))\cap G_{+})\) 
is satisfied (thus the block-entry \(S^r(\Phi(\alpha))_{11}\) is a square
matrix for all \(\alpha\) which are close to
\(\alpha_{\scriptscriptstyle 0}\)). Assume also that
\({\rm det}\,(S^r(\Phi(\alpha))_{11})\not\equiv 0\), but
\({\rm det}\,(S^r(\Phi(\alpha_{\scriptscriptstyle 0}))_{11})= 0\).
According to Theorem {\rm \ref{invFactor}}, applied to the matrix-function
\(\Phi(\alpha)\) with \(\alpha\neq\alpha_{\scriptscriptstyle 0}\)) (and
\(\alpha\) which is close to \(\alpha_{\scriptscriptstyle 0}\)),
the matrix function \(\Phi(z,\alpha)\)\ 
(\(\alpha\neq\alpha_{\scriptscriptstyle 0}\))
 admits the factorization 
of the form \(\Phi(z,\alpha)=\Phi_{+}(z,\alpha)\,\Phi_{-}(z,\alpha)\),
in other words, \(D(z,\alpha)\equiv I\) for \(\alpha\neq\alpha_{0}\).
Moreover, the factors \(\Phi_{+}(z,\alpha)\), \(\Phi_{-}(z,\alpha)\)
depend on \(\alpha\) holomorphically for 
\(\alpha\neq\alpha_{\scriptscriptstyle 0}\) (this may be established from the
explicit formulas {\rm (\ref{forPhi+}), (\ref{forPhi-})}
 for the factors).
However, according to Theorem {\ref{WNFdirTheo}}, the factorization of
this form is impossible for \(\alpha=\alpha_{\scriptscriptstyle 0}\)
(otherwise the Lyapunov-Sylvester matrix 
\(S^r(\Phi(\alpha_{\scriptscriptstyle 0}))\)
would be invertible). Thus, \(D(\alpha_{\scriptscriptstyle 0}))\neq I\),
i.e. not all \(\kappa_{j}(\alpha_{\scriptscriptstyle 0}), j=1, 2, \dots, k\),
vanish (although still
 \(\sum_{1\leq j \leq k}\kappa_{j}(\alpha_{\scriptscriptstyle 0})=0\)) .
So, the factorization {\rm (\ref{genFact})} undergoes a bifurcation
at the value of the parameter \(\alpha\) such that
 \({\rm det}\,S^r(\Phi(\alpha)_{11})\) vanishes.
The natural question arise {\sf how to describe the bifurcation in a
clear way?}
We believe that for rational matrix function \(\Phi\) in general
position it may be done more or less explicitly.
\end{rem}

\vspace*{0.6cm}

\noindent
\vspace*{0.2cm} 
\begin{minipage}{15.0cm}
\section{\hspace{-0.4cm}.\hspace{0.19cm}
SOME HISTORICAL REMARKS
}
\end{minipage}\\[-0.5cm]
\setcounter{equation}{0}

There are many different ways to specify and represent analytic
functions: for instance, Taylor series, decomposition in continuous
fractions, representations by Cauchy integrals or by Fourier integrals,
etc. (Of course, the distinction between different representation
methods is often artificial and hard to make.) 
In the first half of 1970's this toolkit was enriched by an additional
representation method: so called system realizations of analytic functions.
The sources of the theory of system realizations belong to several
different domains, in particular synthesis theory of linear electrical networks,
the theory of linear control systems, and the theory of
operator colligations (or nodes) and their characteristic functions.
Investigations in these theories were carried out by representatives
of different scientific disciplines. 
The investigations done by the mathematicians
have their root in the pioneering work of M. S. Liv\v{s}ic, who
is the forefather of the theory of system realizations.

In the middle of 1940's M. S. Liv\v{s}ic has introduced the notion
of the characteristic function of a linear operator. 
This notion was first introduced for nonselfadjoint extensions
of isometric operators with defect indices $(1,1)$ \cite{L1}
and more generally $(n,n)$ \cite{L2,L3} (for $n>1$ the characteristic
function is matrix valued), and later for general operators with
finite nonhermitian (or nonunitary) rank \cite{L4}.
M. S. Liv\v{s}ic has discovered the following properties of
the characteristic function \footnote{
we formulate these properties, on purpose, in a somewhat rough
and therefore slightly imprecise form}.
\begin{enumerate}
\item
The characteristic function determines the corresponding operator
essentially uniquely up to unitary equivalence
(first results of this kind are contained already in \cite{L1}).
\item
For each invariant subspace of the operator there is a decomposition
of the characteristic function into a product of two factors:
one of these factors is the characteristic function of the restriction
of the original operator onto the given invariant subspace, and
the other factor is the characteristic function of the compression
of the original operator onto the corresponding coinvariant subspace.
\end{enumerate}
Property 2 is the so called ``multiplication theorem'' for
characteristic functions. Initially the multiplication theorem
was established under various additional technical restrictions.
The final formulation of the theorem has been obtained in the framework
of the theory of operator colligations and their characteristic functions.

The theory of operator colligations that was created by M. S. Liv\v{s}ic
(and that was further developped by his collaborator M. S. Brodski\i{i},
see the book \cite{Br}) was a natural development of the theory
of nonselfadjoint operators and their characteristic functions.
This development was also intimately related with applications
of the theory of commuting nonselfadjoint operators to physical problems,
in particular to the problems of scattering and to elementary particles
physics, and later to synthesis problems for electrical networks
\cite{L5,L6,LF}. The theory of open systems took a definitive shape
in the works of M. S. Liv\v{s}ic in the early 1960's; the contents
of these works were incorporated in the monograph \cite{L9}.

A theory parallel to the theory of open systems of M. S. Liv\v{s}ic
has been developped by several other authors under the name of
the theory of linear stationary dinamical systems. (For an exposition
of system theory see \cite{Z}, \cite{Fu}, \cite{KFA}; we especially
recommend the survey \cite{Kaas}.) The transfer matrix function of
such a system is precisely the characteristic function of the corresponding
operator colligation. As M. S. Liv\v{s}ic has shown in several important
examples, for scattering systems the characteristic function coincides
with the scattering matrix \cite{L5,L6}.

A different line of investigation leading to the theory of system
realizations is connected with control theory and with the name of R. Kalman
(see \cite{KFA}). Here one also encounters the notion of the transfer
matrix function. Let us emphasize in this connection one important 
circumstance. In physical problems that M. S. Liv\v{s}ic was motivated by
there appeared always an ``energy balance'' condition implying
$J$-contractiveness of the characteristic (transfer) function and
its ``symmetry'' with respect to the unit circle (or the real axis).
It also imposed considering the adjoint operator $A^*$ together with
the state operator $A$ of the system.
Energy balance condition does not play an important role in control theory
and this leads naturally to considering a general pair of operators
$A,B$ on the state space instead of the pair $A,A^*$.
(In fact, R. Kalman develops system theory over arbitrary fields rather
than over the field of complex numbers.)

Many results of the analytic theory of electrical networks can be
considered as realization results for matrix functions of various classes.
(There are many expositions of the theory of electrical networks;
we recommend especially the monograph of V. Belevitch \cite{Bel} that
seems as if written for a mathematician. See also the survey \cite{EfPo}.)
The well known Darlington's synthesis method for passive networks
has been formulated for mathematicians by V. P. Potapov in 1966 
as a realization problem for passive rational matrix functions \cite{Pot3};
this problem was considered in details by Potapov's Ph. D. student
E. Melamud \cite{Me}.

Already in \cite{L4} M. S. Liv\v{s}ic used the theory of characteristic
functions for the reduction of a nonselfadjoint operator to a triangular
form, generalizing the theorem of I. Schur that an arbitrary matrix
can be brought to an upper (or lower) triangular form by unitary equivalence.
In this approach a multiplicative decomposition of the characteristic
function corresponds to an ``additive'' decomposition of the operator
itself over a linearly ordered chain of its invariant subspaces.
This correspondence was used by M. S. Liv\v{s}ic in both directions.
Using invariant subspaces of a finite-dimensional approximation of
a given operator he decomposes the characteristic function of
the approximating operator into factors, then passing to a limit he obtains
a multiplicative decomposition of the characteristic function of
the given operator, and finally using this multiplicative decomposition
he constructs a triangular model of the given operator which is unitary
equivalent to the operator itself.

A multiplicative decomposition of a meromorphic $J$-contractive matrix
function on the unit disk (or on the upper half plane) has been obtained
by V. P. Potapov in \cite{Pot3} in a purely function theoretic way 
(the simpler case $J=I$ was handled much earlier in \cite{Pot1}).
However an important special case of V. P. Potapov's theorem has already been
obtained by M. S. Liv\v{s}ic in \cite{L4} using operator theoretic methods;
V. P. Potapov used the investigations of M. S. Liv\v{s}ic as a guideline
in his own function theoretic approach.

It is natural to ask whether one can reduce an operator to a diagonal
form by a linear --- no longer unitary --- transformation.
In the middle of 1950's M. S. Liv\v{s}ic has posed this problem to his then
Ph. D. student L. A. Sakhnovich who has obtained numerous
results in this direction. In particular the problem of the reduction
of a nonselfadjoint operator to a diagonal form has lead L. A. Sakhnovich
to a relation of the form
\begin{equation} \label{5.1}
AX-XA^* = GJG^*
\end{equation}
considered as an equation for $X$ (formula (3) in \cite{S1}).
This relation appeared in a hidden form also in the theory
of $J$-contractive matrix functions
(see formula (22) and the following unnumbered formula in Chapter 2,
Section 4 of \cite{EfPo}) and in the theory of classical interpolation
problems (formula (12) in Section 1 of \cite{Kov}).
A relation of the form (\ref{5.1}) in \cite{Kov} appears exactly
in connection with what we called the chain identity; however
the chain identity plays there a secondary role and is not highlighted.
It should be mentionned that all the main ideas of the paper \cite{Kov},
published in 1983, have been suggested by V. P. Potapov some 10--12 years
earlier; unfortunately, V. P. Potapov's contribution is not adequately
reflected there. An identity of the form (\ref{5.1}) has also been
considered by L. de Branges in the framework of a certain generalized
moment problem (in the language of Hilbert spaces of entire functions;
see Theorem 27 of \cite{Bran} and also formula (6.10) in \cite{GolM}).
A. A. Nudelman has used an identity of the form (\ref{5.1}) at the basis
of an abstract scheme that he developped for considering classical
interpolation problems.
It is interesting to note that an identity of the form (\ref{5.1})
(and a related identity
$X - AXA^* = GJG^*$)
appears in an entirely different context as well, namely
in connection with fast inversion algorithms for structured matrices
(Toeplitz, Hankel, Vandermonde, etc.); see the survey \cite{KS}.

Everything needed for the theory of system realizations in its present
form has thus been available by the middle of 1970's.
An important step was taken by L. A. Sakhnovich in \cite{S2}
(a detailed exposition of these results is contained in \cite{S3}):
he studied the spectral factorization of a rational matrix function
$R$ with both $R$ and $R^{-1}$ given as transfer functions of
the corresponding linear systems (operator colligations).
The spectra of the state space operators $A$ and $B$ of these systems do not
intersect. One considers a pair of Sylvester--Lyapunov equations
$AT-TB=F_1G_1$, $SA-BS=F_2G_2$,
where $F_1$, $G_1$, $F_2$, $G_2$ are the input and the output operators
of the systems realizing $R$, $R^{-1}$.
It is shown that if the corresponding blocks of the solutions $T$ and $S$
of these equations are invertible then the matrix function $R$ admits
a spectral factorization, and formulas for the spectral factors
(having the form of the formulas (4.53--4.54) of the present paper)
are obtained. N. M. Kostenko (a Ph. D. student of L. A. Sakhnovich)
has shown in \cite{Kos} that the invertibility of these blocks is
necessary for the existence of a spectral factorization.
Relations which are analogous to zero-pole coupling relations were
considered already in \cite{S1} and used in \cite{S2}.
Let us notice that the spectral factors constructed in \cite{S2}
are simply the chacteristic (transfer) matrix functions of subsystems
arising by restricting the respective state space operators to
the corresponding spectral subspaces. {\em This is the form taken by
the multiplication theorem of M. S. Liv\v{s}ic in the current situation}
(the matrices $T$ and $S$ define ``metrics'' which are now just
bilinear functionals, neither positive definite nor even hermitian).

Unfortunately the paper \cite{S2} did not have the impact it deserved.
The subsequent development of the theory of system realizations is connected
with the name of I. Gohberg. I. Gohberg has also lead and inspired
a coherent work of many mathematicians in the theory of system realizations
and its applications and this theory experienced a fast growth from
the late 1970's onward. Already in 1979 there appeared the monograph 
\cite{BGK1} dealing with spectral factorizations of rational matrix functions
given as transfer functions of linear systems (operator colligations).
The paper \cite{GKLR} considers the realization problem for matrix functions
$R$ and $R^{-1}$ as transfer functions starting with ``local data''
(principal parts of Laurent series for $R$ and $R^{-1}$ at each pole);
see also \cite{BGR1}. In the beginning it was assumed that the pole sets
for $R$ and $R^{-1}$ do not intersect; later the general case when these
sets may intersect was considered as well. These and many other questions
are considered in great details in the monograph \cite{BGR2};
see also \cite{KRR}. The collection of papers \cite{CoMe} is
dedicated to the spectral factorization for rational matrix functions
based on the theory of system realizations; the papers \cite{BGK2}
and \cite{BGK3} are especially related with our exposition in Section 4.
There are results on system realization for rectangular (non-square)
rational matrix functions given by local data \cite{BGRa}.
There are also realization results for matrix functions on
a Riemann surface \cite{BV},
using deep new ideas of M. S. Liv\v{s}ic and his collaborators 
on characteristic functions
for {\em commuting tuples} of nonselfadjoint operators \cite{LKMV}.

Factorization of matrix functions is a tool for many other problems,
e.g., the theory of inverse problems for differential equations and
prediction theory of stationary stochastic processes.
If the corresponding matrix function is rational,
this factorization (which is a technical tool for the original problem)
may be carried out using system realizations which then become involved
in the solution of the original problem as well.
See, e.g., \cite{AG1,AG2}. It is clear that the theory of system realizations
can be successfully used also for the solution of the problems considered
in \cite{Yag}. It would be interesting to connect the questions
considered in \cite{Dei} with the theory of system realizations.

\end{document}